\documentclass[11pt]{amsart}
\pdfoutput=1
\usepackage[colorlinks, linkcolor=blue, citecolor=blue, pagebackref]{hyperref}
\usepackage{subcaption}
\usepackage{tabularray}
\usepackage{makecell}
\usepackage{amssymb}
\usepackage{amsrefs}
\usepackage{mathtools}
\usepackage{ytableau}
\usepackage[indent]{parskip}
\usepackage{stmaryrd}
\usepackage{tikz}
\usetikzlibrary{decorations.pathmorphing, shadows.blur, fit, positioning}
\usepackage{enumitem}
\usepackage{fullpage}
\usepackage{halloweenmath}
\usepackage{multirow}
\newcounter{x}
\newcounter{y}
\newcounter{z}
\newcounter{h}
\newcommand\xaxis{210}
\newcommand\yaxis{-30}
\newcommand\zaxis{90}
\newcommand\topside[3]{
  \fill[fill=white, draw=black,shift={(\xaxis:#1)},shift={(\yaxis:#2)},
  shift={(\zaxis:#3)}] (0,0) -- (30:1) -- (0,1) --(150:1)--(0,0);
}
\newcommand\leftside[3]{
  \fill[fill=black!25, draw=black,shift={(\xaxis:#1)},shift={(\yaxis:#2)},
  shift={(\zaxis:#3)}] (0,0) -- (0,-1) -- (210:1) --(150:1)--(0,0);
}
\newcommand\rightside[3]{
  \fill[fill=black!50, draw=black,shift={(\xaxis:#1)},shift={(\yaxis:#2)},
  shift={(\zaxis:#3)}] (0,0) -- (30:1) -- (-30:1) --(0,-1)--(0,0);
}
\newcommand\cube[3]{
  \topside{#1}{#2}{#3} \leftside{#1}{#2}{#3} \rightside{#1}{#2}{#3}
}
\newcommand\floorside[3]{
  \fill[fill=black!15, draw=black,shift={(\xaxis:#1)},shift={(\yaxis:#2)},
  shift={(\zaxis:#3)}] (0,0) -- (30:1) -- (0,1) --(150:1)--(0,0);
}
\newcommand\planepartition[1]{
\setcounter{h}{0}
 \foreach \a in {#1}{
    \foreach \b in \a {
        \ifnum \b>\value{h} \setcounter{h}{\b}\fi
        }
 }
 \setcounter{x}{-1}
  \foreach \a in {#1} {
        \addtocounter{x}{1}
        \setcounter{y}{-1}
        \foreach \b in \a {
            \addtocounter{y}{1}
            \setcounter{z}{-1}
            \ifnum \b>0 {
                \pgfmathtruncatemacro\suivant{\b+1}
                \foreach \c in {1,...,\value{h}} {
                    \addtocounter{z}{1}
                    \ifnum \c<\suivant             
                        \cube{\value{x}}{\value{y}}{\value{z}}
                    \else {
                        \ifnum\value{x}=0 \leftside{-1}{\value{y}}{\value{z}}\fi
                        \ifnum\value{y}=0 \rightside{\value{x}}{-1}{\value{z}}\fi
                        }\fi
                    }
             }
            \else {\floorside{\value{x}}{\value{y}}{\value{z}}}
            \fi
    }
  }
}

\tikzstyle{dot}=[circle,fill=black, minimum size = 2.25pt, inner sep=0pt]

\theoremstyle{plain}
\newtheorem{theorem}{Theorem}[section]
\newtheorem{prop}[theorem]{Proposition}
\newtheorem{lemma}[theorem]{Lemma}

\newtheorem{conj}[theorem]{Conjecture}

\theoremstyle{definition}
\newtheorem{dfn}[theorem]{Definition}
\newtheorem{example}[theorem]{Example}
\newtheorem{problem}[theorem]{Problem}

\theoremstyle{remark}
\newtheorem{remark}[theorem]{Remark}
\numberwithin{equation}{section}

\newcommand{\Qsigma}{\mathcal{Q}_k(\sigma)}
\newcommand{\Tsigma}{\mathcal{T}(\sigma)}
\newcommand{\PP}{\mathcal{P}_k}
\newcommand{\C}{\mathbb C}
\newcommand{\R}{\mathbb R}
\DeclareMathOperator{\GL}{GL}
\DeclareMathOperator{\U}{U}
\DeclareMathOperator{\SO}{SO}
\newcommand{\gl}{\mathfrak{gl}}
\DeclareMathOperator{\Sp}{Sp}
\DeclareMathOperator{\Mp}{Mp}
\renewcommand{\O}{\operatorname{O}}
\renewcommand{\sp}{\mathfrak{sp}}
\newcommand{\so}{\mathfrak{so}}
\newcommand{\g}{\mathfrak{g}}
\renewcommand{\k}{\mathfrak{k}}
\newcommand{\p}{\mathfrak{p}}
\newcommand{\scrO}{\mathcal{O}}
\newcommand{\Deg}{\operatorname{Deg}}

\DeclareMathOperator{\SSYT}{SSYT}
\newcommand{\la}{\lambda}

\ytableausetup{smalltableaux}

\begin{document}

\title{A combinatorial interpretation of the Bernstein degree \\ of unitary highest weight modules}

\author{William Q.~Erickson}
\address{
Department of Mathematics\\
Baylor University \\ 
One Bear Place \#97328\\
Waco, TX 76798} 
\email{william.q.erickson@gmail.com}

\author{Markus Hunziker}
\address{
Department of Mathematics\\
Baylor University \\ 
One Bear Place \#97328\\
Waco, TX 76798} 
\email{Markus\_Hunziker@baylor.edu}

\begin{abstract}
    The Bernstein degree ($\Deg$) is a fundamental invariant of admissible representations of a real reductive Lie group $G_\R$.
    Our main result concerns the classical dual pairs $(G_\R, H_\R(k))$, namely $(\U(p,q), \: \U(k))$, $(\Mp(2n, \R), \: \O(k))$, and $(\O^*(2n), \: \Sp(k))$, where $k$ is any positive integer.
    In this setting, via Howe duality, each irreducible representation $\sigma$ of $H_\R(k)$ corresponds to a unitary highest weight module $L_{\la(\sigma)}$ for $G_\R$.
    A landmark result of Nishiyama–Ochiai–Taniguchi (2001) expressed $\Deg L_{\la(\sigma)}$ as a product of two quantities: 
    the dimension of $\sigma$ and the degree of the associated variety.
    However, this result was limited to a specific range of the parameter $k$ (namely $k \leq r$, the real rank of $G_\R$).
    The present paper resolves this limitation by introducing, for all~$k$, the combinatorial interpretation $\Deg L_{\la(\sigma)} = \# ( \Qsigma \times \PP)$, where $\Qsigma$ is a certain set of semistandard tableaux and $\PP$ is a set of plane partitions.
    (The result remains partly conjectural in the $\Mp(2n, \R)$ case.)
    Beyond the dual pair setting, we generalize the set $\PP$ to all groups $G_\R$ of Hermitian type, and we exhibit analogues of the Nishiyama–Ochiai–Taniguchi result for certain families of unitary highest weight modules of $\operatorname{E}_6$ and $\operatorname{E}_7$.
\end{abstract}

\subjclass[2020]{Primary 22E46; Secondary 05E10, 14L30, 32M15}

\keywords{Unitary highest weight modules, Bernstein degree, associated variety, Hermitian symmetric space, Howe duality, Stanley decompositions}

\maketitle

\setcounter{tocdepth}{1}
\tableofcontents

\section{Introduction}

\subsection{Bernstein degree of unitary highest weight modules}

Let $X$ be the $(\g, K)$-module of an admissible representation of a connected noncompact real reductive Lie group $G_\R$ with finite center.
Two invariants arising from the asymptotic growth of a good $U(\g)$-filtration on $X$ are its \emph{Gelfand--Kirillov dimension} ${\operatorname{Dim}}(X)$ and its \emph{Bernstein degree} $\Deg X$.
By a result of Vogan~\cite{Vogan}*{Thm.~1.1}, for any good filtration $\{ X_n \}$,
\[
\dim X_n \sim \frac{\Deg X}{{\operatorname{Dim}}(X)!} n^{{\operatorname{Dim}}(X)}
\]
as $n \rightarrow \infty$.
Equivalently, there is a unique polynomial $Q_X(t)$ with integer coefficients, such that the reduced Hilbert series of $X$ is given by
\begin{equation}
    \label{Deg = Q(1)}
    H_X(t) = \frac{Q_X(t)}{(1-t)^{{\operatorname{Dim}}(X)}}, \quad \text{and} \quad \Deg X = Q_X(1).
\end{equation}
Uniting these two invariants is the notion of the associated variety $\mathcal{AV}(X) \subset \g^*$, introduced by Vogan~\cites{Vogan,Vogan81}.
In particular, we have ${\operatorname{Dim}}(X) = \dim \mathcal{AV}(X)$, while $\Deg X$ encodes the multiplicity of $\mathcal{AV}(X)$ as an algebraic cycle.

Harish-Chandra~\cites{HC55,HC56} showed that a noncompact simple Lie group $G_\R$ admits nontrivial unitary highest weight modules if and only if $G_\R$ is of Hermitian type, that is, $G_\R / K_\R$ is a Hermitian symmetric space (where $K_\R$ is a maximal compact subgroup).
From now on, suppose that $G_\R / K_\R$ is an irreducible Hermitian symmetric space of noncompact type.
Let $\g_\R = \k_\R \oplus \p_\R$ be the corresponding Cartan decomposition of $\g_\R$, and let $\g = \p^- \oplus \k \oplus \p^+$
be the triangular decomposition of its complexification.
The action of the complexified group $K$ on the abelian nilradical $\p^+$ has finitely many orbits, whose closures form a chain of algebraic varieties
\begin{equation}
    \label{O bars in intro}
    \{0\} = \overline{\scrO}_0 \subset \overline{\scrO}_1 \subset \cdots \subset \overline{\scrO}_r = \p^+,
\end{equation}
where $r$ denotes the real rank of $G_\R$.

If $L$ is an irreducible unitary highest weight module of $G_\R$, then its associated variety $\mathcal{AV}(L)$ equals one of the orbit closures $\overline{\scrO}_k$ in~\eqref{O bars in intro}, and its associated cycle~\cite{Vogan91} takes the form
\[
\mathcal{AC}(L) = m_L \, [\overline{\scrO}_k]
\]
for some positive integer multiplicity $m_L$.
This yields the following expression for the Bernstein degree:
\begin{equation}
    \label{NOT general Bdeg m deg}
    \Deg L = m_L \cdot \deg \overline{\scrO}_k,
\end{equation}
where $\deg$ denotes the degree as an algebraic variety.

\subsection{The Nishiyama--Ochiai--Taniguchi result}

The main result obtained by Nishiyama--Ochiai--Taniguchi in their paper~\cite{NOT} is an elegant formula for Bernstein degree in the following \emph{dual pair} setting.
Given a positive integer $k$, consider one of the three pairs 
\begin{equation}
    \label{dual pair setting}
    (G_\R, H_\R(k)) = (\U(p,q), \: \U(k)), \quad (\operatorname{Mp}(2n , \R), \: \O(k)), \quad \text{or} \quad (\O^*(2n), \: \Sp(k)),
\end{equation}
where $\operatorname{Mp}(2n, \R)$ denotes the metaplectic double cover of $\Sp(2n, \R)$.
In each dual pair setting~\eqref{dual pair setting}, let $\widehat{H}(k)$ denote the set of irreducible representations of the compact group $H_\R(k)$ (equivalently, of its complexification $H(k)$).
The irreducible representations $\sigma \in \widehat{H}(k)$ parametrize a family of irreducible unitary highest weight modules of $G_\R$, in the sense that there is a map
\begin{equation}
    \label{sigma to lambda}
    \sigma \mapsto L_{\la(\sigma)}
\end{equation}
where $L_{\la(\sigma)}$ is the simple $\g$-module with highest weight $\la(\sigma)$.
(See Table~\ref{table:main info} for details.)
Upon restricting to those $\sigma$'s for which $L_{\la(\sigma)} \neq 0$, the map~\eqref{sigma to lambda} is injective.
Moreover, for $G_\R = \U(p,q)$ or $\operatorname{Mp}(2n, \R)$, all irreducible unitary highest weight modules arise in this way as some $L_{\la(\sigma)}$.

For $\sigma \in \widehat{H}(k)$, it is known (see~\cite{EW}*{Thm.~6} or~\cite{Joseph}) that $\mathcal{AV}(L_{\la(\sigma)}) = \overline{\scrO}_{\min\{k,r\}}$, and hence by~\eqref{NOT general Bdeg m deg} there is some positive integer $m_k(\sigma)$ such that
\begin{equation}
    \label{Deg mk sigma}
    \Deg L_{\la(\sigma)} = m_k(\sigma) \cdot \deg \overline{\scrO}_{\min\{k,r\}}.
\end{equation}
The main result of Nishiyama--Ochiai--Taniguchi~\cite{NOT}*{Thm.~B} states that for sufficiently small $k$, this multiplicity $m_k(\sigma)$ is given by the dimension of the irreducible representation $U_\sigma$ of $H_\R(k)$
labeled by $\sigma$:
\begin{equation}
    \label{NOT main result}
    \text{If $k \leq r$, then } \Deg L_{\la(\sigma)} = \dim U_\sigma \cdot \deg \overline{\scrO}_k.
\end{equation}

\subsection{Main result of this paper}

The Bernstein degree formula~\eqref{NOT main result} is limited to the range $k \leq r$.
In this paper, we generalize~\eqref{NOT main result} by introducing a uniform combinatorial interpretation that is valid for all $k$.
In particular, we replace $\dim U_\sigma$ in~\eqref{NOT main result} by the cardinality of a certain set $\Qsigma$ of semistandard tableaux,
and we replace $\deg \overline{\scrO}_k$ by the cardinality of a set $\PP$ of plane partitions.
(See Definitions~\ref{def:T(sigma)} and~\ref{def:R}.)
Our main result is the following theorem.

\begin{theorem}
    \label{thm:main result in intro}
    Let $G_\R$ be one of the groups in the dual pair setting~\eqref{dual pair setting}, where $k$ is any positive integer.
    \textup{(}In the $G_\R = \operatorname{Mp}(2n, \R)$ case, the result remains conjectural within a certain range of $k$-values; see details in Theorem~\ref{thm:main result}.\textup{)}
    For all $\sigma \in \widehat{H}(k)$, we have
    \[
    \Deg L_{\la(\sigma)} = \# \Big( \Qsigma \times \PP \Big).
    \]
\end{theorem}

In particular, the cardinality of $\Qsigma$ equals the multiplicity $m_k(\sigma)$ in~\eqref{Deg mk sigma}, and the tableaux in $\Qsigma$ satisfy certain constraints on the entries in their initial columns (see Definition~\ref{def:T(sigma)}).
A key feature of our construction is that as $k$ increases, these constraints on $\Qsigma$ progressively relax.
As a result, the sets $\Qsigma$ interpolate between linear bases of two distinguished representations:
\begin{equation}
    \label{size Q}
    \# \Qsigma = 
    \begin{cases}
    \dim U_\sigma & \text{if $k \leq r$},\\
    \dim F_{\la(\sigma)} & \text{if $k \geq s$},
    \end{cases}
\end{equation}
where $F_{\la(\sigma)}$ is the simple $K$-module with highest weight $\la(\sigma)$, and where $s$ is the threshold beyond which every $L_{\la(\sigma)}$ is free as a $\C[\p^+]$-module.
(See Table~\ref{table:main info}.)
Indeed, when $k \leq r$, our set $\Qsigma$ (up to a constant shift in tableau entries) recovers a well-known tableau set modeling a basis for $U_\sigma$: 
in the $\U(p,q)$ case, this is the set of \emph{rational tableaux} given by Stembridge~\cite{Stembridge}, and in the other two cases, this is the set of \emph{orthogonal} or \emph{symplectic} tableaux given by Proctor~\cite{ProctorRSK}.
On the other hand, when $k \geq s$, the initial column constraints on $\Qsigma$ vanish, meaning that the sets $\Qsigma$ become full sets of semistandard tableaux that model linear bases for irreducible representations of general linear groups (the group $K$ in our setting).
In this way, our construction of $\Qsigma$ sheds light on the mysterious gap $r < k < s$ left open in \eqref{NOT main result}.

The cardinality of $\PP$ equals the degree of the orbit closure $\overline{\scrO}_{\min\{k,r\}}$ in~\eqref{Deg mk sigma}, and the plane partitions in $\PP$ are bounded by $k$ and contained inside a certain diagram $D_k$ (see Definition~\ref{def:R}).
In the dual pair setting, $\overline{\scrO}_{\min\{k,r\}}$ is a determinantal variety inside a space of generic, symmetric, or skew-symmetric matrices~\eqref{O bar det var};
consequently, this fact $\#\PP = \deg \overline{\scrO}_{\min\{k,r\}}$ is likely known to experts, due to the extensive work in the 1990s on determinantal rings and lattice paths.
(The bijection $\Theta$ in Definition~\ref{def:Theta} connects the plane partition perspective with that of lattice paths.)
It is worth mentioning that the classical plane partition formulas for $\#\PP$ (see Proposition~\ref{prop:enumerate R}) are much simpler than the original degree formulas due to Giambelli--Thom--Porteous~\cites{Giambelli,HarrisTu,Porteous}.
Within the vast literature concerning plane partitions, we highlight Proctor's work~\cite{ProctorPPs} relating various families of bounded plane partitions to weight bases for representations of classical Lie groups.

\subsection{Method of proof}

Via Howe duality (see Section~\ref{sec:Howe}), the unitary highest weight modules $L_{\la(\sigma)}$ can be viewed as \emph{modules of covariants} $M_\sigma$ (see~\eqref{MOCs} below) for the complex classical groups $H(k)$ from~\eqref{dual pair setting}.
The identification $L_{\la(\sigma)} \cong M_\sigma$, as graded modules over the ring $\C[\p^+]$, follows from viewing $L_{\la(\sigma)}$ as the $U_\sigma$-isotypic multiplicity space in the Howe duality decomposition~\eqref{Howe decomp}.
In our previous work~\cite{EricksonHunzikerMOC2024}*{Thm.~3.2}, we introduced a combinatorial model for $M_\sigma$
by organizing Jackson's~\cite{Jackson} linear basis of standard monomials into a \emph{Stanley decomposition} (see Definition~\ref{def:Stanley decomp}).
In particular, we expressed $M_\sigma$ as a direct sum of Stanley spaces indexed by a set $\mathcal{J}_k(\sigma)$ consisting of combinatorial objects we call \emph{jellyfish} (see Definition~\ref{def:jellyfish}).
Very roughly speaking, a jellyfish is an ordered pair $(T, \mathbf{F})$ where $T$ is a semistandard tableau of shape $\sigma$, and $\mathbf{F}$ is a union of $k$ nonintersecting lattice paths whose endpoints are compatible (in a technical sense) with the initial column of $T$.
Section~\ref{sec:jellyfish} in the present paper is a self-contained summary of our results from~\cite{EricksonHunzikerMOC2024}.

The two invariants in the opening paragraph, namely the Gelfand--Kirillov dimension $\operatorname{Dim} X$ and the Bernstein degree $\Deg X$, have analogues in the setting of graded modules:
these are the \emph{Krull dimension} $\operatorname{Kdim} M$ and the \emph{multiplicity} $e(M)$, respectively.
Thus, in translating between this paper and our previous paper~\cite{EricksonHunzikerMOC2024}, we have $\Deg L_{\la(\sigma)} = e(M_\sigma)$.
In turn, it is a general fact (see Lemma~\ref{lemma:Bdeg Stanley decomp}) that for any Stanley decomposition of a graded module $M$, the multiplicity $e(M)$ equals the number of Stanley spaces of maximum Krull dimension.
Thus by identifying the subset $\widehat{\mathcal{J}}_k(\sigma)$ of jellyfish with maximum Krull dimension, we obtain the following (see Lemma~\ref{lemma:Deg equals J hat}):
\[
\Deg L_{\la(\sigma)} = e(M_\sigma) = \# \widehat{\mathcal{J}}_k(\sigma).
\]
The heart of the argument, leading immediately to Theorem~\ref{thm:main result in intro}, occurs in Lemma~\ref{lemma:Theta tilde}, where we establish a bijection
\begin{equation}
    \label{bijection intro}
    {\rm id} \times \Theta : \Qsigma \times \PP \longrightarrow \widehat{\mathcal{J}}_k(\sigma),
\end{equation}
where $\Theta$ takes a plane partition in $\PP$ to its corresponding union of lattice paths.
In this way, upon restriction to the $k \leq r$ case, we produce a new and very different proof of the Nishiyama--Ochiai--Taniguchi result~\eqref{NOT main result}, in addition to extending it to arbitrary values of $k$.

We mentioned above that in the case $G_\R = \Mp(2n, \R)$, Theorem~\ref{thm:main result in intro} remains partially conjectural, namely in the range $r < k < s$.
(See Conjecture~\ref{conj:main}, which we have verified extensively by computer.)
In the $\U(p,q)$ and $\O^*(2n)$ cases, our construction of the jellyfish in~\cite{EricksonHunzikerMOC2024} was guided by the standard monomial results in Jackson~\cite{Jackson}, whereas the $\operatorname{Mp}(2n, \R)$ case (being considerably more delicate) fell outside the scope of Jackson's methods.
Hence in the $\operatorname{Mp}(2n, \R)$ case, our proof of Theorem~\ref{thm:main result in intro} for $k \leq r$ and $k \geq s$ appeals to~\eqref{NOT main result}--\eqref{size Q}.

\subsection{Beyond the dual pair setting}

The ultimate goal of this research program is to generalize Theorem~\ref{thm:main result in intro} uniformly outside the dual pair setting.
This goal requires finding analogues for the sets $\Qsigma$ and $\PP$, for all groups $G_\R$ of Hermitian type (including the exceptional groups).
In this paper, we take the first step toward this goal by generalizing the set $\PP$ as follows.

Recall that in the dual pair setting detailed above (see also Definition~\ref{def:R}), we defined $\PP$ as the set of plane partitions bounded by $k$ and contained in a certain diagram $D_k$ (either a rectangle or a staircase).
By redefining these underlying diagrams $D_k$ in terms of the Hasse diagram of the poset $\Phi(\p^+)$ of positive noncompact roots (see Definition~\ref{def:Dk}), we obtain analogous diagrams $D_k$ across all Hermitian types.
We refer the reader to Table~\ref{table:Dk} (page~\pageref{table:Dk}) for the depictions of these $D_k$'s in all types.
In this way, we generalize not only the motivating property $\# \PP = \deg \overline{\scrO}_k$, but also the combinatorial interpretation of the Hilbert series of the coordinate ring $\C[ \overline{\scrO}_k ]$ in terms of $\PP$.
In the following proposition (which we prove in Section~\ref{sec:open}, restated as Proposition~\ref{prop:Dk general}), given a plane partition $P \in \PP$, the statistic $c(P)$ is obtained by summing local differences over all the entries in~$P$; see~\eqref{c(P)}.
We write $|D_k|$ to denote the number of boxes in the diagram $D_k$.

\begin{prop}
    \label{prop:Dk in intro}
    Let $(G_\R, K_\R)$ be an irreducible Hermitian symmetric pair of noncompact type, and let $r$ denote the real rank of $G_\R$.
    For all $1 \leq k \leq r$, we have the Hilbert series
    \[
    H_{\C[\overline{\scrO}_k]}(t) = \frac{\sum_{P \in \PP} t^{c(P)}}{(1-t)^{\dim \p^+ - |D_k|}}.
    \]
\end{prop}

\noindent (As mentioned above, in the dual pair setting, Proposition~\ref{prop:Dk in intro} can be easily aligned with the equivalent result obtained in the 1990s~\cites{Sturmfels,Conca94,Herzog,GhorpadeKrattenthalerPfaffians,BrunsHerzog} for generic, symmetric, and skew-symmetric determinantal rings, via well-known bijections between plane partitions and families of nonintersecting lattice paths.)
We conclude Section~\ref{sec:open} by observing analogues of the Nishiyama--Ochiai--Taniguchi result~\eqref{NOT main result} where $G_\R$ is one of the exceptional groups $\operatorname{E}_6$ or $\operatorname{E}_7$;
see Table~\ref{table: NOT exceptional} (page~\pageref{table: NOT exceptional}).
These striking examples suggest the existence of a new analogue of Howe duality for the exceptional groups.



\section{Howe duality}
\label{sec:Howe}

\subsection{Hermitian symmetric pairs}

\label{sub:HS}

Suppose $G_\R$ is a connected real reductive group, with real rank $r$.
Further suppose that $K_\R\subset G_\R$ is a maximal compact subgroup such that $G_\R/K_\R$ is an irreducible Hermitian symmetric space of noncompact type.
Let $\g_{\R} = \k_{\R} \oplus \p_{\R}$ be a 
Cartan decomposition of the Lie algebra of $G_\R$, and let 
$\g = \k \oplus \p$ be the corresponding decomposition
of the complexified Lie algebra. 
From the general theory of Hermitian symmetric pairs, there exists a distinguished element $h_0 \in \mathfrak{z}(\k)$ such that $\operatorname{ad} h_0$ acts on $\g$ with eigenvalues $0$ and $\pm 1$.  
We thus have a triangular decomposition 
 \[
    \g = \p^- \oplus \k \oplus \p^+,
 \]
where $\p^{\pm} = \{ x \in \g : [h_0, x] = \pm x\}$.
The subalgebra $\mathfrak{q} = \k \oplus \p^+$ is a maximal parabolic subalgebra of $\g$, with Levi subalgebra $\k$ and abelian nilradical $\p^+$.
Upon fixing a Cartan subalgebra $\mathfrak{h}$ of both $\g$ and $\k$, we write $\Phi$ for the root system of $(\g,\mathfrak{h})$, and $\g_\alpha$ for the root space corresponding to the root $\alpha \in \Phi$.
Then we put
\begin{align}
\label{Phi p+}
\begin{split}
    \Phi(\p^+) &\coloneqq \{\alpha \in \Phi : \g_\alpha \subseteq \p^+\},\\
    S & \coloneqq \C[\p^+] = \C[z_{ij}],
\end{split}
\end{align}
where the $z_{ij}$'s are the matrix coordinate functions on $\p^+$ (see Table~\ref{table:GR}).

Let $K$ be the complexification of $K_\R$.
Then $\k$ is the Lie algebra of $K$, and the adjoint action of $\k$ on
$\g$ exponentiates to a $K$-action.
In particular, $K$ acts on $\p^+$.
The $K$-orbits in $\p^+$ are  
$\scr{O}_0:=\{0\}$ and
$\scrO_k:=K\cdot(e_{\gamma_1}+\dots +e_{\gamma_k})$ for $1\leq k\leq r$, where $e_{\gamma_i}\in \g_{\gamma_i}$ denotes a root vector corresponding to $\gamma_i$, and where $\gamma_1, \ldots, \gamma_r$ are Harish-Chandra's set of strongly orthogonal roots.
Recall from~\eqref{O bars in intro} that the closures of the $K$-orbits in $\p^+$ form a chain of algebraic varieties
\[
    \{0\} = \overline{\scrO}_0 \subset \overline{\scrO}_1 \subset \cdots \subset \overline{\scrO}_r = \p^+.
\]

\begin{table}[t]
    \begin{center}
\begin{tblr}{colspec={|Q[m,c]|Q[m,c]|Q[m,c]|Q[m,c]|Q[m,c]|Q[m,c]|Q[m,c]|Q[m,c]|},stretch=1.5}

\hline

$G_\R$ & $\g_\R$ & $\g$ & $K$ & $\g = \left\{\left[\begin{smallmatrix}A&B\\C&D\end{smallmatrix}\right] : \ldots \right\}$ & $\k = \left\{\left[\begin{smallmatrix}A&0\\0&D\end{smallmatrix}\right]\right\}$ & $\p^+ =  \left\{\left[\begin{smallmatrix}
    0 & B \\ 0 & 0
\end{smallmatrix}\right]\right\}$ \\ \hline[2pt]

$\U(p,q)$ & $\mathfrak{u}(p,q)$ & $\gl_{p+q}$ & $\GL_p \times \GL_q$ & no constraints & $\gl_p \oplus \gl_q$ & $\operatorname{M}_{p,q}$ \\ \hline

$\operatorname{Mp}(2n,\R)$ & $\sp(2n,\R)$ & $\sp_{2n}$ & $\widetilde{\GL}_n$ & {$A=-D^t$, \\ $B=B^t$, \\ $C=C^t$} & $\gl_n$ & $\operatorname{SM}_n$ \\ \hline

$\O^*(2n)$ & $\so^*(2n)$ & $\so_{2n}$ & $\GL_n$ & {$A=-D^t$, \\ $B=-B^t$, \\ $C=-C^t$} & $\gl_n$ & $\operatorname{AM}_n$ \\ \hline

\end{tblr}
\end{center}
    \caption{The real reductive groups $G_\R$ arising in the dual pair setting; see Table~\ref{table:main info} for the corresponding complex classical groups $H(k)$.
    Note that $\widetilde{\GL}_n$ denotes the double cover $\widetilde{\GL}_n \coloneqq \{(g,s)\in \GL_n \times  \, \C^{\times} : \det(g) = s^2\}$.}
    \label{table:GR}
\end{table}

\subsection{Howe duality}
\label{sub:Howe duality}

In this paper, we focus on the three groups in Table~\ref{table:GR}, namely $G_\R = \U(p,q)$, $\operatorname{Mp}(2n, \R)$, or $\O^*(2n)$, defined as follows:
\begin{align*}
\U(p,q)& \coloneqq \left\{g\in \GL(p+q,\C) : g \begin{pmatrix}  I_p &0 \\   0 & -I_q\end{pmatrix}g^* =  \begin{pmatrix}  I_p &0 \\   0 & -I_q\end{pmatrix}\right\},\\[5pt]
\Sp(2n,\R)&\coloneqq\left\{g\in \GL(2n,\C) : g \begin{pmatrix}  0 & I_n \\  -I_n& 0\end{pmatrix}g^t =  \begin{pmatrix}  0 & I_n \\  -I_n & 0\end{pmatrix}\right\}\cap\U(n,n),\\[5pt]
\operatorname{Mp}(2n, \R) & \coloneqq \text{the metaplectic double cover of $\Sp(2n, \R)$},\\[5pt]
\O^*(2n) & \coloneqq \left\{g\in \GL(2n,\C) : g \begin{pmatrix}  0 & I_n \\  I_n & 0\end{pmatrix}g^t =  \begin{pmatrix}  0 & I_n \\   I_n & 0\end{pmatrix}\right\}\cap\U(n,n),
\end{align*}
where $I_n$ denotes the $n \times n$ identity matrix.
Note that we have 
$\U(p,q) \cap \U(p+q) \cong \U(p)\times \U(q)$, and $\Sp(2n,\R)\cap \U(2n) \cong \U(n)$, and $\O^*(2n)\cap \U(2n)\cong \U(n)$,
where in the last two cases $\U(n)$ 
is embedded block-diagonally as follows: 
\[
\left\{\begin{pmatrix}
    a & 0 \\
    0 & (a^{-1})^{t}
\end{pmatrix}
: a \in \U(n)\right\}\cong \U(n).
\]
It follows that
$\Sp(2n,\R)\subseteq \operatorname{SL}(2n,\C)$ and $\O^*(2n) \subseteq \operatorname{SL}(2n,\C)$. 
For this reason, many authors write $\operatorname{SO}^*(2n)$ to denote $\O^*(2n)$.
From now on, we denote the \emph{complex} classical groups by $\GL_n \coloneqq \GL(n, \C)$ and $\O_n \coloneqq \O(n, \C)$, and $\Sp_{2n} \coloneqq \Sp(2n, \C)$.

\begin{remark} 
Our definition of $\Sp(2n,\R)$ above
differs from the more standard definition 
by a Cayley transform as follows.
Let $\mathbf{c}\in \operatorname{SL}(2n,\C)$ be the matrix given by
\[
\mathbf{c} \coloneqq \frac{1}{\sqrt{2}}
\begin{pmatrix}I_n &\sqrt{-1}I_n\\\sqrt{-1}I_n& I_n
\end{pmatrix}.
\]
Then $\mathbf{c} \Sp(2n,\R)\mathbf{c}^{-1}\subseteq \operatorname{SL}(2n,\R)$, with equality if and only if $n=1$.
\end{remark}

Let $k$ be a positive integer.
Corresponding to each group $G_\R$ in Table~\ref{table:GR} is a compact group $H_\R(k)$, giving rise to the following \emph{compact dual pairs} $(G_\R, H_\R(k))$:
\[
(\U(p,q), \: \U(k)), \qquad (\operatorname{Mp}(2n, \R), \: \O(k)), \qquad (\O^*(2n), \: \Sp(k)).
\]
Let $H(k)$ denote the complexification of $H_\R(k)$, given in~\eqref{W table} below.
The group $H(k)$ acts naturally on the following space $W$ of complex matrices, as follows:
\begin{equation}
\label{W table}
\begingroup
\renewcommand{\arraystretch}{1.5}
\begin{array}{|c|c|c|c|}
\hline
G_\R & H(k) & W & \text{$H(k)$-action on $W$} \\ \Xhline{2pt}
\U(p,q) & \GL_k & \operatorname{M}_{p,k} \oplus \operatorname{M}_{k,q} & h \cdot (X,Y) = (Xh^{-1}, \: hY) \\ \hline
\operatorname{Mp}(2n,\R) & \O_k & \operatorname{M}_{k,n} & h \cdot X = hX \\ \hline
\O^*(2n) & \Sp_{2k} & \operatorname{M}_{2k,n} & h \cdot X = hX \\ \hline
\end{array}
\endgroup
\end{equation}
This induces the usual action by $H(k)$ on $\C[W]$, via $(h \cdot f)(w) = f(h^{-1} \cdot w)$.
Let $\mathcal{D}(W)$ denote the Weyl algebra of polynomial-coefficient differential operators on $\C[W]$, and let $\mathcal{D}(W)^{H(k)}$ denote the subalgebra of $H(k)$-invariant operators. 
Then the complexified Lie algebra $\g$ of $G_\R$ embeds in $\mathcal{D}(W)^{H(k)}$ as a Lie subalgebra, and is spanned by a generating set of the associative algebra $\mathcal{D}(W)^{H(k)}$.
Let 
\[
\omega: \g \longrightarrow \mathcal{D}(W)^{H(k)}
\]
denote the injective homomorphism;
then $\g$ acts on $\C[W]$ via its image $\omega(\g)$ in $\mathcal{D}(W)^{H(k)}$.
For details, see the appendix in~\cite{EricksonHunzikerMOC2024}, where we write down explicitly the homomorphism $\omega$ and the action of $\g$ by differential operators.

Let $\widehat{H}(k)$ denote the set of irreducible rational representations of $H(k)$, up to equivalence.
It is standard to view the elements $\sigma \in \widehat{H}(k)$ as partitions (or pairs of partitions, for $\GL_k$), as given in Table~\ref{table:main info}.
We recall that a \emph{partition} is a finite, weakly decreasing sequence of positive integers.
Given a partition $\sigma = (\sigma_1, \ldots, \sigma_m)$, we write $\ell(\sigma) \coloneqq m$ for its \emph{length}.
We write $\sigma'_j \coloneqq \#\{ i : \sigma_i \geq j \}$.
When $H(k) = \GL_k$, we will generalize the notion of a partition by considering decreasing $k$-tuples $\sigma$ of (possibly nonpositive) integers.
Any such $k$-tuple can be expressed uniquely as a pair $\sigma = (\sigma^+, \sigma^-)$, where $\sigma^+$ and $\sigma^-$ are true partitions: 
in particular, $\sigma^+$ consists of the positive coordinates of $\sigma$, and $\sigma^-$ is the partition obtained by negating and reversing the negative coordinates of $\sigma$.
For example, if $\sigma = (6,3,3,2,0,0,-1,-3,-5)$, then $\sigma^+ = (6,3,3,2)$ and $\sigma^- = (5,3,1)$.

Let $U_\sigma$ denote a model for the irreducible representation $\sigma \in \widehat{H}(k)$.
If $H(k) = \Sp_{2k}$ or $\GL_{k}$, then $U_\sigma$ is the irreducible representation of $H(k)$ with highest weight $\sigma$ (in standard coordinates).
For $H(k) = \O_k$, which is not connected, the situation is more subtle; see~\cite{GW}*{pp.~438--9} for a detailed construction of the representations $U_\sigma$.
If $\sigma \in \widehat{H}(k)$,
then $\mathcal{D}(W)^{H(k)}$ and hence $\g$
acts on the multiplicity space
\begin{equation}
    \label{M sigma new}
    M_\sigma \coloneqq {\rm Hom}_{H(k)}(U_\sigma^*, \C[W]).
\end{equation}
The action of an operator $D \in \mathcal{D}(W)^{H(k)}$ on $\varphi \in M_\sigma$ is given by
$D \cdot \varphi = D \circ \varphi$.
The algebra $\mathcal{D}(W)^{H(k)}$ acts irreducibly on $M_\sigma$, and hence $\g$ does as well.
Due to the canonical isomorphisms
\begin{equation}
    \label{MOCs}
    M_\sigma \cong (\C[W] \otimes U_\sigma)^{H(k)} \cong \Big\{ \text{$H(k)$-equivariant polynomial functions from $W$ to $U_\sigma$} \Big\},
\end{equation}
one can identify $M_\sigma$ with the \emph{module of covariants} on the right-hand side of~\eqref{MOCs}; see~\cite{EricksonHunzikerMOC2024}*{\S4}.

Roughly speaking, a \emph{$(\g,K)$-module} is a 
complex vector space carrying representations of both $\g$ and $K$, such that $K$ acts locally finitely and the actions of $\g$
and $K$ are compatible; see~\cite{SchmidNotes}*{Def.~3.2.3} for details.
In the dual pair setting, the ``classical'' $\k$-action on $M_\sigma$ is given by the differential of the right matrix multiplication by $K$ on $W$, but this classical action is \emph{not} the action obtained by restricting the $\g$-action to~$\k$.
Instead, we should view $M_\sigma$ as a $\k$-module on which the classical $\k$-action is ``twisted'' by multiplying by a suitable power of the determinant.
Hence, writing $M_\sigma^{\text{class.}}$ to denote the $\k$-module with the classical action, we will view $M_\sigma$ as the $\k$-module $M_\sigma^{\text{class.}} \otimes F_{-kc\zeta}$, where $F_{-kc\zeta}$ is a certain one-dimensional $\k$-module;
in particular, $c$ is a certain constant intrinsic to $G_\R / K_\R$ (see~\cite{BaiHunziker}*{\S2}), and $\zeta$ is the unique fundamental weight of $\g$ that is orthogonal to the compact roots.
In this way, the twisted $\k$-action on $M_\sigma$ integrates to a $K$-action in a manner compatible with the $\g$-action, and hence 
$M_\sigma$ can be viewed as a $(\g, K)$-module.
Further, $M_\sigma$ is a highest weight $\g$-module, which follows from work of Harish-Chandra~\cites{HC55,HC56} since $G_\R/K_\R$ is Hermitian symmetric.
We denote this highest weight by $\la(\sigma)$.
To write down $\la(\sigma)$ in coordinates, we introduce the following shorthand for a weight of $\g$:
\begin{align}
\label{arrow notation}
\begin{split}
    \overscriptrightarrow{\sigma} &\coloneqq \text{vector obtained from $\sigma$ by padding with 0's on the right},\\
    \overscriptleftarrow{\sigma} &\coloneqq \text{vector obtained by reversing $\overscriptrightarrow{\sigma}$},\\
    \mu - k &\coloneqq \text{vector obtained from $\mu$ by subtracting $k$ from every coordinate}.
    \end{split}
\end{align}
In Table~\ref{table:main info}, we write weights for $\g$ in standard epsilon coordinates on the Cartan subalgebra consisting of diagonal matrices.
In particular, a weight for $\g = \gl_{p+q}$ is written as $(a_1, \ldots, a_p \mid b_1, \ldots, b_q)$;
a weight for $\g = \sp_{2n}$ or $\so_{2n}$ is written as $(a_1, \ldots, a_n)$.

\begin{theorem}[Howe duality; see~\cite{Howe89}, \cite{KV}]
\label{thm:Howe duality}
Assume one of the three settings in Tables~\ref{table:GR} and~\ref{table:main info}, with $W$ as in~\eqref{W table}.
We have the following multiplicity-free decomposition as a $(\g,K) \times H(k)$-module:
\begin{equation}
    \label{Howe decomp}
    \C[W] \cong \bigoplus_{\sigma \in \Sigma} 
    M_\sigma
   \otimes U_{\sigma}^*,
\end{equation}
where $M_\sigma$ is the multiplicity space in~\eqref{M sigma new}, and where $\Sigma \coloneqq \{ \sigma \in \widehat{H}(k) : M_\sigma \neq 0\}$.
Furthermore, the map $\sigma \mapsto \lambda(\sigma)$ is an injective map from $\Sigma$ into the set of $\k$-dominant integral weights, such that as a $(\g, K)$-module
\begin{equation}
\label{covariants = L lambda}
 M_\sigma \cong L_{\la(\sigma)}
 \coloneqq \textup{simple $\g$-module with highest weight $\la(\sigma)$},
\end{equation}
and $L_{\la(\sigma)}$ is unitarizable with respect to $\g_\R$.
\end{theorem}

\begin{table}[t]
    \centering
    \resizebox{\linewidth}{!}
{\begin{tblr}{colspec={|Q[m,c]|Q[m,c]|Q[m,c]|Q[m,c]|Q[m,c]|Q[m,c]|Q[m,c]|},stretch=1.5}

\hline

$G_\R$ & $H(k)$ & $\widehat{H}(k) = \{ \sigma : \ldots \}$ & {$\Sigma =$ \\ $\{ \sigma \in \widehat{H}(k) : \ldots \}$} & $\lambda(\sigma)$ & $r$ & $s$ \\

\hline[2pt]

$\U(p,q)$ & $\GL_k$ & $\ell(\sigma^+) + \ell(\sigma^-) \leq k$ & {$\ell(\sigma^+) \leq q$, \\ $\ell(\sigma^-) \leq p$\phantom{,}} & $\left( -\overscriptleftarrow{\sigma^{-}} - k \; \big| \; \overscriptrightarrow{\sigma^+} \, \right)$ & $\min\{p,q\}$ & $p+q-1$ \\

\hline

$\operatorname{Mp}(2n, \R)$ & $\O_k$ & $\sigma'_1 + \sigma'_2 \leq k$ & $\ell(\sigma) \leq n$ & $- \overscriptleftarrow{\sigma} -\frac{k}{2}$ & $n$ & $2n-1$ \\

\hline

$\O^*(2n)$ & $\Sp_{2k}$ & $\ell(\sigma) \leq k$ & $\ell(\sigma) \leq n$ & $- \overscriptleftarrow{\sigma} -k$ & $\lfloor n/2 \rfloor$ & $n-1$ \\

\hline

\end{tblr}}
    \caption{Companion table to Table~\ref{table:GR}; see Theorem~\ref{thm:Howe duality}.
    The parameter $r$ is the real rank of $G_\R$, and $s$ is the sharp bound such that $k \geq s$ implies that $L_{\la(\sigma)}$ is free for all $\sigma \in \widehat{H}(k)$.
    The arrow notation for $\la(\sigma)$ is the shorthand given in~\eqref{arrow notation}.
    }
    \label{table:main info}
\end{table}

\section{The set \texorpdfstring{$\Qsigma$}{Q sigma}}
\label{sec:Q sigma}

\subsection{Semistandard tableaux}

Recall the definition of a \emph{partition} from Section~\ref{sub:Howe duality}.
We adopt the shorthand $(a^m) \coloneqq (a, \ldots, a)$, and we write $0$ to denote the empty partition, which has length~$0$.
The \emph{Young diagram} associated to a partition $\sigma$ is an arrangement of boxes in left-justified rows, such that the $i$th row from the top contains $\sigma_i$ boxes.
Note that $\sigma'_j \coloneqq \#\{ i : \sigma_i \leq j \}$ can be viewed as the number of boxes in column $j$ of the Young diagram of $\sigma$.

A \emph{semistandard Young tableau} is a Young diagram in which the boxes are filled with positive integer entries, weakly increasing across each row and strictly increasing down each column.
We write closed intervals of integers using the standard notation $[m,n] \coloneqq \{x \in \mathbb{Z} : m \leq x \leq n \}$, and we adopt the usual shorthand $[n] \coloneqq [1, n]$.
The \emph{shape} of a tableau is the partition corresponding to its underlying Young diagram.
The empty tableau $\varnothing$ is the tableau with shape $0$.
We write
\[
\SSYT(\sigma, n) \coloneqq \Big\{\text{semistandard Young tableaux of shape $\sigma$, with entries taken from $[n]$} \Big\}.
\]
Below, we show an example of a partition $\sigma$ along with a tableau $T \in \SSYT(\sigma, n)$:
\begin{equation}
    \label{T example}
    \ytableausetup{smalltableaux,centertableaux}
    \sigma = (3,2,1,1) \quad \leadsto \quad \ytableaushort{126,33,5,6} \in \SSYT(\sigma,n), \text{ where $n \geq 6$.}
\end{equation}
Note that if $\ell(\sigma) > n$, then we automatically have $\SSYT(\sigma, n) = \varnothing$.
On the other hand, if $\sigma = 0$ is the empty partition, then $\SSYT(0, n) = \{ \varnothing \}$ is the singleton consisting of the empty tableau.
Given a tableau $T$, we write
\begin{align*}
    T_1 & \coloneqq \text{set of entries in column $1$ of $T$},\\
    T_1 \uplus T_2 & \coloneqq \text{multiset of entries in columns $1$ and $2$ of $T$},
\end{align*}
where we count columns from the left.
For the tableau $T$ displayed in~\eqref{T example}, we have $T_1 = \{1,3,5,6\}$ and $T_1 \uplus T_2 = \{\!\!\{1, 2, 3, 3, 5, 6 \}\!\!\}$.
In the $G_\R = \U(p,q)$ setting, where $\widehat{H}(k)$ consists of pairs $\sigma = (\sigma^+, \sigma^-)$, a tableau $T$ of shape $\sigma$ is an ordered pair $T = (T^+, T^-)$ of tableaux such that $T^\pm$ has shape $\sigma^\pm$.

\subsection{The set \texorpdfstring{$\Qsigma$}{Q sigma}}

In Definition~\ref{def:T(sigma)} below, we define $\Qsigma$ as a subset of a full set $\Tsigma$ of semistandard tableaux of shape $\sigma$.
As mentioned above, for $G_\R = \U(p,q)$, each ``tableau'' $T \in \Tsigma$ is actually an ordered pair $T = (T^+, T^-)$.

\begin{dfn}
    \label{def:T(sigma)}
    Assume one of the three settings in Table~\ref{table:main info}.
    For $\sigma \in \widehat{H}(k)$, define
    \[
    \Qsigma \coloneqq \left\{ T \in \Tsigma : \begin{array}{l} 
\alpha_i(T) < i, \\
\text{for all positive integers $i$ such that} \\
k-r < i \leq k
\end{array}
\right\},
    \]
 where $\Tsigma$ and $\alpha_i(T)$ are defined as follows:

\begin{center} 
\begin{tblr}{colspec={|Q[m,c]|Q[m,c]|Q[m,c]|},stretch=1.5}

\hline

$G_\R$ & $\Tsigma$ & $\alpha_i(T)$ \\

\hline[2pt]

$\U(p,q)$ &  {$\phantom{\times}\SSYT(\sigma^+, q)$\\$ \times \SSYT(\sigma^-, p)$} & {$\phantom{+}\;\#\{x \in T^+_1 : x < q-k+i \}$\\[1ex]$+ \;\#\{y \in T^-_1 : y < p - k + i \}$}  \\

\hline

$\operatorname{Mp}(2n, \R)$ & $\SSYT(\sigma, n)$ & $\#\{x \in T_1 \uplus T_2 : x < n - k + i \}$ \\

\hline

$\O^*(2n)$ & $\SSYT(\sigma, n)$ & $\#\{x \in T_1 : x < n-1-2k+2i \}$ \\

\hline

\end{tblr}

\end{center}
\end{dfn}

Comparing the $\Sigma$ column in Table~\ref{table:main info} with Definition~\ref{def:T(sigma)}, we observe that $\sigma \in \Sigma$ if and only if $\mathcal{T}(\sigma) \neq \varnothing$.

\begin{remark}
    \label{rem:i range}
    For $G_\R = \operatorname{Mp}(2n, \R)$ or $\O^*(2n)$, it is straightforward to verify by Definition~\ref{def:T(sigma)} that the constraint $\alpha_i(T) < i$ is trivial if $i  \leq k-r$;
    hence there is no harm in omitting the lower bound on $i$, so that $1 \leq i \leq k$.
    By contrast, for $G_\R = \U(p,q)$, the constraint $\alpha_i(T) < i$ is not necessarily trivial for $i \leq k-r$, so we cannot omit this lower bound without changing the definition of $\Qsigma$.
    On the other hand, if one wishes to \emph{avoid} checking trivial constraints (see Example~\ref{ex:Qsigma} below), then one need only check positive integers $i$ such that
    \[
    \begin{cases}
        \min\{2k - s, \: k-r\} < i \leq \ell(\sigma^+) + \ell(\sigma^-) & \text{if } G_\R = \U(p,q),\\
        2k - s < i \leq \sigma'_1 + \sigma'_2 & \text{if } G_\R = \operatorname{Mp}(2n, \R),\\
        2k - s < i \leq \ell(\sigma) & \text{if } G_\R = \O^*(2n).
    \end{cases}
    \]
\end{remark}  

\begin{example}\
\label{ex:Qsigma}
    \begin{enumerate}
        
    \item Let $G_\R = \U(p,q)$, where $(p,q) = (3,5)$ and $k=5$.
    By Table~\ref{table:main info}, we have $r = 3$, so $k-r=2$.
    In $\alpha_i(T)$ in Definition~\ref{def:T(sigma)}, the conditions are $x < q - k + i = i $ and $y < p - k + i = i - 2$.
    Thus, the set $\Qsigma$ consists of those tableaux $T \in \Tsigma$ such that $
    \#\{x \in T^+_1 : x < i \} + \#\{y \in T^-_1 : y < i-2 \} < i$, for $i=3,4,5$:
       \begin{align*}
            \#\{x \in T^+_1 : x < 3 \} + \#\{y \in T^-_1 : y < 1\} & < 3,\\
            \#\{x \in T^+_1 : x < 4 \} + \#\{y \in T^-_1 : y < 2\} & < 4,\\
            \#\{x \in T^+_1 : x < 5 \} + \#\{y \in T^-_1 : y < 3\} & < 5.
        \end{align*}
    (Note that the first constraint is trivial, since it is satisfied by any semistandard tableau.)
    Thus we have
    \ytableausetup{boxsize=1.2em}
    \[
    \Qsigma = \left\{ T \in \Tsigma : (T^+, T^-) \text{ has initial columns } \left( \ytableaushort[\scriptstyle]{\ast,\ast,{\geq 5}}, \ytableaushort[\scriptstyle]{1,2,\none} \right) \text{ or } \left( \ytableaushort[\scriptstyle]{\ast,\ast,{\geq 4}}, \ytableaushort[\scriptstyle]{1,3,\none} \right) \text{ or } \left( \ytableaushort[\scriptstyle]{\ast,\ast,\ast}, \ytableaushort[\scriptstyle]{2,3,\none} \right) \right\},
    \]
    where the stars denote arbitrary entries, and where we keep only the topmost $\ell(\sigma^+)$ and $\ell(\sigma^-)$ many boxes in the left- and right-hand column, respectively.

    \item Let $G_\R = \operatorname{Mp}(2n, \R)$, where $n=4$ and $k=5$.
    By Table~\ref{table:main info}, we have $r = 4$, so $k-r=1$.
    In $\alpha_i(T)$ in Definition~\ref{def:T(sigma)}, the condition is $x < n - k + i = i-1$.
    Thus, $\Qsigma$ consists of those tableaux $T \in \Tsigma$ such that $\# \{ x \in T_1 \uplus T_2 : x < i-1 \} < i$, for $i=2,3,4,5$:
    \begin{align*}
        \#\{x \in T_1 \uplus T_2 : x < 1 \} & < 2, \\
        \#\{x \in T_1 \uplus T_2 : x < 2 \} & < 3, \\
        \#\{x \in T_1 \uplus T_2 : x < 3 \} & < 4, \\
        \#\{x \in T_1 \uplus T_2 : x < 4 \} & < 5.
    \end{align*}
    (Note that the first two constraints are trivial.)
    Thus we have
    \[
    \Qsigma = \left\{ T \in \Tsigma : \text{$T$ has initial columns } \, \ytableaushort[\scriptstyle]{1 \ast,2 {\geq 4},3} \,\text{ or }\, \ytableaushort[\scriptstyle]{1 \ast,2 {\geq 3},4} \,\text{ or }\, \ytableaushort[\scriptstyle]{\ast \ast,3 \ast,4} \right\},
    \]
    where we keep only the topmost $\sigma'_1$ and $\sigma'_2$ many boxes in the first and second column, respectively.

    \item Let $G_\R = \O^*(2n)$, where $n = 7$ and $k=4$.
        By Table~\ref{table:main info}, we have $r=3$, so $k-r=1$.
        In $\alpha_i(T)$ in Definition~\ref{def:T(sigma)}, the condition on $x \in T_1$ is $x < n-1-2k+2i = 2i-2$.
        Thus, the set $\Qsigma$ consists of those tableaux $T \in \Tsigma$ such that $\#\{x \in T_1 : x < 2i-2 \} < i$, for $i=2,3,4$:
    \begin{align*}
        \#\{ x \in T_1 : x < 2 \} & < 2,\\
        \#\{ x \in T_1 : x < 4 \} & < 3,\\
        \#\{ x \in T_1 : x < 6 \} & < 4.
    \end{align*}
    (Note that the first constraint is trivial.)
    Thus we have
        \[
        \Qsigma = \left\{ T \in \Tsigma : T \text{ has initial column \,  $\ytableaushort[\scriptstyle]{\ast, \ast, {\geq4}, {\geq6}}$} \, \right\},
        \]
        where we keep only the topmost $\ell(\sigma)$ many boxes.
    \end{enumerate}
    
\end{example}
    
    \begin{lemma}
        \label{lemma:U sigma F lambda}
        Assume one of the three settings in Tables~\ref{table:GR} and~\ref{table:main info}, where $\sigma \in \widehat{H}(k)$.
        Let $F_{\la(\sigma)}$ denote the simple $\k$-module with highest weight $\la(\sigma)$.

        \begin{enumerate}[label=\textup{(\alph*)}]
            \item If $k \leq r$, then $\# \Qsigma = \dim U_\sigma$.

            \item If $k \geq s$, then $\# \Qsigma = \dim F_{\la(\sigma)}$.
        \end{enumerate}
    \end{lemma}

    \begin{proof}\

    (a) Suppose $k \leq r$.
        Then by Definition~\ref{def:T(sigma)}, we have $T \in \Qsigma$ if and only if $\alpha_i(T) < i$ for all $1 \leq i \leq k$.
        Equivalently, in the definition of $\alpha_i(T)$ and in the condition $\alpha_i(T) < i$, one can replace all the strict inequalities $<$ by weak inequalities $\leq$, and shift the range of $i$-values so that $0 \leq i \leq k-1$.
        But at $i=0$, the resulting constraint is trivial, so it may be omitted;
        at $i=k$, the resulting constraint is automatically satisfied by the fact that $\sigma \in \widehat{H}(k)$, so there is no harm in including it.
        Hence we may still use the range $1 \leq i \leq k$.
        
        \begin{itemize}
            \item Let $G_\R = \U(p,q)$.
            Then $U_\sigma$ is the irreducible representation of $H(k) = \GL_k$ with highest weight $\sigma$.
            Let $\overline{\mathcal{Q}}(\sigma)$ denote the set of all pairs $(U^+, U^-) \in \SSYT(\sigma^+, k) \times \SSYT(\sigma^-, k)$ such that $\#\{x \in U^+_1 : x \leq i\} +  \#\{y \in U^-_1 : y \leq i\} \leq i$, for all $ 1 \leq i \leq k$.
            Then $\overline{\mathcal{Q}}(\sigma)$ is Stembridge's set of \emph{rational tableaux} in~\cite{Stembridge}*{Def.~2.2}.
            By~\cite{Stembridge}*{Prop.~2.4(a)}, we have $\dim U_\sigma = \# \overline{\mathcal{Q}}(\sigma)$.
            To complete the proof, we exhibit the bijection $\Qsigma \longrightarrow \overline{\mathcal{Q}}(\sigma)$ given by $(T^+, T^-) \mapsto (\overline{T}^+, \overline{T}^-)$, where $\overline{T}^+$ (resp., $\overline{T}^-$) is obtained by subtracting $q-k$ (resp., $p-k$) from every entry in $T^+$ (resp., $T^-$).
            
            \item Let $G_\R = \operatorname{Mp}(2n, \R)$.
            Then $U_\sigma$ is the irreducible representation of $H(k) = \O_k$ labeled by $\sigma$.
            Let $\overline{\mathcal{Q}}(\sigma)$ denote the set of all $U \in \SSYT(\sigma, k)$ such that $\#\{x \in U_1 \uplus U_2 : x \leq i \} \leq i$, for all $1 \leq i \leq k$.
            Then $\overline{\mathcal{Q}}(\sigma)$ is Proctor's set of \emph{orthogonal tableaux} in~\cite{ProctorRSK}*{p.~30}.
            By~\cite{Proctor94}*{Thm.~3.1}, we have $\dim U_\sigma = \# \overline{\mathcal{Q}}(\sigma)$.
            To complete the proof, we exhibit the bijection $\Qsigma \longrightarrow \overline{\mathcal{Q}}(\sigma)$ given by $T \mapsto \overline{T}$, where $\overline{T}$ is obtained by subtracting $n-k$ from every entry in $T$.
            
            \item Let $G_\R = \O^*(2n)$.
            Then $U_\sigma$ is the irreducible representation of $H(k) = \Sp_{2k}$ with highest weight $\sigma$.
            Let $\overline{\mathcal{Q}}(\sigma)$ denote the set of all $U \in \SSYT(\sigma, 2k)$ such that $\#\{x \in U_1 : x \leq 2i \} \leq i$, for all $1 \leq i \leq k$.
            Then $\overline{\mathcal{Q}}(\sigma)$ is Proctor's set of \emph{symplectic tableaux} in~\cite{ProctorRSK}*{p.~30}.
            By~\cite{Proctor94}*{Thm.~3.1}, we have $\dim U_\sigma = \# \overline{\mathcal{Q}}(\sigma)$.
            To complete the proof, we exhibit the bijection $\Qsigma \longrightarrow \overline{\mathcal{Q}}(\sigma)$ given by $T \mapsto \overline{T}$, where $\overline{T}$ is obtained by subtracting $n-2k$ from every entry in $T$.
        \end{itemize}

        (b) Suppose $k \geq s$.
        For each group $G_\R$, it is straightforward to check directly from Definition~\ref{def:T(sigma)} that the conditions $\alpha_i(T) < i$ are all trivial for $k-r < i \leq k$.
        Hence we have
        \begin{equation}
            \label{Q = T}
            \Qsigma = \mathcal{T}(\sigma).
        \end{equation}
        It remains to show that $\# \mathcal{T}(\sigma) = \dim F_{\la(\sigma)}$.
        If $G_\R = \U(p,q)$, then by Table~\ref{table:GR} we have $\k = \gl_p \oplus \gl_q$.
        Letting $F_p(\mu)$ denote the simple $\gl_p$-module with highest weight $\mu$ (and similarly for $\gl_q$), we have
        \begin{align*}
            \#\mathcal{T}(\sigma) = \#\SSYT(\sigma^-, p) \cdot \# \SSYT(\sigma^+, q)) &= \dim F_p(\overscriptrightarrow{\sigma^-}) \cdot \dim F_q(\overscriptrightarrow{\sigma^+}) \\
            &= \dim \left(F_p(\overscriptrightarrow{\sigma^-}) \otimes F_p(k^p) \right)^* \cdot \dim F_q(\overscriptrightarrow{\sigma^+}) \\
            &= \dim \left(F_p(- \overscriptleftarrow{\sigma^-} - k) \otimes F_q(\overscriptrightarrow{\sigma^+}) \right) \\
            &= \dim F_{(- \overscriptleftarrow{\sigma^-} - k \; \mid \; \overscriptrightarrow{\sigma^+})} = \dim F_{\la(\sigma)}.
        \end{align*}
        If $G_\R = \operatorname{Mp}(2n, \R)$ or $\O^*(2n)$, then by Table~\ref{table:GR} we have $\k = \gl_n$, and since $\mathcal{T}(\sigma) = \SSYT(\sigma, n)$, the proof is just a simpler version of the $\U(p,q)$ case above.
    \end{proof}
    
    The following  lemma gives an equivalent (case-by-case) criterion for $T \in \Qsigma$ which will be useful in proving Proposition~\ref{prop:enumerate Qsigma}.
    We write $T_{j,\ell}$ to denote the entry in row $j$ and column $\ell$ of~$T$.
    For a set $S$ of integers, we write $S_{(j)}$ to denote the $j$th smallest element of $S$.
    Hence we have equivalent notation $T_{j,\ell} = (T_\ell)_{(j)}$.

    \begin{lemma}
    \label{lemma:explicit criteria}
        Let $\sigma \in \widehat{H}(k)$ in one of the settings in Table~\ref{table:main info}. 
        Let $T \in \Tsigma$ as in Definition~\ref{def:T(sigma)}.

        \normalfont

        \begin{enumerate}
            \item Let $G_\R = \U(p,q)$. 
            Let $\overline{T}^+$ (resp., $\overline{T}^-$) be the tableau obtained by adding $k-q$ (resp., $k-p$) to each entry of $T^+$ (resp., $T^-$).
            Define the sign
            \begin{equation*}
                \epsilon \coloneqq \begin{cases}
                    + & \text{if } p \leq q,\\
                    - & \text{if } p > q.
                \end{cases}
            \end{equation*}
            We have $T \in \Qsigma$ if and only if
            \[
            \overline{T}^{\, -\epsilon}_{1,1} \geq 1, \text{ and } \overline{T}^{\,\epsilon}_{j,1} \geq ([k] \setminus \overline{T}^{\, -\epsilon}_1)_{(j)}
            \]
            for all positive integers $j$ such that $k- \min\{p,q\} < j \leq \ell(\sigma^{\epsilon})$.

        \item Let $G_\R = \operatorname{Mp}(2n, \R)$.
        Define the interval $I \coloneqq [ n-k+1, \: n]$.

        We have $T \in \Qsigma$ if and only if
        \[
        T_{1,1} \geq n-k+1, \text{ and } T_{j,2} \geq (I \setminus T_1)_{(j)}
        \]
        for all positive integers $j$ such that $1 \leq j \leq \sigma'_2$.

        \item Let $G_\R = \O^*(2n)$.     
            We have $T \in \Qsigma$ if and only if
            \[
            T_{j,1} \geq n + 2(j-k) - 1 \quad \text{for all } 1 \leq j \leq \ell(\sigma).
            \]
        \end{enumerate}
    \end{lemma}

    \begin{proof}\

        (1) Let $G_\R = \U(p,q)$.
        Suppose $p \leq q$; then we have $r=p$, and  $T^{\,\epsilon} = T^+$ and $T^{- \epsilon} = T^-$.
        By Definition~\ref{def:T(sigma)}, we have $T \in \Qsigma$ if and only if $\#\{ x \in \overline{T}^+_1 : x < i \} + \# \{ y \in \overline{T}^-_1 : y < i \} < i$ for all $\max\{0,k-p\} < i \leq k$.
        Equivalently, shifting the index $i$ by 1, we have $T \in \Qsigma$ if and only if
        \begin{equation}
            \label{equivalent i condition GL}
            \#\{ x \in \overline{T}^+_1 : x \leq i \} + \# \{ y \in \overline{T}^-_1 : y \leq i \} \leq i, \quad \text{for all } \max\{0, \: k-p\} \leq i \leq k-1.
        \end{equation}
        If $\overline{T}^-_{1,1} \leq 0$, then $k < p$, and so $i=0$ lies in the interval $\max\{0, k-p\} \leq i \leq k$;
        but then~\eqref{equivalent i condition GL} is violated for $i=0$, and thus $T \notin \Qsigma$.
        Therefore, $T \in \Qsigma$ implies $\overline{T}^-_{1,1} \geq 1$, as claimed in this lemma.
        Consequently, we have $\overline{T}^-_1 \subseteq [k]$, and so~\eqref{equivalent i condition GL} is equivalent to the following:
        \begin{align}
        \label{y to z GL}
        \begin{split}
            \#\{ x \in \overline{T}^+_1 : x \leq i \} & \leq i - \# \{ y \in \overline{T}^-_1 : y \leq i \} \\
            &= \# \{ z \in ([k] \setminus \overline{T}^-_1) : z \leq i \}, \quad \text{for all $\max\{0, \: k-p\} \leq i \leq k-1$.}
            \end{split}
        \end{align}
        If $k \leq p$, then~\eqref{y to z GL} occurs if and only if
        \begin{equation}
            \label{final condition GL}
            \overline{T}^+_{j,1} \geq ([k] \setminus \overline{T}^-_1)_{(j)}
        \end{equation}
        for all $1 \leq j \leq \#\overline{T}^+_1 = \ell(\sigma^+)$, and we are done.
        If $k > p$, then $\overline{T}^-_1 \subseteq [k-p+1, \: k]$, and thus $[k] \setminus \overline{T}^-_1$ contains the initial interval $[k-p]$;
        in this case,~\eqref{y to z GL} occurs if and only if~\eqref{final condition GL} holds for all $k-p < j \leq \ell(\sigma^+)$, so we are done.
        This completes the proof for the $p \leq q$ case;
        the $p > q$ case is identical, upon interchanging the roles of $p$ and $q$, and the roles of $\overline{T}^+$ and $\overline{T}^-$.

        (2) Let $G_\R = \operatorname{Mp}(2n,\R)$.
        By Definition~\ref{def:T(sigma)} and by Remark~\ref{rem:i range}, we have $T \in \Qsigma$ if and only if $\#\{x \in T_1 \uplus T_2 : x < n-k+i\} < i$ for all $1 \leq i \leq k$.
        Equivalently, $T \in \Qsigma$ if and only if
        \begin{equation}
        \label{equivalent O}
            \#\{x \in T_1 : x \leq n-k+i \} + \#\{y \in T_2 : y \leq n-k+i \} \leq i, \quad \text{for all } 0 \leq i \leq k-1. 
        \end{equation}
        If $T_{1,1} \leq n-k$, then $T$ violates the $i=0$ case in~\eqref{equivalent O}, and we have $T \notin \Qsigma$; hence $T \in \Qsigma$ implies $T_{1,1} \geq n-k+1$, as claimed in this lemma.
        Thus we may assume that both $T_1$ and $T_2$ are subsets of $I = [n-k+1, \: n]$, and so~\eqref{equivalent O} is equivalent to the following:
        \begin{align*}
            \#\{ y \in T_2 : y \leq n-k+i \} & \leq i - \# \{ x \in T_1 : x \leq n-k+i \} \\
            &= \# \{ z \in (I \setminus T_1) : z \leq n - k + i \}, \quad \text{for all $1 \leq i \leq k-1$,}
        \end{align*}
        which occurs if and only if $T_{j,2} \geq (I \setminus T_1)_{(j)}$ for all $1 \leq j \leq \# T_2 = \sigma'_2$.

        (3) Let $G_\R = \O^*(2n)$.
        In Definition~\ref{def:T(sigma)}, replacing $i$ by $j$, we have $\alpha_j(T) < j$ if and only if $\#\{x \in T_1 : x < n - 1 - 2k + 2j \} < j$, if and only if
        \begin{equation}
            \label{Tj1 Sp}
            T_{j,1} \geq n - 1 - 2k + 2j = n + 2(j-k) - 1.
        \end{equation}
        Therefore, by Definition~\ref{def:T(sigma)}, we have $T \in \Qsigma$ if and only if~\eqref{Tj1 Sp} is true for all positive integers $j$ such that $k-r < j \leq \#T_1 = \ell(\sigma)$.
        Equivalently, by Remark~\ref{rem:i range}, we can take $ 1\leq j \leq \ell(\sigma)$.
    \end{proof}

In the following proposition, we adopt the convention that $\binom{a}{b} = 0$ if $a$ or $b$ is negative.

\begin{prop}
    \label{prop:enumerate Qsigma}

    Assume one of the three settings in Table~\ref{table:main info}, where $\sigma \in \widehat{H}(k)$.
    Recall the set $\Qsigma$ from Definition~\ref{def:T(sigma)}.

    \normalfont

    \begin{enumerate}
    
    \item Let $G_\R = \U(p,q)$.
    We have
    \[
    \#\Qsigma = \det \left[ \binom{\sigma_i - i + j + c_j + d_j}{\sigma_i - i + j + c_j} \right]_{1 \leq i,j \leq k},
    \]
    where 
    \[
    c_j \coloneqq \begin{cases}
        0 & \textup{if } j \leq k-r,\\
        \sigma_1^{-\epsilon} & \textup{if }j>k-r
        \end{cases} 
        \quad
        \textup{and}
        \quad
        d_j \coloneqq -1 + \begin{cases}
        R & \textup{if } j \leq k-r,\\
        \min\{k,r\} & \textup{if } j > k-r,
    \end{cases} \vspace{1ex}
    \]
    with $r \coloneqq \min\{p,q\}$ and $R \coloneqq \max \{p,q\}$ and the sign $\epsilon$ as in Lemma~\ref{lemma:explicit criteria}(1).

    \item Let $G_\R = \operatorname{Mp}(2n, \R)$.
    We have
    \[
    \#\Qsigma = \sum_C \det \left[ \binom{\sigma_i - i + j + n - a_j(C) - 1}{\sigma_i - i + j - 1} \right]_{1 \leq i,j \leq \sigma'_2} 
    \]
    with the sum ranging over all $\sigma'_1$-element subsets $C \subseteq [\max\{1, \: n-k+1\}, \: n]$, and
    \[
    a_j(C) \coloneqq \max\{ C_{(j)}, (I \setminus C)_{(j)} \},
    \]
    where $I = [n-k+1, n]$.

    \item Let $G_\R = \O^*(2n)$.
    We have
    \[
    \#\Qsigma = \det \left[ \binom{\sigma_i - i + j + n - a_j}{\sigma_i - i + j} \right]_{1 \leq i,j \leq k},
    \]
    where $a_j := \max\{1, \: n+2(j-k) - 1\}$. 
    \end{enumerate}
\end{prop}

\begin{proof}
    We appeal to a well-known result of Gessel--Viennot~\cite{GV85}*{\S6}; see also the exposition by Krattenthaler~\cite{Krattenthaler}*{\S10.13}.
    In particular, let $\la$ and $\mu$ be partitions of length at most $k$, such that $\mu_i \leq \la_i$ for all $1 \leq i \leq k$.
    Let $\la / \mu$ denote the skew shape obtained by removing the Young diagram of $\mu$ from that of $\la$, where the upper-left corners of the diagrams coincide.
    Let $\mathbf{a} = (a_1, \ldots, a_k)$ and $\mathbf{b} = (b_1, \ldots, b_k)$ be weakly increasing sequences of integers such that each $a_i \leq b_i$.
    Define
    \[
    \SSYT(\la/\mu, \mathbf{a}, \mathbf{b}) \coloneqq \left\{ \begin{array}{l} \text{semistandard Young tableaux of shape $\la/\mu$,}\\
    \text{where all entries in row $i$ lie in $[a_i, b_i]$}
    \end{array}\right\}.
    \]
    Note that $\SSYT(\la,n)$ is the special case where $\mu = 0$ and $\mathbf{a} = (1,\ldots,1)$ and $\mathbf{b} = (n, \ldots, n)$.
    By~\cite{Krattenthaler}*{Thm.~10.13.3}, we have
    \begin{equation}
        \label{Krattenthaler result}
        \#\SSYT(\la / \mu, \mathbf{a}, \mathbf{b}) = \det \left[ \binom{\lambda_i - \mu_j - i + j + b_i - a_j}{\lambda_i - \mu_j - i + j} \right]_{1 \leq i,j \leq k}.
    \end{equation}

    (1) Let $G_\R = \U(p,q)$, and suppose $p \leq q$.
    Thus $r=p$, $R=q$, and $\sigma^{-\epsilon} = \sigma^-$.
    (Just as in Lemma~\ref{lemma:explicit criteria}(1), the $p > q$ case is obtained by switching the roles of $p$ and $q$, and of $T^+$ and $T^-$.)
    Let $T \in \Tsigma$, and define $\overline{T}^\pm$ as in Lemma~\ref{lemma:explicit criteria}(1).
    We construct a single (possibly skew) semistandard tableau $\overline{T}'$ as follows:
    \begin{itemize}
        \item Start with $\overline{T}^+$, observing that
        \begin{equation}
        \label{T bar plus}
            \overline{T}^+ \in \SSYT(\sigma^+/0, \mathbf{a}, \mathbf{b}), \quad \text{where each } a_i = k-q+1 \text{ and } b_i = k.
        \end{equation}

        \item Adjoin the top-justified columns $C_{\sigma^-_1}, \ldots, C_3, C_2, C_1$ to the left of $\overline{T}^+$, where each
        \[
        C_\ell \coloneqq [k] \setminus \overline{T}^-_\ell.
        \]

        \item If $k > p$, then remove the topmost $k-p$  rows from each column $C_\ell$ above.
    \end{itemize}
    Note that the columns $C_{\sigma^-}, \ldots, C_1$ form a semistandard tableau $C$ with entries in $[k]$; this follows from the same arguments used in the proof of Lemma~\ref{lemma:explicit criteria}(1), or from a result of Stembridge~\cite{Stembridge}*{Lemma 2.5(b)}.
    Moreover, since (by Lemma~\ref{lemma:explicit criteria}) $\overline{T}^-$ consists of entries in the interval $[\max\{1, \: k-p+1\}, k]$, if $j \leq k-p$ then the $j$th entry in every column $C_\ell$ equals $j$;
    this is the reason for having removed the rectangular block of $k-p$  rows and $\sigma^-_1$  columns from the upper-left (since these entries are all forced, and have no effect on enumerating tableaux).
    On the other hand, if $j > k-p$, then by Lemma~\ref{lemma:explicit criteria}(1), the entry in row $j$ of $\overline{T}^+_1$ is no less than the entry in row $j$ of $C_1$, and therefore the entire tableau $\overline{T}'$ is semistandard.
    By construction, $\overline{T}'$ has shape $\la / \mu$, where $\la = \sigma + ((\sigma^-_1)^k)$, and where $\mu = ((\sigma^-_1)^{k-p})$ if $k > p$ and $\mu = 0$ otherwise.
    Note that if $j \leq k-p$, then row $j$ in $\overline{T}'$ is the same as row $j$ in $\overline{T}^+$, whose entries are bounded below by $a_j = k-q+1$ and above by $b_j = k$, by~\eqref{T bar plus}.
    If $j > k-p$, then the initial entries in row $j$ come from the $C_\ell$'s, whose entries are bounded below by $a_j = \max\{1, \: k-p+1\}$ and above by $b_j = k$.
    Hence the map $T \mapsto \overline{T}'$ gives a bijection
    \[
    \Qsigma \longrightarrow \SSYT(\la/\mu, \mathbf{a}, \mathbf{b}),
    \]
    where, for all $1 \leq i,j \leq k$,
    \[
    \la_i = \sigma_i + \sigma^-_1, 
    \qquad 
    \mu_j = \begin{cases}
        \sigma^-_1 & \text{if } j \leq k-p,\\
        0 & \text{otherwise},
    \end{cases}
    \qquad 
    a_j = \begin{cases}
        k-q+1 & \text{if } j \leq k-p,\\
        \max\{1, \: k-p+1\} & \text{otherwise},
    \end{cases}
    \]
    and $b_i = k$.
    The result follows from~\eqref{Krattenthaler result}, upon setting $c_j \coloneqq \sigma^-_1 - \mu_j$ and $d_j \coloneqq k - a_j$.

    (2) Let $G_\R = \operatorname{Mp}(2n, \R)$.
    Given $T \in \Qsigma$, we decompose $T$ into its initial column $T_1$ and the tableau $T_{>1}$ consisting of its remaining columns $T_2, \ldots, T_{\sigma_1}$.
    By Lemma~\ref{lemma:explicit criteria}(2), this map $T \mapsto (T_1, T_{>1})$ gives a bijection
    \[
    \Qsigma \longrightarrow \coprod_{T_1} \left( \{T_1\} \times \left\{T_{>1} \in \SSYT(\sigma-(1^{\sigma'_1}), \: n) : \begin{array}{l} (T_{>1})_{j,1} \geq \max\{(T_1)_{(j)}, (I \setminus T_1)_{(j)} \\[.5ex] \text{for all } 1 \leq j \leq \sigma'_2 \end{array}
    \right\} \right),
    \]
    where the disjoint union ranges over all possible initial columns $T_1$, that is, where $T_1$ is a $\sigma'_1$-element subset of $[\max\{1, \: n-k+1\}, n]$.
    Upon replacing $T_1$ by $C$, we see that the cardinality of each disjoint component in the right-hand side equals $\#\SSYT(\sigma-(1^{\sigma'_1}), \mathbf{a}, \mathbf{b})$, where $a_j = \max\{C_{(j)}, (I\setminus C)_{(j)}\}$ and where every $b_i = n$.
    The result now follows from~\eqref{Krattenthaler result}. 

    (3) Let $G_\R = \O^*(2n)$.
    By Lemma~\ref{lemma:explicit criteria}(3), we have $\Qsigma = \SSYT(\sigma/0, \mathbf{a}, \mathbf{b})$, where $\mathbf{a} = (a_1, \ldots, a_k)$ is given by $a_j = \max\{1, \: n + 2(j-k) - 1\}$, and where $\mathbf{b} = (n, \ldots, n)$.
    The result follows from~\eqref{Krattenthaler result}.
    \end{proof}

    \begin{remark}
        In light of Lemma~\ref{lemma:U sigma F lambda}, we conclude that the determinantal formulas in Proposition~\ref{prop:enumerate Qsigma} must collapse to a simple product formula when $k \leq r$ or $k \geq s$.
        In particular, if $k \leq r$, then $\# \Qsigma = \dim U_\sigma$ is given by the well-known hook-content formulas for $H(k) = \GL_k$, $\O_k$, and $\Sp_{2k}$.
        Likewise, if $k \geq s$, then $\# \Qsigma = \dim F_{\la(\sigma)}$ is given by the hook-content formula for $K = \GL_p \times \GL_q$ or $\GL_n$.
    \end{remark}

\section{The set \texorpdfstring{$\PP$}{Pk}}

\subsection{Plane partitions and the set \texorpdfstring{$\PP$}{Pk}}

We use the term \emph{diagram} to describe any edge-connected union $D$ of unit squares (\emph{boxes}), not necessarily the Young diagram of a partition.
In particular, the following diagrams play an important role.
A \emph{$p \times q$ rectangle} is the Young diagram with shape $(q^p)$.
Following Proctor~\cite{ProctorPPs}, we define an \emph{$n$-staircase} to be the Young diagram with shape $(n,\ldots,3,2,1)$, and a \emph{shifted $n$-staircase} to be the horizontal reflection of an $n$-staircase.

Let $k$ be a positive integer.
In this paper, a \emph{plane partition} bounded by $k$, and contained in a diagram $D$, is a filling of $D$ with entries taken from $[0,k]$, such that the entries weakly increase from left to right, and from bottom to top.
\textbf{NB:} This means that our two-dimensional pictures of plane partitions are 90-degree clockwise rotations of the usual convention in the literature.
Nonetheless, when using the typical \emph{three}-dimensional depiction of a plane partition --- where an entry $t$ is represented by a stack of $t$  unit cubes --- this difference in convention vanishes.
See Figure~\ref{fig:Theta} for an example of the plane partitions in the following definition.

\begin{dfn}
    \label{def:R}
    Assume one of the three settings in Table~\ref{table:main info}, and define the diagram
    \[
    D_k \coloneqq \begin{cases}
       \text{a $(p-k) \times (q-k)$ rectangle} & \text{if $G_\R = \U(p,q)$},\\
       \text{an $(n-k)$-staircase} & \text{if $G_\R = \operatorname{Mp}(2n, \R)$},\\
       \text{a shifted $(n-2k-1)$-staircase} & \text{if $G_\R = \O^*(2n)$}.
    \end{cases}
    \]
    We define the set
    \[
    \PP \coloneqq \{ \text{plane partitions bounded by $k$ and contained in the diagram $D_k$} \}.
    \]
\end{dfn}

See the top row of Figure~\ref{fig:Theta} for examples of the plane partitions in $\PP$.
We observe from Definition~\ref{def:R} that if $k \geq r$ then $D_k = \varnothing$, and thus $\PP = \{ \varnothing \}$, where $\varnothing$ is the empty plane partition.
Thus,
\begin{equation}
    \label{P when k geq r}
    \text{if $k \geq r$ then $\# \PP = 1$}.
\end{equation}

\begin{prop}
    \label{prop:enumerate R}
    Let $\PP$ be as defined in Definition~\ref{def:R}.

    \normalfont

    \begin{enumerate}
    
    \item If $G_\R = \U(p,q)$, then
    \[
    \# \PP = \prod_{i=1}^{p-k} \;\; \prod_{j =1}^{q-k} \;\; \frac{k +i+j-1}{i+j-1}.
    \]

    \item If $G_\R = \operatorname{Mp}(2n, \R)$, then
    \[
    \# \PP = \prod_{1 \leq i \leq j \leq n-k} \frac{k + i + j - 1}{i + j - 1}.
    \]

    \item If $G_\R = \O^*(2n)$, then
    \[
    \# \PP = \prod_{1 \leq i \leq j \leq n-2k-1} \frac{2k + i + j}{i+j}.
    \]
    \end{enumerate}
\end{prop}

\begin{proof}
    Part (1) is the formula for boxed $(p-k) \times (q-k) \times k$ plane partitions due to MacMahon~\cite{MacMahon}*{\S IX}.
    Part (2) is the formula for symmetric boxed $(n-k) \times (n-k) \times k$ plane partitions, conjectured by MacMahon and ultimately proved by Andrews~\cite{Andrews}*{\S2}.
    Part (3) is given by Proctor~\cite{ProctorPPs}*{Thm.~1}, upon setting $q=1$ in his formula (CGH) or (CGI).
\end{proof}

\subsection{Plane partitions, lattice paths, and order complexes}

We next highlight a well-known bijection between plane partitions and families of nonintersecting lattice paths.
To make this precise, we introduce a poset $(\mathbf{P}, \leq)$ whose underlying set $\mathbf{P} \subset \mathbb{N}^2$ is given as follows:
\begin{equation}
    \label{P}
    \mathbf{P} \coloneqq 
    \begin{cases}
        \{(i,j) : 1 \leq i \leq p \text{ and } 1 \leq j \leq q \} & \text{ if $G_\R = \U(p,q)$},\\
        \{(i,j) : 1 \leq i \leq j \leq n \} & \text{ if $G_\R = \operatorname{Mp}(2n, \R)$},\\
        \{(i,j) : 1 \leq i < j \leq n\} & \text{ if $G_\R = \O^*(2n)$}.
    \end{cases}
\end{equation}
For $G_\R = \U(p,q)$ or $\O^*(2n)$, we equip $\mathbf{P}$ with the product order, whereby $(i,j) \leq (i', j')$ if and only if $i \leq i'$ and $j \leq j'$.
For $G_\R = \operatorname{Mp}(2n, \R)$, we reflect this partial order so that $(i,j) \leq (i', j')$ if and only if $i \leq i'$ and $j \geq j'$.
In all cases, we will depict $\mathbf{P}$ so that the minimal element is in the northwest corner; thus each element is covered by the element to its east and the element to its south.
Thus, adapting the terminology from Definition~\ref{def:R}, we observe that our depiction of $\mathbf{P}$ is a $p \times q$ rectangle, an $n$-staircase, and a shifted $(n-1)$-staircase, respectively:

\[
    \begin{tikzpicture}[scale=.2,baseline]
        \foreach \x in {1,...,12}{\foreach \y in {2,...,9}{\node [dot] at (\x,\y) {};}}
        \node[below,align=center] at (6.5,1) {$G_\R = \U(p,q)$,\\$(p,q) = (8,12)$};
    \end{tikzpicture}
    \qquad\qquad
    \begin{tikzpicture}[scale=.2,baseline]
        \foreach \x in {2,...,9}{\foreach \y in {\x,...,9}{\node [dot] at (2+\x,\y) {};}}
        \node[below,align=center] at (6.5,1) {$G_\R = \operatorname{Mp}(2n, \R)$,\\$n=8$};
    \end{tikzpicture}
    \qquad\qquad
    \begin{tikzpicture}[scale=.2,baseline]
        \foreach \x in {2,...,9}{\foreach \y in {\x,...,9}{\node [dot] at (10-\x,\y) {};}}
        \node[below,align=center] at (4.5,1) {$G_\R = \O^*(2n)$,\\$n=9$};
    \end{tikzpicture}
\]

\begin{remark}
    Recall the set $\Phi(\p^+)$ from~\eqref{Phi p+}. 
    When $G_\R = \U(p,q)$ or $\O^*(2n)$, the Lie algebra $\g$ is simply laced, and we have an isomorphism of posets
    $\mathbf{P} \longrightarrow \Phi(\p^+)$ given by
\[
    \begin{cases}
    (i,j) \longmapsto 
    (\overscriptleftarrow{\epsilon_i} \mid - \overscriptrightarrow{\epsilon_j}) & \text{if } G_\R = \U(p,q), \\[2ex]
    (i,j) \longmapsto \overscriptleftarrow{\epsilon_i + \epsilon_j} & \text{if } G_\R = \O^*(2n).
    \end{cases}
\]
We further develop this observation in Section~\ref{sec:open}.
\end{remark}

Recall that a \emph{chain} in $\mathbf{P}$ is a totally ordered subset of $\mathbf{P}$.
An \emph{antichain} is a subset of $\mathbf{P}$ in which any two elements are incomparable.
The \emph{width} of a subset $\mathbf{S} \subseteq \mathbf{P}$ is the cardinality of the largest antichain contained in $\mathbf{S}$.
For a positive integer $k$, the \emph{$k$th order complex} on $\mathbf{P}$ is the abstract simplicial complex
\begin{equation}
    \label{Delta_k}
    \Delta_k(\mathbf{P}) \coloneqq \{ \mathbf{S} \subseteq \mathbf{P} : \mathbf{S} \text{ has width at most $k$} \}.
\end{equation}
A \emph{facet} of $\Delta_k(\mathbf{P})$ is an element of $\Delta_k(\mathbf{P})$ that is maximal with respect to inclusion.
We denote the set of facets by
\begin{equation}
    \label{F hat general}
    \widehat{\mathcal{F}}_k \coloneqq \{ \mathbf{F} : \text{$\mathbf{F}$ is a facet of $\Delta_k(\mathbf{P})$}\}.
\end{equation}
Observe that for all three cases in~\eqref{P}, the width of $\mathbf{P}$ equals the parameter $r$ (see Table~\ref{table:main info}).
Therefore, if $k \geq r$ then $\Delta_k(\mathbf{P})$ is just the power set of $\mathbf{P}$, and $\widehat{\mathcal{F}}_k = \{ \mathbf{P} \}$.

Note that a saturated chain in $\mathbf{P}$ is a \emph{lattice path}, obtained by taking steps in the direction $(1,0)$ or $(0,1)$.
Accordingly, in our diagrams the steps in a lattice path are either south or east.
If the starting point and endpoint coincide, then the lattice path contains only that single point.
Assuming $k \leq r$, it is well known that
\begin{equation}
    \label{F hat as paths}
    \widehat{\mathcal{F}}_k = \{ \mathbf{F} : \mathbf{F} \text{ is a maximal union of $k$  nonintersecting lattice paths in $\mathbf{P}$} \}.
\end{equation}
Note that in all three settings (still assuming $k \leq r$), there are certain points common to every $\mathbf{F} \in \widehat{\mathcal{F}}_k$:
\begin{itemize}
       
    \item Let $G_\R = \U(p,q)$.
    Let $a_1, \ldots, a_k \in \mathbf{P}$ be the points lying on the $k$th antidiagonal $i+j = k+1$, listed from northeast to southwest.
    Let $b_1, \ldots, b_k \in \mathbf{P}$ be obtained by translating $a_1, \ldots, a_k$ (respectively) by $(p-k, \: q-k)$.
    Then every $\mathbf{F}$ contains the disjoint union $\coprod_{t=1}^k (\mathbf{K}_t \sqcup \mathbf{M}_t)$, where $\mathbf{K}_t$ is the horizontal lattice path ending at $a_t$, and $\mathbf{M}_t$ is the vertical lattice path starting at $b_t$.

    \item Let $G_\R = \operatorname{Mp}(2n, \R)$.
    Let $a_1, \ldots, a_k \in \mathbf{P}$ be the points lying on the $k$th antidiagonal, listed from northeast to southwest.
    Then every $\mathbf{F}$ contains the disjoint union $\coprod_{t=1}^k \mathbf{K}_t$, where $\mathbf{K}_t$ is defined in the same way as the $\U(p,q)$ case above.

    \item Let $G_\R = \O^*(2n)$.
    Let $a_1, \ldots, a_k \in \mathbf{P}$ be the points lying on the $2k$th antidiagonal $i+j = 2(k+1)$, listed from northeast to southwest.
    Let $b_1, \ldots, b_k \in \mathbf{P}$ be obtained by translating $a_1, \ldots, a_k$ (respectively) by $(n-2k-1, \: n-2k-1)$ (equivalently, reflecting about the line $i+j=n+1$).
    Then every $\mathbf{F}$ contains the disjoint union $\coprod_{t=1}^k (\mathbf{K}_t \sqcup \mathbf{M}_t)$, defined in the same way as the $\U(p,q)$ case above.
    
    \end{itemize}
    See the bottom row of Figure~\ref{fig:Theta}, where we indicate the $a_t$'s (resp., $b_t$'s) with enlarged dots in the northwest (resp., southeast) corner of $\mathbf{P}$.
    The $\mathbf{K}_t$'s (resp., $\mathbf{M}_t$'s) are the gray paths ending at the $a_t$'s (resp., starting at the $b_t$'s).
    In the following definition, for each plane partition $P \in \PP$ we construct a subset $\Theta(P) \subseteq \mathbf{P}$.

\begin{figure}
    \centering
    \ytableausetup{smalltableaux}
\[
\begingroup
\setlength{\arraycolsep}{3ex}
\begin{array}{rcccc}
    G_\R: & \U(p,q) & \operatorname{Mp}(2n, \R) & \O^*(2n) & \\
    \text{Parameters :} &  (p,q) = (7,9), \; k=3 & n=7, \; k=3 & n=11, \; k=3 & \\[3ex]
    P \in \PP: 
    &
    \ytableaushort{122333,012233,001122,000001}
    &
    \ytableaushort{1233,022,01,0}
    & 
    \ytableaushort{1233,\none 122,\none\none 02,\none\none\none 1}
    &
    \\
    \parbox{3.5cm}{\flushright \phantom{.} \\ \phantom{.} \\ \phantom{.} \\ 3D depiction of $P$ :}
    &
    \begin{tikzpicture}[scale=.3,baseline=(current bounding box.north)]
    \planepartition{{3,3,2,1},{3,3,2,0},{3,2,1,0},{2,2,1,0},{2,1,0,0},{1,0,0,0}}
    \end{tikzpicture} 
    &
    \begin{tikzpicture}[scale=.3, baseline=(current bounding box.north)]
    \planepartition{{3},{3,2},{2,2,1},{1,0,0,0}}
    \end{tikzpicture}
    &
    \begin{tikzpicture}[scale=.3,baseline=(current bounding box.north)]
    \planepartition{{3,2,2,1},{3,2,0},{2,1},{1}}
    \end{tikzpicture}
    &
    \\
    \parbox{3.5cm}{\flushright \phantom{.} \\ \phantom{.} \\ \phantom{.} \\ With $\mathbf{L}_1, \ldots, \mathbf{L}_k$ :}
    &
    \begin{tikzpicture}[scale=.3,baseline=(current bounding box.north)]
    \planepartition{{3,3,2,1},{3,3,2,0},{3,2,1,0},{2,2,1,0},{2,1,0,0},{1,0,0,0}}
    \draw [ultra thick] 
    (0,-.5) ++ (150:3) ++ (210:3) node [circle,fill=black, minimum size = 4pt, inner sep=0pt] {} -- ++ (30:3) -- ++(-30:1) -- ++(30:1) -- ++ (-30:1) -- ++(30:2) -- ++ (-30:2) node [circle,fill=black, minimum size = 4pt, inner sep=0pt] {}
    (0,-1.5) ++ (150:3) ++ (210:3) node [circle,fill=black, minimum size = 4pt, inner sep=0pt] 
    {} -- ++ (30:1) -- ++(-30:1) -- ++(30:1) -- ++ (-30:1) -- ++(30:2) -- ++ (-30:1) -- ++(30:2) -- ++ (-30:1) node [circle,fill=black, minimum size = 4pt, inner sep=0pt] {}
    (0,-2.5) ++ (150:3) ++ (210:3) node [circle,fill=black, minimum size = 4pt, inner sep=0pt] 
    {} -- ++ (-30:1) -- ++(30:1) -- ++(-30:1) -- ++ (30:1) -- ++(-30:1) -- ++ (30:3) -- ++(-30:1) -- ++ (30:1) node [circle,fill=black, minimum size = 4pt, inner sep=0pt] {}
    ;
    \end{tikzpicture} 
    &
    \begin{tikzpicture}[scale=.3, baseline=(current bounding box.north)]
    \planepartition{{3},{3,2},{2,2,1},{1,0,0,0}}
    \draw [ultra thick] 
    (0,-.5) ++ (150:3) ++ (210:1) node [circle,fill=black, minimum size = 4pt, inner sep=0pt] {} -- ++ (30:2) -- ++(-30:1) -- ++(30:1)
    (0,-1.5) ++ (150:3) ++ (210:1) node [circle,fill=black, minimum size = 4pt, inner sep=0pt] 
    {} -- ++ (30:1) -- ++(-30:2) -- ++(30:1)
    (0,-2.5) ++ (150:3) ++ (210:1) node [circle,fill=black, minimum size = 4pt, inner sep=0pt] 
    {} -- ++ (-30:1) -- ++(30:1) -- ++(-30:2) 
    ;
    \end{tikzpicture}
   &
   \begin{tikzpicture}[scale=.3,baseline=(current bounding box.north)]
    \planepartition{{3,2,2,1},{3,2,0},{2,1},{1}}
    \draw [ultra thick] 
    (0,-.5) ++ (150:3) ++ (210:1) node [circle,fill=black, minimum size = 4pt, inner sep=0pt] {} -- ++ (30:2) -- ++(-30:1) -- ++(30:2) -- ++ (-30:3) node [circle,fill=black, minimum size = 4pt, inner sep=0pt] {}
    (0,-1.5) ++ (150:3) ++ (210:1) node [circle,fill=black, minimum size = 4pt, inner sep=0pt] 
    {} -- ++ (30:1) -- ++(-30:1) -- ++(30:1) -- ++ (-30:1) -- ++(30:1) -- ++ (-30:1) -- ++(30:1) -- ++ (-30:1) node [circle,fill=black, minimum size = 4pt, inner sep=0pt] {}
    (0,-2.5) ++ (150:3) ++ (210:1) node [circle,fill=black, minimum size = 4pt, inner sep=0pt] 
    {} -- ++ (-30:1) -- ++(30:1) -- ++(-30:1) -- ++ (30:2) -- ++(-30:2) -- ++ (30:1) node [circle,fill=black, minimum size = 4pt, inner sep=0pt] {}
    ;
    \end{tikzpicture}
    &
   \\
   \parbox{3.5cm}{\flushright \phantom{.} \\ \phantom{.} \\ \phantom{.} \\ $\Theta(P) \in \widehat{\mathcal{F}}_k$ :} 
    &
        \begin{tikzpicture}[scale=.25, baseline=(current bounding box.north)]
        \draw [ultra thick, lightgray] 
        (1,7) -- ++(2,0) (1,6) -- ++(1,0) (9,3) -- ++(0,-2) (8,2) -- ++(0,-1) ;
        \draw[ultra thick] 
        (3,7) -- ++(3,0) -- ++(0,-1) -- ++(1,0) -- ++(0,-1) -- ++(2,0) -- ++(0,-2)
        (2,6) -- ++(1,0) -- ++(0,-1) -- ++(1,0) -- ++(0,-1) -- ++(2,0) -- ++(0,-1) -- ++(2,0) -- ++(0,-1)
        (1,5) -- ++(0,-1) -- ++(1,0) -- ++(0,-1) -- ++(1,0) -- ++(0,-1) -- ++(3,0) -- ++(0,-1) -- ++(1,0);
        \node [circle,fill=black, minimum size = 4pt, inner sep=0pt] at (3,7) {};
        \node [circle,fill=black, minimum size = 4pt, inner sep=0pt] at (2,6) {};
        \node [circle,fill=black, minimum size = 4pt, inner sep=0pt] at (1,5) {};
        \node [circle,fill=black, minimum size = 4pt, inner sep=0pt] at (9,3) {};
        \node [circle,fill=black, minimum size = 4pt, inner sep=0pt] at (8,2) {};
        \node [circle,fill=black, minimum size = 4pt, inner sep=0pt] at (7,1) {};
        \foreach \x in {1,...,9}{\foreach \y in {1,...,7}{\node [dot] at (\x,\y) {};}}
    \end{tikzpicture} 
    &
    \begin{tikzpicture}[scale=.25, baseline=(current bounding box.north)]
        \draw [ultra thick, lightgray] 
        (1,7) -- ++(2,0) (1,6) -- ++(1,0);
        \draw[ultra thick] 
        (3,7) -- ++(2,0) -- ++(0,-1) -- ++(1,0)
        (2,6) -- ++(1,0) -- ++(0,-2) -- ++(1,0)
        (1,5) -- ++(0,-1) -- ++(1,0) -- ++(0,-2);
        \node [circle,fill=black, minimum size = 4pt, inner sep=0pt] at (3,7) {};
        \node [circle,fill=black, minimum size = 4pt, inner sep=0pt] at (2,6) {};
        \node [circle,fill=black, minimum size = 4pt, inner sep=0pt] at (1,5) {};
        \foreach \x in {1,...,7}{\foreach \y in {\x,...,7}{\node [dot] at (\x,\y) {};}}
    \end{tikzpicture}
    &
    \begin{tikzpicture}[scale=.25, baseline=(current bounding box.north)]
        \draw [ultra thick, lightgray] 
        (1,11) -- ++(5,0) (2,10) -- ++(3,0) (3,9) -- ++(1,0) (10,7) -- ++(0,-5) (9,6) -- ++(0,-3) (8,5) -- ++(0,-1);
        \draw[ultra thick] 
        (6,11) -- ++(2,0) -- ++(0,-1) -- ++(2,0) -- ++(0,-3)
        (5,10) -- ++(1,0) -- ++(0,-1) -- ++(1,0) -- ++(0,-1) -- ++(1,0) -- ++(0,-1) -- ++(1,0) -- ++(0,-1)
        (4,9) -- ++(0,-1) -- ++(1,0) -- ++(0,-1) -- ++(2,0) -- ++(0,-2) -- ++(1,0);
        \node [circle,fill=black, minimum size = 4pt, inner sep=0pt] at (6,11) {};
        \node [circle,fill=black, minimum size = 4pt, inner sep=0pt] at (5,10) {};
        \node [circle,fill=black, minimum size = 4pt, inner sep=0pt] at (4,9) {};
        \node [circle,fill=black, minimum size = 4pt, inner sep=0pt] at (10,7) {};
        \node [circle,fill=black, minimum size = 4pt, inner sep=0pt] at (9,6) {};
        \node [circle,fill=black, minimum size = 4pt, inner sep=0pt] at (8,5) {};
        \foreach \x in {2,...,11}{\foreach \y in {\x,...,11}{\node [dot] at (12-\x,\y) {};}}
    \end{tikzpicture} 
    &
\end{array}
\endgroup
\]
    \caption{The bijection $\Theta : \PP \longrightarrow \widehat{\mathcal{F}}_k$ in Definition~\ref{def:Theta} and Lemma~\ref{lemma:Theta}.}
    \label{fig:Theta}
\end{figure}

\begin{dfn} 
    \label{def:Theta}
Let $P \in \PP$, as in Definition~\ref{def:R}.
Recall that $k \geq r$ implies $\PP = \{ \varnothing \}$; in this case, we define $\Theta(\varnothing) \coloneqq \mathbf{P}$.
If $k < r$, then construct $\Theta(P) \subseteq \mathbf{P}$ as follows:
\begin{itemize}
    \item For $1 \leq t \leq k$:

    \begin{itemize}
    
    \item Align the upper-left corner of $P$ with the point $a_t$.
    
    \item Draw the lattice path $\mathbf{L}_t \subseteq \mathbf{P}$ from $a_t$ to $b_t$ (or if $H(k) = \O_k$, from $a_t$ to the main antidiagonal of $\mathbf{P}$) such that $\mathbf{L}_t$ is the southwest boundary of the subdiagram $P_t$ containing those boxes with entries strictly greater than $k - t$.

    \end{itemize}

    \item Define $\Theta(P) \coloneqq \coprod_{t=1}^k \mathbf{N}_t$, where $\mathbf{N}_t \coloneqq \mathbf{K}_t \cup \mathbf{L}_t \cup \mathbf{M}_t$;
    if $G_\R = \operatorname{Mp}(2n, \R)$, then $\mathbf{N}_t \coloneqq \mathbf{K}_t \cup \mathbf{L}_t$.
    \end{itemize}
    
\end{dfn}

For examples of the map $P \mapsto \Theta(P)$ in Definition~\ref{def:Theta}, see the top and bottom row of Figure~\ref{fig:Theta}; 
in the bottom row, the paths $\mathbf{L}_1, \ldots, \mathbf{L}_k$ (from northeast to southwest) are drawn in black.
In Figure~\ref{fig:Theta}, we also show the three-dimensional depiction of $P$, which makes it easy to trace out $\mathbf{L}_1, \ldots, \mathbf{L}_k$ (from top to bottom).

The following lemma is closely related to the bijections given by Krattenthaler~\cite{KrattenthalerMem95}, who used two-row arrays (equivalently, nonnegative integer matrices) instead of plane partitions.
(See in particular Propositions 28, 30, and 33 in~\cite{KrattenthalerMem95}, corresponding to our cases $G_\R = \U(p,q)$, $\O^*(2n)$, and $\operatorname{Mp}(2n, \R)$, respectively.)

\begin{lemma}
    \label{lemma:Theta}
    The map $P \mapsto \Theta(P)$ in Definition~\ref{def:Theta} gives a bijection $\Theta : \PP \longrightarrow \widehat{\mathcal{F}}_k$.
\end{lemma}

\begin{proof}
    (Throughout the proof, suppress any mention of $b_t$ and $\mathbf{M}_t$ in the $G_\R = \operatorname{Mp}(2n, \R)$ case.)
    By construction, $\Theta(P)$ is a family of $k$  lattice paths in $\mathbf{P}$, since each $\mathbf{K}_t$ ends at the starting point $a_t$ of $\mathbf{L}_t$, which ends at the starting point $b_t$ of $\mathbf{M}_t$.
    Since $P_1 \subseteq \cdots \subseteq P_k$, the southwest boundary of $P_{t+1}$ lies weakly southwest of the southwest boundary of $P_{t}$;
    thus, since the starting point $a_{t+1}$ of $\mathbf{L}_{t+1}$ is strictly southwest of the starting point $a_{t}$ of $\mathbf{L}_{t}$, the paths $\mathbf{L}_1, \ldots, \mathbf{L}_t$ are indeed nonintersecting.
    This justifies the disjoint union $\coprod_t \mathbf{N}_t$ written in the last step of Definition~\ref{def:Theta}.
    Maximality follows immediately from the starting and ending points of the $\mathbf{N}_t$'s.
    Hence by~\eqref{F hat as paths}, we have $\Theta(P) \in \widehat{\mathcal{F}}_k$.
    It is clear that distinct plane partitions in $\PP$ yield distinct families $\{ \mathbf{L}_1, \ldots, \mathbf{L}_k\}$, so $\Theta$ is one-to-one.
    
    In the other direction, given $\mathbf{F} \in \widehat{\mathcal{F}}_k$, one can decompose $\mathbf{F}$ into a maximal union of $k$  nonintersecting lattice paths $\mathbf{N}_1, \ldots, \mathbf{N}_k$, by~\eqref{F hat as paths}, such that each $\mathbf{N}_t$ contains $\mathbf{K}_t$ and $\mathbf{M}_t$.
    This leaves nonintersecting lattice paths $\mathbf{L}_t$ from $a_t$ to $b_t$, for $1 \leq t \leq k$.
    Now by superimposing these paths $\mathbf{L}_t$ so that they all start in the upper-left corner of the diagram $D_k$, we obtain the plane partition $\Theta^{-1}(\mathbf{F}) \in \PP$ as follows:
    the entry in each box of $D_k$ equals the number of $\mathbf{L}_t$'s lying to the southwest of that box.
    Hence $\Theta$ is onto, which completes the proof.    
\end{proof}

\section{Stanley decompositions of unitary highest weight modules}
\label{sec:jellyfish}

This section is a self-contained summary of the main result in~\cite{EricksonHunzikerMOC2024}; see Lemma~\ref{lemma:Stanley decomp via jellyfish} below.
This lemma plays a key role in proving our main result in Section~\ref{sec:main result}.

\subsection{Multiplicity and Stanley decompositions}
\label{sub:Stanley decomp}

We recall some fundamental results from commutative algebra concerning finitely generated graded modules over polynomial rings.

\begin{dfn}
    \label{def:Hilbert series}
    Let $M = \bigoplus_{n = 0}^\infty M_n$ be a finitely generated graded $S$-module with Krull dimension $d$, where $S$ is a polynomial ring over $\C$ in determinates of degree 1.
    The \emph{Hilbert--Poincar\'e series} of $M$ is the formal power series
    \[
    H_M(t) \coloneqq \sum_{n = 0}^\infty (\dim M_n) t^n.
    \]

    \end{dfn}

    Let $M$ be as in Definition~\ref{def:Hilbert series}.
    By the Hilbert--Serre theorem (see also~\cite{BrunsHerzog}*{Cor.~4.1.8}), there exists a unique $Q_M(t) \in \mathbb{Z}[t]$ with $Q(1) \neq 0$ such that
    \begin{equation}
        \label{Hilbert series general}
    H_M(t) = \frac{Q_M(t)}{(1-t)^d}.
    \end{equation}

\begin{dfn}
    \label{def:e(M)}
    Let $M$ be as in Definition~\ref{def:Hilbert series}, with reduced Hilbert series given by~\eqref{Hilbert series general}.
    The \emph{multiplicity} of $M$ is the positive integer $e(M) \coloneqq Q_M(1)$.
\end{dfn}

The multiplicity $e(M)$ equals $(d-1)!$ times the leading coefficient of the Hilbert polynomial $P_M(t)$; see the exposition in~\cite{BrunsHerzog}, particularly page 150 and Proposition 4.1.9.
Comparing Definition~\ref{def:e(M)} with~\eqref{Deg = Q(1)}, we see that multiplicity is synonymous with Bernstein degree, in the following sense.
Given a good filtration $\{ X_n \}$ in the context of~\eqref{Deg = Q(1)}, the graded object
\[
\operatorname{gr} X \coloneqq \bigoplus_{n=0}^\infty X_n / X_{n-1}
\]
is a graded module over $S(\mathfrak{g})$, and following~\cite{NOT}*{Def.~1.3}, we define
\begin{equation}
    \label{Deg equals e}
    \Deg X \coloneqq e( \operatorname{gr} X).
\end{equation}

\begin{dfn}[see~\cite{BKU}*{Def.~2.1}]
    \label{def:Stanley decomp}
    Let $M$ be as in Definition~\ref{def:Hilbert series}.
    A \emph{Stanley decomposition} of $M$ is a finite family $(S_i, \eta_i)_{i \in \mathcal{I}}$
    where $\eta_i \in M$ is homogeneous, and $S_i$ is a graded $\C$-algebra retract of $S$ such that $S_i \cap \operatorname{Ann} \eta_i = 0$, and 
    \begin{equation}
    \label{Stanley decomp module}
        M = \bigoplus_{i \in \mathcal{I}} S_i \eta_i
    \end{equation}
    as a graded vector space.
    Each $S_i \eta_i$ in~\eqref{Stanley decomp module} is called the \emph{Stanley space} corresponding to $i \in \mathcal{I}$.
\end{dfn}

By~\cite{HVZ}*{Lemma~1.1}, any module $M$ satisfying the hypotheses of Definition~\ref{def:Hilbert series} admits a Stanley decomposition.
(The terminology refers to  Stanley's work~\cite{Stanley82} on commutative monoids.)
Note that a Stanley decomposition~\eqref{Stanley decomp module} yields the following closed form for the Hilbert series:
\begin{equation}
    \label{HS from SD}
    H_M(t) = \sum_{i \in \mathcal{I}} \frac{t^{\deg \eta_i}}{(1-t)^{{\rm Kdim}(S_i)}},
\end{equation}
where ${\rm Kdim}$ denotes the Krull dimension.
Comparing this with the reduced form of the Hilbert series in~\eqref{Hilbert series general}, we make the following key observation.

\begin{lemma}[\cite{Jahan}*{p.~1018}]
     \label{lemma:Bdeg Stanley decomp}
     Let $M$ be as in Definition~\ref{def:Hilbert series}, and let $(S_i, \eta_i)_{i \in \mathcal{I}}$ be a Stanley decomposition of $M$.
     Then
     \[
     d \coloneqq {\rm Kdim}(M) = \max_{i \in \mathcal{I}} \{ {\rm Kdim}(S_i) \}.
     \]
     Furthermore, defining
     \[
     \widehat{\mathcal{I}} \coloneqq \{ i \in \mathcal{I} : {\rm Kdim}(S_i) = d \},
     \]
     we have
     \[
     e(M) = \# \widehat{\mathcal{I}}.
     \]
 \end{lemma}

 \subsection{Stanley decompositions of unitary highest weight modules}

In this paper, the role of $M$ in~\eqref{Stanley decomp module} is played by the unitary highest weight modules $L_{\la(\sigma)}$ arising in the dual pair setting~\eqref{covariants = L lambda}.
The role of $S$ is played by the ring defined in~\eqref{Phi p+}, namely
\[
S = \C[\p^+] = \C[z_{ij} : (i,j) \in \mathbf{P} ],
\]
where $\mathbf{P}$ is the poset in~\eqref{P}.
Each $L_{\la(\sigma)}$ has the structure of a graded $S$-module, where the grading is induced from the degree filtration coming from $U(\g)$, normalized so that the weight space corresponding to the highest weight $\la$ has degree 0.
In particular, the degree of the $\mu$-weight space is $(\la-\mu)(h_0)$, where $h_0 \in \mathfrak{z}(\k)$ is the distinguished element described at the beginning of Section~\ref{sub:HS}.
Thus by~\eqref{Deg equals e}, the Bernstein degree of $L_{\la(\sigma)}$ is synonymous with its multiplicity as a graded $S$-module:
\begin{equation}
    \label{Deg = e for L}
    \Deg L_{\la(\sigma)} = e( L_{\la(\sigma)}).
\end{equation}

In previous work~\cite{EricksonHunzikerMOC2024}, we obtained Stanley decompositions of all modules $L_{\la(\sigma)}$, for the groups $\U(p,q)$ and $\O^*(2n)$.
The role of the index set $\mathcal{I}$ in~\eqref{Stanley decomp module} is played by a set $\mathcal{J}_k(\sigma)$ consisting of combinatorial objects we call \emph{jellyfish of shape~$\sigma$}.
Note that in~\cite{EricksonHunzikerMOC2024}, we stated this result in terms of the modules of covariants $M_\sigma$ of the classical groups $H(k)$; see~\eqref{M sigma new}--\eqref{MOCs} and~\eqref{Howe decomp} above.
Recall from~\eqref{covariants = L lambda} that (upon viewing $M_\sigma$ as the $\k$-module $M_\sigma^{\text{class.}} \otimes F_{-kc \zeta}$) we have $M_\sigma \cong L_\la(\sigma)$ as $(\g,K)$-modules;
moreover, in~\cite{EricksonHunzikerMOC2024} we observed that, upon shifting the natural grading on $M_\sigma$ via $d \mapsto (d-|\sigma|)/2$, we also have $M_\sigma \cong L_{\la(\sigma)}$ as graded $S$-modules.
Therefore, in stating Lemma~\ref{lemma:Stanley decomp via jellyfish} below (where we record our main result from~\cite{EricksonHunzikerMOC2024}), we apply to $L_{\la(\sigma)}$ our Stanley decompositions of $M_\sigma$.
After the lemma, we will explain the notation in detail.
(In light of Lemma~\ref{lemma:Bdeg Stanley decomp}, however, we will omit here those details that do not pertain to the Krull dimension of the Stanley spaces; we refer the reader to the exposition in~\cite{EricksonHunzikerMOC2024}*{\S3}.)

\begin{lemma}[\cite{EricksonHunzikerMOC2024}*{Thm.~3.2}]
    \label{lemma:Stanley decomp via jellyfish}
    Assume the setting in Table~\ref{table:main info} where $G_\R = \U(p,q)$ or $\O^*(2n)$.
    
    If $k < s$, then the graded $S$-module $L_{\la(\sigma)}$ admits a Stanley decomposition
    \begin{equation}
        \label{Stanley decomp in lemma}
        L_{\la(\sigma)} \cong M_\sigma = \bigoplus_{(T, \mathbf{F}) \in \mathcal{J}_k(\sigma)}  \C[ z_{ij} : (i,j) \in \mathbf{F}] \left(\prod_{(i,j) \in {\rm cor}(\mathbf{F})} \hspace{-2ex} z_{ij} \right) \varphi_T,
    \end{equation}
    where
    {\normalfont
    \begin{itemize}
        \item $\mathcal{J}_k(\sigma)$ is the set of jellyfish to be defined below (Definition~\ref{def:jellyfish});
        \item each $T$ is a tableau in the set $\Tsigma$ from Definition~\ref{def:T(sigma)};
        \item each $\mathbf{F}$ is a (not necessarily maximal) union of $k$  nonintersecting lattice paths in $\mathbf{P}$ (see~\eqref{P}), and thus is an element of the complex $\Delta_k(\mathbf{P})$ defined above in~\eqref{Delta_k}; see also~\eqref{F} below;
        \item ${\rm cor}(\mathbf{F})$ is the set of ``corners'' of the lattice paths determining $\mathbf{F}$, defined in~\cite{EricksonHunzikerMOC2024}*{eqn.~(3.4)};
        the details are not pertinent to the result in the present paper.
        \item $\varphi_T$ is an $H(k)$-equivariant map canonically induced by $T$, and defined in~\cite{EricksonHunzikerMOC2024}*{eqn.~(2.12)};
        the details are not pertinent to the result in the present paper.
    \end{itemize}
    }
    
    If $k \geq s$, then $L_{\la(\sigma)}$ is a free $S$-module with Stanley decomposition
    \begin{equation}
        \label{jellyfish decomp k geq s}
        L_{\la(\sigma)} \cong M_\sigma =  \bigoplus_{T \in \Tsigma} S \, \varphi_T.
    \end{equation}
\end{lemma}

\subsection{Jellyfish: uniform overview}
\label{sub:jellyfish overview}

Sections~\ref{sub:jellyfish overview}--\ref{sub:jellyfish GL} will be a summary of our exposition in~\cite{EricksonHunzikerMOC2024}*{\S3}.
In this subsection, we give a uniform definition of the set $\mathcal{J}_k(\sigma)$ parametrizing the Stanley spaces in~\eqref{Stanley decomp in lemma}.
Then in the following two subsections, we will give the details specific to $G_\R = \O^*(2n)$ and $\U(p,q)$, respectively;
we choose this order because the combinatorial details in the $\O^*(2n)$ case are considerably more straightforward.

From now on, assume that $k < s$.
We recall the set $\mathbf{P}$ from~\eqref{P}.
We will define a region $\mathbf{A} \subseteq \mathbf{P}$ lying above a certain antidiagonal.
Then $\delta(\mathbf{A})$ or $\delta(\mathbf{P})$ will denote the set of points lying on a certain boundary of $\mathbf{A}$ or $\mathbf{P}$, respectively.
We define the set
\renewcommand{\arraystretch}{1}
\begin{equation}
    \label{F}
    \mathcal{F}_k \coloneqq \Big\{\big(\mathbf{A} \setminus \delta(\mathbf{A})\big) \sqcup \coprod_{t=1}^k \mathbf{L}_t  :  
    \text{each $\mathbf{L}_t$ is a lattice path from $\delta(\mathbf{A})$ to $\delta(\mathbf{P})$}
    \Big\},
\end{equation}
which contains the facet set $\widehat{\mathcal{F}}_k$ defined above in~\eqref{F hat general}.
Given $\mathbf{F} \in \mathcal{F}_k$, a \emph{corner} of $\mathbf{F}$ is an east-to-south turn in one of its lattice paths, such that the string of points to its southwest contains at least one point not in $\mathbf{F}$.
(Although in the $\U(p,q)$ case the decomposition of $\mathbf{F}$ into lattice paths $\mathbf{L}_t$ is in general not unique, nonetheless the corners of $\mathbf{F}$ are well defined.)
The symbol ${\rm cor}(\mathbf{F})$ in Lemma~\ref{lemma:Stanley decomp via jellyfish} denotes the set of corners of $\mathbf{F}$;
this will not, however, be relevant to the main result in this paper.

We define the collection
\begin{equation}
    \label{E}
    \mathcal{E}_k \coloneqq \left\{ 
    \{ \text{endpoint of $\mathbf{L}_t$} \}_{t=1}^k : \begin{array}{l} 
    \{\mathbf{L}_t\}_{t=1}^k \text{ is a family of nonintersecting}\\ \text{lattice paths from $\delta(\mathbf{A})$ to $\delta(\mathbf{P})$}
    \end{array}
    \right\},
\end{equation}
consisting of all possible sets $\mathbf{E} \subseteq \delta(\mathbf{P})$ of endpoints for the lattice paths appearing in~\eqref{F}.
We then partition $\mathcal{F}_k$ and $\Tsigma$ (from Definition~\ref{def:T(sigma)}) into subsets $\mathcal{F}_{\mathbf{E}}$ and $\mathcal{T}_{\mathbf{E}}(\sigma)$ as follows.

For $\mathbf{E} \in \mathcal{E}$, define
\begin{equation}
    \label{F_E}
    \mathcal{F}_{\mathbf{E}} \coloneqq \left\{ \mathbf{F} \in \mathcal{F}_k : \begin{array}{l} 
    \mathbf{F} = ( \mathbf{A} \setminus \delta(\mathbf{A})) \sqcup \coprod_{t=1}^k \mathbf{L}_t, \\
    \text{where each $\mathbf{L}_t$ is a lattice path from $\delta(\mathbf{A})$ to $\mathbf{E}$}
    \end{array}
    \right\}.
\end{equation}
For $G_\R = \O^*(2n)$ we will have $\mathcal{F}_k = \coprod_{\mathbf{E} \in \mathcal{E}} \mathcal{F}_{\mathbf{E}}$, but for $G_\R = \U(p,q)$ it is possible that $\mathcal{F}_{\mathbf{E}} = \mathcal{F}_{\mathbf{E}'}$.
Note that $\# \mathbf{F}$ is determined by the endpoints of its lattice paths;
hence, for each $\mathbf{E} \in \mathcal{E}$, every element of $\mathcal{F}_{\mathbf{E}}$ has the same cardinality.
In light of Lemma~\ref{lemma:Bdeg Stanley decomp}, we define
\begin{align}
    \label{F hat E hat}
    \begin{split}
    d_k & \coloneqq \max_{\mathbf{F} \in \mathcal{F}_k} \{ \# \mathbf{F} \},\\
    \widehat{\mathcal{F}}_k &\coloneqq \{ \mathbf{F} \in \mathcal{F}_k : \# \mathbf{F} = d_k \},\\
    \widehat{\mathcal{E}}_k &\coloneqq \{ \mathbf{E} \in \mathcal{E}_k : \text{if } \mathbf{F} \in \mathcal{F}_{\mathbf{E}} \text{ then } \mathbf{F} \in \widehat{\mathcal{F}}_k \}.
    \end{split}
\end{align}
Note that this definition of $\widehat{\mathcal{F}}_k$ is equivalent to the one given in~\eqref{F hat general}.

Given $\sigma \in \widehat{H}(k)$, and recalling the tableau set $\Tsigma$ from Definition~\ref{def:T(sigma)}, we will define a map
\begin{equation}
    \label{end map}
    {\rm end} : \Tsigma \longrightarrow \mathcal{E}_k,
\end{equation}
so that (very roughly speaking) ${\rm end}(T)$ is the set of endpoints most closely aligning with the initial column of $T$.
From now on, we write $T_{t,1}$ to denote the entry of $T$ in row $t$ and column $1$.
For each $\mathbf{E} \in \mathcal{E}_k$, we define the fiber
\begin{equation}
    \label{T_E}
    \mathcal{T}_{\mathbf{E}}(\sigma) \coloneqq \{ T \in \Tsigma : {\rm end}(T) = \mathbf{E} \}.
\end{equation}
Clearly $\Tsigma = \coprod_{\mathbf{E} \in \mathcal{E}} T_{\mathbf{E}}(\sigma)$, although we may have $\mathcal{T}_{\mathbf{E}}(\sigma) = \varnothing$ for certain $\mathbf{E} \in \mathcal{E}_k$.

\begin{dfn}
    \label{def:jellyfish}
    Let $G_\R = \U(p,q)$ or $\O^*(2n)$, with $k<s$ and $\sigma \in \widehat{H}(k)$, as in Table~\ref{table:main info}.
    A \emph{jellyfish of shape~$\sigma$} is an element of the set
    \begin{align*}
        \mathcal{J}_k(\sigma) &\coloneqq \coprod_{\mathbf{E} \in \mathcal{E}_k} \mathcal{T}_{\mathbf{E}}(\sigma) \times \mathcal{F}_{\mathbf{E}}, 
    \end{align*}
    with $\mathcal{T}_{\mathbf{E}}(\sigma)$ and $\mathcal{F}_{\mathbf{E}}$ as in~\eqref{T_E} and~\eqref{F_E}, and with case-by-case details to be described below.
    Recalling the set $\widehat{\mathcal{F}}_k$ from~\eqref{F hat E hat}, we say a jellyfish is \emph{maximal} if it belongs to the subset
    \[
    \widehat{\mathcal{J}}_k(\sigma) \coloneqq \Big\{ (T, \mathbf{F}) \in \mathcal{J}_k(\sigma) : \mathbf{F} \in \widehat{\mathcal{F}}_k \Big\}.
    \]
    
\end{dfn}

\begin{lemma}
    \label{lemma:Deg equals J hat}
    Let $G_\R = \U(p,q)$ or $\O^*(2n)$, with $k < s$ and $\sigma \in \widehat{H}(k)$, as in Table~\ref{table:main info}.
    We have
    \[
    \Deg L_{\la(\sigma)} = \# \widehat{\mathcal{J}}_k(\sigma).
    \]
\end{lemma}

\begin{proof}
Consider the Stanley decomposition of $L_{\la(\sigma)}$ in~\eqref{Stanley decomp in lemma}, which is the $k<s$ case of Lemma~\ref{lemma:Stanley decomp via jellyfish}.
Each $(T, \mathbf{F}) \in \mathcal{J}_k(\sigma)$ corresponds to the Stanley space $S_{(T, \mathbf{F})} \eta_{(T, \mathbf{F})}$, where $S_{(T, \mathbf{F})} = \C[ z_{ij} : (i,j) \in \mathbf{F}]$ has Krull dimension equal to $\# \mathbf{F}$.
Now we can apply Lemma~\ref{lemma:Bdeg Stanley decomp}, where the role of $M$ is played by $L_{\la(\sigma)}$, the role of $d$ is played by $d_k$ from~\eqref{F hat E hat}, the role of $\mathcal{I}$ is played by $\mathcal{J}_k(\sigma)$, and the role of $\widehat{\mathcal{I}}$ is played by $\widehat{\mathcal{J}}_k(\sigma)$ from Definition~\ref{def:jellyfish}.
In this special case, Lemma~\ref{lemma:Bdeg Stanley decomp} yields
\[
    e\! \left(L_{\la(\sigma)} \right) = \# \widehat{\mathcal{J}}_k(\sigma).
\]
The result follows from~\eqref{Deg = e for L}.
\end{proof}

\subsection{Details for \texorpdfstring{$\O^*(2n)$}{O*(2n)}}

Let $G_\R = \O^*(2n)$.
Then $\mathbf{P} = \{(i,j) : 1 \leq i < j \leq n \}$ is depicted as an $(n-1)$-staircase, with $(1,2)$ in the northwest corner, and $(n-1,n)$ in the southeast corner.
Let
\begin{align*}
    \mathbf{A} &\coloneqq \{ (i,j) \in \mathbf{P} : i + j \leq 2(k + 1) \},\\
    \delta(\mathbf{A}) & \coloneqq \{ (i,j) \in \mathbf{A} : i + j = 2(k+1) \text{ or } j = n \},\\
    \delta(\mathbf{P}) &\coloneqq \{ (i,j) \in \mathbf{P} : j=n \}.
\end{align*}
Note that $\delta$ denotes the eastern boundary.
Explicitly, we have 
\begin{equation}
    \label{delta A Sp}
    \delta(\mathbf{A}) = \{(t, b_t)\}_{t=1}^k, \quad \text{where} \quad b_t
= \begin{cases}
    n & \text{if }t < 2(k+1) - n,\\
    2(k+1)-t & \text{otherwise}.
\end{cases}
\end{equation}
Referring to~\eqref{F}, we observe that each $\mathbf{F} \in \mathcal{F}_k$ admits a unique decomposition $\mathbf{F} = (\mathbf{A} \setminus \delta(\mathbf{A})) \sqcup \coprod_{t=1}^k \mathbf{L}_t$, where each $\mathbf{L}_t$ is a southeast lattice path starting at $(t, b_t)$ and ending in $\delta(\mathbf{P})$.
In~\eqref{Sp F} below, where $n=11$, we show examples of an element $\mathbf{F} \in \mathcal{F}_k$ for two different values of $k$.
Following~\eqref{F}, we shade the points in $\mathbf{F}$, and we draw squares to indicate the elements of ${\rm cor}(\mathbf{F})$:
\begin{equation} 
    \label{Sp F} 
    \tikzstyle{dot}=[circle,fill=black, minimum size = 4.5pt, inner sep=0pt]
    \tikzstyle{corner}=[rectangle,draw=black,thin, minimum size = 8pt, inner sep=2pt]
\begin{tikzpicture}          [scale=.3,baseline=(current bounding box.center),every node/.style={scale=.7}]
\draw [white,fill=lightgray] (1.5,11.5) -- ++ (0,-1) -- ++ (1,0) -- ++ (0,-1) -- ++ (1,0) -- ++ (0,-1) -- ++ (1,0) -- ++ (0,-1) -- ++ (1,0) -- ++ (0,1)  -- ++ (1,0) -- ++ (0,1)  -- ++ (1,0) -- ++ (0,1) -- ++ (1,0) -- ++ (0,1) -- cycle;
\draw[line width=3pt, lightgray] (9,11) -- ++(1,0) node [corner] {} -- ++(0,-1) -- ++(1,0) (8,10) -- ++(0,-1) -- ++(1,0) -- ++(0,-1) -- ++(2,0) (7,9) -- ++(0,-1) -- ++(1,0) -- ++(0,-1) -- ++(2,0) node [corner] {} -- ++(0,-1) -- ++(1,0) node [corner] {} -- ++(0,-1) (6,8) -- ++(0,-1) -- ++(1,0) -- ++(0,-1) -- ++(2,0) node [corner] {} -- ++(0,-2) -- ++(2,0);
\foreach \x in {2,...,11}{\foreach \y in {\x,...,11}{\node [dot] at (13-\x,\y) {};}}
\node[scale=1.3] at (6,2) {$k=4$};
\draw [densely dotted] 
(1.5, 11.5) -- ++(10,0) -- ++(0,-10) -- ++(-1,0) -- ++(0,1) -- ++(-1,0) -- ++(0,1) -- ++(-1,0) -- ++(0,1) -- ++(-1,0) -- ++(0,1) -- ++(-1,0) -- ++(0,1) -- ++(-1,0) -- ++(0,1) -- ++(-1,0) -- ++(0,1) -- ++(-1,0) -- ++(0,1) -- ++(-1,0) -- ++(0,1) -- ++(-1,0) -- ++(0,1);
\end{tikzpicture}
\qquad\qquad\qquad\qquad
\begin{tikzpicture}          [scale=.3,baseline=(current bounding box.center),every node/.style={scale=.7}]
\draw[white, fill=lightgray] (1.5,11.5) --++ (9,0) -- ++(0,-5) -- ++(-1,0) -- ++(0,-1) -- ++(-1,0) -- ++(0,-1) -- ++(-1,0) -- ++(0,1) -- ++(-1,0) -- ++(0,1) -- ++(-1,0) -- ++(0,1) -- ++(-1,0) -- ++(0,1) -- ++(-1,0) -- ++(0,1) -- ++(-1,0) -- ++(0,1) -- ++(-1,0) -- cycle;
\draw[line width=6pt, lightgray, line cap=round] (11,11) -- ++(0,0)
(11,10) -- ++(0,0)
(11,9) -- ++(0,0)
(11,8) -- ++(0,0);
\draw[line width=3pt, lightgray] (11,7) -- ++(0,-1)
(10,6) -- ++(0,-1) -- ++(1,0)
(9,5) -- ++(0,-1) -- ++(2,0) node [corner] {} -- ++(0,-2);
\foreach \x in {2,...,11}{\foreach \y in {\x,...,11}{\node [dot] at (13-\x,\y) {};}}
\node[scale=1.3] at (6,2) {$k=7$};
\draw [densely dotted] 
(1.5, 11.5) -- ++(10,0) -- ++(0,-10) -- ++(-1,0) -- ++(0,1) -- ++(-1,0) -- ++(0,1) -- ++(-1,0) -- ++(0,1) -- ++(-1,0) -- ++(0,1) -- ++(-1,0) -- ++(0,1) -- ++(-1,0) -- ++(0,1) -- ++(-1,0) -- ++(0,1) -- ++(-1,0) -- ++(0,1) -- ++(-1,0) -- ++(0,1) -- ++(-1,0) -- ++(0,1);
\end{tikzpicture}
\end{equation}

    \noindent Let $\mathbf{E} = \{(i_t, n)\}_{t=1}^k \in \mathcal{E}_k$ where the $i_t$'s are listed in increasing order.
    By~\eqref{F_E}, if $\mathbf{F} \in \mathcal{F}_{\mathbf{E}}$, then each $\mathbf{L}_t$ is a lattice path from $(t, b_t)$ to $(i_t, n)$.
    Hence
    \begin{equation}
        \label{size F Sp proof}
        \#\mathbf{F} = \#(\mathbf{A} \setminus \delta(\mathbf{A})) + \sum_{t=1}^k \# \mathbf{L}_t, \quad \text{where} \quad 
        \#\mathbf{L}_t = (i_t - t) + (n-b_t) + 1.
    \end{equation}
    The $b_t$'s in~\eqref{delta A Sp} are independent of $\mathbf{F}$, and thus $\# \mathbf{F}$ depends only on the sum $\sum_{t=1}^k i_t$.
    Crucially, due to the nonintersecting property of the lattice paths, the maximum attainable value $\hat{\imath}_t$ for $i_t$ is given by
\begin{equation}
    \label{i hat t}
   \hat{\imath}_t = \begin{cases}
    t & \text{if } t < 2(k+1) - n,\\
    n + 2(t-k) - 1 & \text{otherwise}.
\end{cases}
\end{equation}
Note that the topmost case in~\eqref{i hat t}
corresponds to those $t$'s (if any) for which $\mathbf{L}_t = \{(t, b_t)\}$ is necessarily a singleton;
for example, in the right-hand example in~\eqref{Sp F}, this occurs for $t < 2(k+1) - n = 5$, and indeed the four topmost lattice paths are singletons.
Putting $\widehat{\mathbf{E}} \coloneqq \{(\hat{\imath}_t, n)\}_{t=1}^k$, it follows from the above discussion that $\mathbf{F}$ attains the maximum possible cardinality $d_k$ if and only if $\mathbf{F} \in \mathcal{F}_{\widehat{\mathbf{E}}}$;
that is, recalling the sets $\widehat{\mathcal{E}}_k$ and $\widehat{\mathcal{F}}_k$ defined in~\eqref{F hat E hat},
\begin{equation}
    \label{F hat E hat Sp}
    \widehat{\mathcal{E}}_k = \{ \widehat{\mathbf{E}} \}, \text{ and therefore } \widehat{\mathcal{F}}_k = \mathcal{F}_{\widehat{\mathbf{E}}}.
\end{equation}

We now define the ``end'' map~\eqref{end map} as follows, for $T \in \Tsigma = \SSYT(\sigma, n)$:
\begin{equation}
    \label{end T details Sp}
    {\rm end}(T) \coloneqq \{ (i_t, n)\}_{t=1}^k, \quad \text{where } i_t = \min\{ T_{t,1}, \: \hat{\imath}_t \}.
\end{equation}
Since $T_{t,1}$ is undefined for $t > \ell(\sigma)$, in this range the equation~\eqref{end T details Sp} reduces to $i_t = \hat{\imath}_t$.
In the running examples below, for a given tableau $T$ we indicate ${\rm end}(T)$ by highlighting its elements in red inside the vertical strip $\delta(\mathbf{P})$:
\begin{equation*} 
    \label{Sp end} 
    \tikzstyle{dot}=[circle,fill=black, minimum size = 4.5pt, inner sep=0pt]
    \tikzstyle{redglow}=[circle,fill=red, minimum size = 7.5pt, inner sep=0pt]
\begin{tikzpicture}          [scale=.3,baseline=(current bounding box.center),every node/.style={scale=.7}]
\node [redglow] at (0,2) {};
\node [redglow] (E3) at (0,4) {};
\node [redglow] (E2) at (0,6) {};
\node [redglow] (E1) at (0,10) {};
\foreach \y in {2,...,11}{\node [dot] at (0,\y) {};}
\draw[densely dotted, thick, lightgray]
(-.5,11.5) -- ++ (1,0) -- ++(0,-10) -- ++(-1,0) -- cycle;
\node [scale=1.3, left] at (-4,6) {$T = \ytableaushort{2355,89,9}$};
\node[left] at (-.5,8) {$\hat{\imath}_1 \rightarrow$} ;
\node[left] at (-.5,6) {$\hat{\imath}_2 \rightarrow$} ;
\node[left] at (-.5,4) {$\hat{\imath}_3 \rightarrow$} ;
\node[left] at (-.5,2) {$\hat{\imath}_4 \rightarrow$} ;
\node[right] (T1) at (2,10) {$\ytableaushort{2} \leftarrow T_{1,1}$};
\node[right] (T2) at (2,4) {$\ytableaushort{8} \leftarrow T_{2,1}$};
\node[right] (T3) at (2,3) {$\ytableaushort{9} \leftarrow T_{3,1}$};
\draw [ultra thick, lightgray] (T1.west) -- (E1);
\draw [ultra thick, lightgray] (T2.west) -- (E2);
\draw [ultra thick, lightgray] (T3.west) -- (E3);
\node[scale=1.3] at (0,0) {$k=4$};
\node[scale=1.3, red] at (0,13) {${\rm end}(T)$};
\end{tikzpicture}
\qquad\qquad\qquad\qquad\begin{tikzpicture}          [scale=.3,baseline=(current bounding box.center),every node/.style={scale=.7}]
\node [redglow] (E1) at (0,11) {};
\node [redglow] (E2) at (0,10) {};
\node [redglow] (E3) at (0,9) {};
\node [redglow] (E4) at (0,8) {};
\node [redglow] (E5) at (0,6) {};
\node [redglow] (E6) at (0,4) {};
\node [redglow] (E7) at (0,3) {};
\foreach \y in {2,...,11}{\node [dot] at (0,\y) {};}
\draw[densely dotted, thick, lightgray]
(-.5,11.5) -- ++ (1,0) -- ++(0,-10) -- ++(-1,0) -- cycle;
\node [scale=1.3, left] at (-4,6) {$T = \ytableaushort{235,468,57,68,79,8,9}$};
\node[left] at (-.5,11) {$\hat{\imath}_1 \rightarrow$} ;
\node[left] at (-.5,10) {$\hat{\imath}_2 \rightarrow$} ;
\node[left] at (-.5,9) {$\hat{\imath}_3 \rightarrow$} ;
\node[left] at (-.5,8) {$\hat{\imath}_4 \rightarrow$} ;
\node[left] at (-.5,6) {$\hat{\imath}_5 \rightarrow$} ;
\node[left] at (-.5,4) {$\hat{\imath}_6 \rightarrow$} ;
\node[left] at (-.5,2) {$\hat{\imath}_7 \rightarrow$} ;
\node[right] (T1) at (2,10) {$\ytableaushort{2} \leftarrow T_{1,1}$};
\node[right] (T2) at (2,8) {$\ytableaushort{4} \leftarrow T_{2,1}$};
\node[right] (T3) at (2,7) {$\ytableaushort{5} \leftarrow T_{3,1}$};
\node[right] (T4) at (2,6) {$\ytableaushort{6} \leftarrow T_{4,1}$};
\node[right] (T5) at (2,5) {$\ytableaushort{7} \leftarrow T_{5,1}$};
\node[right] (T6) at (2,4) {$\ytableaushort{8} \leftarrow T_{6,1}$};
\node[right] (T7) at (2,3) {$\ytableaushort{9} \leftarrow T_{7,1}$};
\draw [ultra thick, lightgray] (T1.west) -- (E1);
\draw [ultra thick, lightgray] (T2.west) -- (E2);
\draw [ultra thick, lightgray] (T3.west) -- (E3);
\draw [ultra thick, lightgray] (T4.west) -- (E4);
\draw [ultra thick, lightgray] (T5.west) -- (E5);
\draw [ultra thick, lightgray] (T6.west) -- (E6);
\draw [ultra thick, lightgray] (T7.west) -- (E7);
\node[scale=1.3] at (0,0) {$k=7$};
\node[scale=1.3, red] at (0,13) {${\rm end}(T)$};
\end{tikzpicture}
\end{equation*}

\noindent 


\subsection{Details for \texorpdfstring{$\U(p,q)$}{U(p,q)}}
\label{sub:jellyfish GL}

Let $G_\R = \U(p,q)$.
Then $\mathbf{P} = \{ (i,j) : 1 \leq i \leq p, \: 1 \leq j \leq q \}$ is a $p \times q$ grid, with $(1,1)$ in the northwest corner, and $(p,q)$ in the southwest corner.
Let
\begin{align*}
    \mathbf{A} &\coloneqq \{(i,j) \in \mathbf{P} : i + j \leq k + 1\},\\
    \delta(\mathbf{A}) &\coloneqq \{(i,j) \in \mathbf{A} : i + j = k + 1 \text{ or } i = p \text{ or } j = q \},\\
    \delta(\mathbf{P}) & \coloneqq \{(i,j) \in \mathbf{P} : i = p \text{ or } j = q \}.
\end{align*}
Note that $\delta$ denotes the southeastern boundary.
Given $\sigma \in \widehat{H}(k)$, it will be useful to define
\begin{align}
    \label{k+ k-}
    \begin{split}
        k^+ & \coloneqq \max\{ \ell(\sigma^+), \: k - p \},\\
        k^- & \coloneqq k - k^+,
    \end{split}
\end{align}
so that $k = k^+ + k^-$.
The purpose of~\eqref{k+ k-} is to divide $\delta(\mathbf{A})$ into two sides:
due to the upcoming definition of ${\rm end}(T)$ in~\eqref{end T GL details}, it will turn out that $(T, \mathbf{F}) \in \mathcal{J}_k(\sigma)$ only if the $k^+$ many southernmost lattice paths in $\mathbf{F}$ have endpoints along the southern edge of $\delta(\mathbf{P})$, and the remaining $k^-$ many easternmost lattice paths have endpoints along the eastern edge of $\delta(\mathbf{P})$.
In particular, if $k-p$ is positive, then it equals the number of starting points along the southern edge of $\delta(\mathbf{A}) \cap \delta(\mathbf{P})$ that are necessarily their own endpoints.
Explicitly, then, we have $\delta(\mathbf{A}) = \{(a_u, u) \}_{u=1}^{k^+} \sqcup \{(t, b_t)\}_{t=1}^{k^-}$, where
\begin{equation}
    \label{a's b's GL}
    a_u = 
    \begin{cases}
        p & \text{if }u \leq k-p,\\
        k-u+1 & \text{otherwise},
    \end{cases} \qquad\text{and}\qquad
    b_t = 
    \begin{cases}
        q & \text{if } t \leq k-q,\\
        k-t+1 & \text{otherwise}.
    \end{cases}
\end{equation}
By~\eqref{F}, for each $\mathbf{F} \in \mathcal{F}_k$ we have $\mathbf{F} = (\mathbf{A} \setminus \delta(\mathbf{A})) \sqcup \mathbf{L}$, for a (not necessarily unique) disjoint union $\mathbf{L} = \coprod_{u=1}^{k^+} \mathbf{L}^+_u \sqcup \coprod_{t=1}^{k^-} \mathbf{L}^-_t$, where each $\mathbf{L}^+_u$ is a southeast lattice path from $(a_u, u)$ to $\delta(\mathbf{P})$, and $\mathbf{L}^-_t$ is a southeast lattice path from $(t,b_t)$ to $\delta(\mathbf{P})$.
Below, for $(p,q) = (7,10)$, we show examples of an element $\mathbf{F} \in \mathcal{F}_k$ for three different values of $k$:
\begin{equation*} 
    \label{GL F} 
    \tikzstyle{dot}=[circle,fill=black, minimum size = 4.5pt, inner sep=0pt]
    \tikzstyle{corner}=[rectangle,draw=black,thin, minimum size = 8pt, inner sep=2pt]
\begin{tikzpicture}          [scale=.3,baseline=(current bounding box.center),every node/.style={scale=.7}]
\draw[white, fill=lightgray]
(.5,7.5) -- ++(3,0) --++(0,-1) -- ++(-1,0) --++(0,-1) -- ++(-1,0) --++(0,-1) -- ++(-1,0) -- cycle;
\draw[line width=3pt, lightgray, line cap=round] 
(1,4) --++(0,-1) --++(1,0) node [corner] {} -- ++(0,-2) -- ++(2,0)
(2,5) -- ++(1,0) node [corner] {} -- ++(0,-3) -- ++(2,0) -- ++(0,-1) -- ++(2,0)
(3,6) -- ++(2,0) node [corner] {} -- ++(0,-1) -- ++(2,0) node [corner] {} -- ++(0,-2) -- ++(1,0) node [corner] {} -- ++(0,-1) -- ++(1,0) node [corner] {} -- ++(0,-1)
(4,7) -- ++(4,0) node [corner] {} -- ++(0,-1) -- ++(1,0) node[corner] {} --++(0,-1) -- ++ (1,0);
\foreach \x in {1,...,10}{\foreach \y in {1,...,7}{\node [dot] at (\x,\y) {};}}

\node[below,scale=1.3] at (5.5,0) {$k=4$};
\draw [densely dotted] (.5,7.5) rectangle (10.5,.5);
\end{tikzpicture}
\qquad\qquad\qquad
\begin{tikzpicture}          [scale=.3,baseline=(current bounding box.center),every node/.style={scale=.7}]
\draw[white, fill=lightgray]
(.5,7.5) -- ++(7,0) --++(0,-1) -- ++(-1,0) --++(0,-1) -- ++(-1,0)--++(0,-1) -- ++(-1,0)--++(0,-1) -- ++(-1,0)--++(0,-1) -- ++(-1,0)--++(0,-1) -- ++(-2,0) -- cycle;
\draw[line width=6pt, lightgray, line cap=round]
(1,1) --++(0,0);
\draw[line width=3pt, lightgray, line cap=round] 
(2,1) --++(2,0)
(3,2) -- ++(2,0) -- ++(0,-1) -- ++(2,0)
(4,3) -- ++(2,0) -- ++(0,-1) -- ++(3,0) node [corner] {} -- ++(0,-1) -- ++(1,0)
(5,4) -- ++(3,0) node [corner] {} -- ++(0,-1) -- ++(2,0) node[corner] {} --++(0,-1)
(6,5) -- ++(4,0)
(7,6) -- ++(3,0)
(8,7) -- ++(2,0);
\foreach \x in {1,...,10}{\foreach \y in {1,...,7}{\node [dot] at (\x,\y) {};}}

\node[below,scale=1.3] at (5.5,0) {$k=8$};
\draw [densely dotted] (.5,7.5) rectangle (10.5,.5);
\end{tikzpicture}
\qquad\qquad\qquad
\begin{tikzpicture}          [scale=.3,baseline=(current bounding box.center),every node/.style={scale=.7}]
\draw[white, fill=lightgray]
(.5,7.5) -- ++(9,0) --++(0,-3) -- ++(-1,0) --++(0,-1) -- ++(-1,0)--++(0,-1) -- ++(-1,0)--++(0,-1) -- ++(-6,0) -- cycle;
\draw[line width=6pt, lightgray, line cap=round]
(1,1) --++(0,0)
(2,1) --++(0,0)
(3,1) --++(0,0)
(4,1) --++(0,0)
(5,1) --++(0,0)
(10,7) --++(0,0)
(10,6) --++(0,0)
(10,5) --++(0,0);
\draw[line width=3pt, lightgray, line cap=round] 
(6,1) -- ++(1,0)
(7,2) -- ++(1,0) --++ (0,-1)
(8,3) --++(2,0) node [corner] {} --++(0,-1)
(9,4) --++(1,0);
\foreach \x in {1,...,10}{\foreach \y in {1,...,7}{\node [dot] at (\x,\y) {};}}
\node[below,scale=1.3] at (5.5,0) {$k=12$};
\draw [densely dotted] (.5,7.5) rectangle (10.5,.5);
\end{tikzpicture}
\end{equation*}

    \noindent Let $\mathbf{E} = \{(p, j_u)\}_{u=1}^{k^+} \sqcup \{(i_t, q)\}_{t=1}^k \in \mathcal{E}_k$ where the $j_u$'s and $i_t$'s are listed in increasing order.
    By~\eqref{F_E}, if $\mathbf{F} \in \mathcal{F}_{\mathbf{E}}$, then each $\mathbf{L}^+_u$ is a lattice path from $(a_u, u)$ to $(p, j_u)$, and each $\mathbf{L}^-_t$ is a lattice path from $(t, b_t)$ to $(i_t, q)$.
    By same argument as in~\eqref{size F Sp proof}, $\# \mathbf{F}$ depends only on the sum $\sum_u j_u + \sum_t i_t$.
    Crucially, due to the nonintersecting condition on the lattice paths $\mathbf{L}^+_u$ and $\mathbf{L}^-_t$, each $\mathbf{E} \in \mathcal{E}_k$ has the following property: 
    the $\ell \times \ell$ square in the lower-right corner of $\mathbf{P}$ contains at most $\ell$  points in $\mathbf{E}$, for all $1 \leq \ell \leq \min\{k,r\}$.
Thus we have the tight inequalities
    \begin{equation} 
        \label{endpoint spacing GL}
        \#\{ u : j_u > q - \ell \} + \#\{ t: i_t > p - \ell \} \leq \ell \quad \text{for all } 1 \leq \ell \leq \min\{k,r\}.
    \end{equation}
    It follows that $\mathbf{F}$ attains the maximum possible cardinality $d_k$ if and only if the inequalities~\eqref{endpoint spacing GL} are all tight at $\mathbf{E}$, that is, $\#\{ u : j_u > q - \ell \} + \#\{ t: i_t > p - \ell \} = \ell$ for all $1 \leq \ell \leq \min\{k,r\}$.
    Equivalently, since $\# \mathbf{E} = k$, we have that $\#\mathbf{F} = d_k$ if and only if
    \begin{equation}
        \label{E max condition GL}
        \#\{u : j_u \leq q - \ell \} + \#\{t : i_t \leq p - \ell\} = k - \ell \quad \text{for all $1 \leq \ell \leq \min\{k,r\}$},
    \end{equation}
    that is, if and only if $\mathbf{F}$ contains the isosceles right triangle in the southeastern corner of $\mathbf{P}$ with side length $\min\{k,r\}$.
    Thus by~\eqref{F hat E hat},
    \begin{equation}
        \label{F hat E hat GL}
        \widehat{\mathcal{F}}_k = \mathcal{F}_{\mathbf{E}} \text{ if and only if  $\mathbf{E} \in \widehat{\mathcal{E}}_k = \{ \mathbf{E} \in \mathcal{E}_k : \mathbf{E} \text{ satisfies } \eqref{E max condition GL} \}$}.
    \end{equation}
    
    Observe from~\eqref{E max condition GL} that upon fixing the $j_u$'s, the maximum attainable value $\hat{\imath}_t$ of $i_t$ is given by $\hat{\imath}_t = t$ if $t < k-q$,
    while if $t \geq k-q$, then the $\hat{\imath}_t$'s are the largest elements $\hat{\imath}_{\min\{1, \: k-q\}} < \cdots < \hat{\imath}_{k^-}$ of the set
    \begin{equation}
        \label{i hat GL}
        \Big\{ i \in [p] : q-p+i \notin \{j_u \}_{u=1}^{k^+} \Big\}.
    \end{equation}
    In other words, for $\min\{1, k-q\} \leq t \leq k^-$, the $\hat{\imath}_t$'s are the distinct elements chosen from $[p]$ to be as large as possible without violating~\eqref{endpoint spacing GL}; see the example~\eqref{GL end} below.

We define the ``end'' map~\eqref{end map} as follows, for $T = (T^+, T^-) \in \Tsigma = \SSYT(\sigma^+, q) \times \SSYT(\sigma^-, p)$:
\[
    {\rm end}(T) \coloneqq \{ (p, j_u) \}_{u=1}^{k^+} \sqcup \{ (i_t, q)\}_{t=1}^{k^-},
\]
where
\begin{equation}
    \label{end T GL details}
    j_u = \begin{cases}
        u & \text{if } u \leq k-p,\\
        T^+_{u,1} & \text{otherwise},
    \end{cases} \qquad \text{and} \qquad 
    i_t = \min\{ T^-_{t,1}, \: \hat{\imath}_t \}.
\end{equation}
Since $T^-_{t,1}$ is undefined for $t > \ell(\sigma^-)$, in this range the right-hand side of~\eqref{end T GL details} reduces to $i_t = \hat{\imath}_t$.
As an example below, we indicate ${\rm end}(T)$ by highlighting its elements in red inside $\delta(\mathbf{P})$:

\begin{equation} 
    \label{GL end}
    \tikzstyle{dot}=[circle,fill=black, minimum size = 4.5pt, inner sep=0pt]
    \tikzstyle{redglow}=[circle,fill=red, minimum size = 7.5pt, inner sep=0pt]
\begin{tikzpicture}          [scale=.3,baseline=(current bounding box.center),every node/.style={scale=.7}]
\node [redglow] (D1) at (1,1) {};
\node [redglow] (D2) at (4,1) {};
\node [redglow] (D3) at (7,1) {};
\node [redglow] (D4) at (9,1) {};
\node [redglow] (D5) at (10,1) {};
\node [redglow] (E1) at (10,7) {};
\node [redglow] (E2) at (10,5) {};
\node [redglow] (E3) at (10,3) {};
\foreach \x in {1,...,10}{\node [dot] at (\x,1) {};}
\foreach \y in {2,...,7}{\node [dot] at (10,\y) {};}
\draw[densely dotted, thick, lightgray]
(9.5,7.5) -- ++ (1,0) -- ++(0,-7) -- ++(-10,0) -- ++(0,1) -- ++(9,0) -- cycle;
\draw [densely dashed] (8,1) -- ++(0,1.5)
(6,1) -- ++(0,4) -- ++(1.25,0)
(5,1) -- ++(0,5) -- ++(2.25,0);
\node [scale=1.3, left] at (0,4) {$\begin{aligned}T &= \left(\:\ytableaushort{348,469,78,99,{10}}, \: \ytableaushort{125,46,7,\none,\none}\:\right) \\[2ex] k &= 8 \\ k^+ &= 5 \\ k^- &= 3 \end{aligned}$};
\node[left] at (9.5,6) {$\hat{\imath}_1 \rightarrow$} ;
\node[left] at (9.5,5) {$\hat{\imath}_2 \rightarrow$} ;
\node[left] at (9.5,3) {$\hat{\imath}_3 \rightarrow$} ;
\node[right] (T1) at (12,7) {$\ytableaushort{1} \leftarrow T^-_{1,1}$};
\node[right] (T2) at (12,4) {$\ytableaushort{4} \leftarrow T^-_{2,1}$};
\node[right] (T3) at (12,1) {$\ytableaushort{7} \leftarrow T^-_{3,1}$};
\draw [ultra thick, lightgray] (T1.west) -- (E1);
\draw [ultra thick, lightgray] (T2.west) -- (E2);
\draw [ultra thick, lightgray] (T3.west) -- (E3);
\node[below,align=center] (U1) at (3,-1) {$\ytableaushort{3}$\\$\uparrow$\\$T^+_{1,1}$};
\node[below] (U2) at (4,-1) {$\ytableaushort{4}$};
\node[below] (U3) at (7,-1) {$\ytableaushort{7}$};
\node[below] (U4) at (9,-1) {$\ytableaushort{9}$};
\node[below,align=center] (U5) at (10,-1) {$\ytableaushort{{10}}$\\$\uparrow$\\$T^+_{5,1}$};
\draw [loosely dotted, thick] (4,-4) -- (9,-4);
\draw [ultra thick, lightgray] (U1.north) -- (D1);
\draw [ultra thick, lightgray] (U2.north) -- (D2);
\draw [ultra thick, lightgray] (U3.north) -- (D3);
\draw [ultra thick, lightgray] (U4.north) -- (D4);
\draw [ultra thick, lightgray] (U5.north) -- (D5);
\node[scale=1.3, red] at (10,9) {${\rm end}(T)$};
\end{tikzpicture}
\end{equation}

\noindent The dashed lines in~\eqref{GL end} are meant to illustrate the definition of the $\hat{\imath}_t$'s in~\eqref{i hat GL}.
Specifically, on the southern edge of $\delta(\mathbf{P})$, take the easternmost points that are \emph{not} in the set $\{(p, j_u) \}_{u=1}^{k^+}$; then the $\hat{\imath}_t$'s are their corresponding points along the eastern edge of $\delta(\mathbf{P})$.

\section{Main result}
\label{sec:main result}

We recall the sets $\Qsigma$ and $\PP$ from Definitions~\ref{def:T(sigma)} and~\ref{def:R}, along with the facet set $\widehat{\mathcal{F}}_k$ from~\eqref{F hat general}, and the map $\Theta : \PP \longrightarrow \widehat{\mathcal{F}}_k$ from Definition~\ref{def:Theta}.
We also recall the set $\widehat{\mathcal{J}}_k(\sigma)$ of maximal jellyfish in Definition~\ref{def:jellyfish}.

\begin{lemma}
    \label{lemma:Theta tilde}
    Let $G_\R = \U(p,q)$ or $\O^*(2n)$, where $k < s$ (see Table~\ref{table:main info}).
    We have a bijection
    \[
        {\rm id} \times \Theta : \Qsigma \times \PP \longrightarrow \widehat{\mathcal{J}}_k(\sigma).
    \]
\end{lemma}


               



\begin{proof}

By Lemma~\ref{lemma:Theta}, ${\rm id} \times \Theta$ is a bijection from $\Qsigma \times \PP$ to $\Qsigma \times \widehat{\mathcal{F}}_k$.
Hence it suffices to prove that
\begin{equation}
    \label{WTS JQF}
    \widehat{\mathcal{J}}_k(\sigma) = \Qsigma \times \widehat{\mathcal{F}}_k.
\end{equation}
Starting from Definition~\ref{def:jellyfish}, we have
\begin{equation}
    \label{JTF}
    \widehat{\mathcal{J}}_k(\sigma) \coloneqq \{ (T, \mathbf{F}) \in \mathcal{J}_k(\sigma) : \mathbf{F} \in \widehat{\mathcal{F}}_k \} = \coprod_{\mathbf{E} \in \widehat{\mathcal{E}}_k} \mathcal{T}_{\mathbf{E}}(\sigma) \times \mathcal{F}_{\mathbf{E}},
\end{equation}  
where the second equality follows from ~\eqref{F hat E hat} and Definition~\ref{def:jellyfish}.

    First let $G_\R = \O^*(2n)$.
    By~\eqref{F hat E hat Sp}, we can rewrite~\eqref{JTF} as
    \begin{equation}
        \label{JTF rewrite Sp}
        \widehat{\mathcal{J}}_k(\sigma) = \mathcal{T}_{\widehat{\mathbf{E}}} \times \mathcal{F}_{\widehat{\mathbf{E}}} = \mathcal{T}_{\widehat{\mathbf{E}}} \times \widehat{\mathcal{F}}_k.
    \end{equation}
    By~\eqref{T_E}, $T \in \mathcal{T}_{\widehat{\mathbf{E}}}$ if and only if ${\rm end}(T) = \widehat{\mathbf{E}} = \{(\hat{\imath}_t, n)\}_{t=1}^k$, with the $\hat{\imath}_t$'s defined in~\eqref{i hat t}.
    In turn, by~\eqref{end T details Sp}, ${\rm end}(T) = \widehat{\mathbf{E}}$ if and only if $T_{t,1} \geq \hat{\imath}_t$ for all $1 \leq t \leq \ell(\sigma)$.
    This is automatic for $t < 2(k+1)-n$, that is for $t \leq 2k-n+1$, so we need only check that $T_{t,1} \geq \hat{\imath}_t$ for all $\max\{0, 2k-n+1\} < j \leq \ell(\sigma)$.
    But by comparing this and~\eqref{i hat t} with Lemma~\ref{lemma:explicit criteria}(3), where $j$ plays the role of $t$, we conclude that $T \in \mathcal{T}_{\widehat{\mathbf{E}}}$ if and only if $T \in \Qsigma$.
    Hence $\mathcal{T}_{\widehat{\mathbf{E}}} = \Qsigma$, and so by~\eqref{JTF rewrite Sp} we obtain~\eqref{WTS JQF}, which completes the proof for $\O^*(2n)$.
    
    Similarly, let $G_\R = \U(p,q)$.
    By~\eqref{F hat E hat GL}, we can rewrite~\eqref{JTF} as
    \begin{equation}
        \label{JTF rewrite GL}
        \widehat{\mathcal{J}}_k(\sigma) = \Bigg(\coprod_{\mathbf{E} \in \widehat{\mathcal{E}}_k} \mathcal{T}_{\mathbf{E}}(\sigma) \Bigg) \times \widehat{\mathcal{F}}_k.
    \end{equation}
    By~\eqref{T_E}, $T \in \coprod_{\mathbf{E} \in \widehat{\mathcal{E}}_k} \mathcal{T}_{\mathbf{E}}(\sigma)$ if and only if ${\rm end}(T) \in \widehat{\mathcal{E}}_k$.
    In turn, by~\eqref{F hat E hat GL}, ${\rm end}(T) \in \widehat{\mathcal{E}}_k$ if and only if ${\rm end}(T) = \{(p, j_u)\}_{u=1}^{k^+} \sqcup \{(i_t, q)\}_{t=1}^{k^-}$ satisfies~\eqref{E max condition GL}, which we reproduce below:
    \[
        \#\{u : j_u \leq q - \ell \} + \#\{t : i_t \leq p - \ell\} = k - \ell \quad \text{for all $1 \leq \ell \leq \min\{k,r\}$}.
    \]
    By~\eqref{end T GL details}, each $T^+_{u,1} \geq j_u$ and $T^-_{t,1} \geq i_t$, and therefore, we have ${\rm end}(T) \in \widehat{\mathcal{E}}_k$ if and only if
    \[
    \#\{u : T^+_{u,1} \leq q - \ell \} + \#\{t : T^-_{t,1} \leq p-\ell \} \leq k - \ell \quad \text{for all $1 \leq \ell \leq \min\{k,r\}$},
    \]
    if and only if
    \begin{equation}
        \label{penultimate in GL proof}
        \#\{x \in T^+_1 : x \leq q - \ell \} + \#\{y \in T^-_1 : y \leq p - \ell\} \leq k - \ell \quad \text{for all $1 \leq \ell \leq \min\{k,r\}$};
    \end{equation}
    then by setting $i \coloneqq k - \ell  + 1$, we rewrite~\eqref{penultimate in GL proof} as 
    \begin{equation}
        \label{ultimate in GL proof}
        \#\{x \in T^+_1 : x < q - k + i\} + \#\{y \in T^-_1 : y < p - k + i \} < i \quad \text{for all $\max\{0, \: k-r\} < i \leq k$},
    \end{equation}
    which coincides with the necessary and sufficient condition for $T \in \Qsigma$ in Definition~\ref{def:T(sigma)}.
    Consequently, ${\rm end}(T) \in \widehat{\mathcal{E}}_k$ if and only if $T \in \Qsigma$.
    Hence $\coprod_{\mathbf{E} \in \widehat{\mathcal{E}}} \mathcal{T}_{\mathbf{E}}(\sigma) = \Qsigma$, and so by~\eqref{JTF rewrite GL} we obtain~\eqref{WTS JQF}, which completes the proof for $\U(p,q)$.
\end{proof}

\begin{theorem}
    \label{thm:main result}
    Assume one of the three settings in Table~\ref{table:main info};
    if $G_\R = \operatorname{Mp}(2n, \R)$, then assume that $k \notin [n+1, \: 2n-2]$.
    Let $L_{\la(\sigma)}$ be the simple $\g$-module with highest weight $\la(\sigma)$.
    We have
    \[
    \Deg L_{\la(\sigma)} = \#(\Qsigma \times \PP),
    \]
    with $\Qsigma$ and $\PP$ as in Definitions~\ref{def:T(sigma)} and~\ref{def:R}.
    In particular, if $k \geq r$, then $\Deg L_{\la(\sigma)} = \#\Qsigma$.
\end{theorem}

\begin{conj}
\label{conj:main}
    Theorem~\ref{thm:main result} remains true without the assumption on $k$ in the $\operatorname{Mp}(2n, \R)$ case.
\end{conj}

\begin{proof}[Proof of Theorem~\ref{thm:main result}]

    Let $G_\R = \U(p,q)$ or $\O^*(2n)$.
    First suppose $k < s$.
    By Lemma~\ref{lemma:Deg equals J hat} we have $\Deg L_{\la(\sigma)} = \# \widehat{\mathcal{J}}_k(\sigma)$, and so the result follows from the bijection ${\rm id} \times \Theta$ in Lemma~\ref{lemma:Theta tilde}.
    In the case $k \geq r$, we have $\# \PP = 1$ by~\eqref{P when k geq r}.

    Suppose $k \geq s$.
    The $k \geq s$ case of Lemma~\ref{lemma:Stanley decomp via jellyfish}, in~\eqref{jellyfish decomp k geq s}, states that $L_{\la(\sigma)} \cong \bigoplus_{T \in \Tsigma} S \, \varphi_T$.
    In this decomposition, the number of Stanley spaces equals $\#\Tsigma$, all of which have the same dimension (i.e., the Krull dimension of $S$, namely $\# \mathbf{P}$).
    Thus by Lemma~\ref{lemma:Bdeg Stanley decomp}, we have $\Deg L_{\la(\sigma)} = \# \mathcal{T}(\sigma)$.
    Since $k \geq s \geq r$, we have $\mathcal{T}(\sigma) = \Qsigma$ by~\eqref{Q = T}, and $\# \PP = 1$ by~\eqref{P when k geq r}.
    This completes the proof for $\U(p,q)$ and $\O^*(2n)$.

    Let $G_\R = \Mp(2n, \R)$.
    Recall from Table~\ref{table:main info} that $r=n$.
    Thus if $k \leq n$, then by the result~\eqref{NOT main result} from~\cite{NOT}, we have
    \[
    \Deg L_{\la(\sigma)} = \dim U_\sigma \cdot \deg \overline{\scrO}_k.
    \]
    But we have $\dim U_\sigma = \# \Qsigma$ by Lemma~\ref{lemma:U sigma F lambda}(a), and $\deg \overline{\scrO}_k = \#\PP$ by Proposition~\ref{prop:Dk general} below.
    This completes the proof for the case $k \leq n$.
    On the other hand, recall from Table~\ref{table:main info} that $s = 2n-1$, and now suppose that $k \geq 2n-1$.
    Then $L_{\la(\sigma)} \cong S \otimes F_{\la(\sigma)}$ is a free $S$-module, and it follows that $\Deg L_{\la(\sigma)} = \dim F_{\la(\sigma)}$.
    By Lemma~\ref{lemma:U sigma F lambda}(b), we have $\dim F_{\la(\sigma)} = \#\Qsigma$, and by~\eqref{P when k geq r}, we have $\#\PP = 1$.
    This completes the proof for the case $k \geq 2n-1$.   
\end{proof}

Note that explicit formulas for the Bernstein degree are obtained by combining Theorem~\ref{thm:main result} with our enumerative Propositions~\ref{prop:enumerate Qsigma} and~\ref{prop:enumerate R}.


\section{Open problems and observations}
\label{sec:open}

The main result in this paper leads naturally to the following problems.

\begin{problem}
    \label{problem:conjecture}
    Prove Conjecture~\ref{conj:main}.
\end{problem}

From our perspective, the ideal way to solve Problem~\ref{problem:conjecture} would be to extend Lemma~\ref{lemma:Stanley decomp via jellyfish} to the $G_\R = \operatorname{Mp}(2n, \R)$ case, by defining analogues of jellyfish for $\operatorname{Mp}(2n, \R)$.
Note that  Conjecture~\ref{conj:main} at least provides a candidate for the subset of maximal jellyfish, namely $\widehat{\mathcal{J}}_k(\sigma) = \Qsigma \times \Theta(\PP)$.
We hope that this may ultimately suggest the correct way to define the full set $\mathcal{J}_k(\sigma)$ of jellyfish, which (upon proof) would then lead to the solution of Problem~\ref{problem:conjecture}.

\begin{problem}
    \label{problem:O 2n outside dual pairs}
    Extend Theorem~\ref{thm:main result} to those unitary highest weight modules of $\O^*(2n)$ arising outside the dual pair setting.
\end{problem}

Unlike the groups $\U(p,q)$ and $\operatorname{Mp}(2n, \R)$, the group $\O^*(2n)$ admits unitary highest weight representations arising outside the dual pair setting; that is, certain unitary highest weight modules are not of the form $L_{\la(\sigma)}$ for some $\sigma \in \widehat{H}(k)$.
For details, see~\cite{DES}*{pp.~48--49}.
Without the group $H(k)$ in the picture, Problem~\ref{problem:O 2n outside dual pairs} amounts to determining a set of combinatorial objects, whether tableaux or otherwise, to play the role of $\Qsigma$ in Theorem~\ref{thm:main result}.

Solving Problems~\ref{problem:conjecture}--\ref{problem:O 2n outside dual pairs} would establish a uniform combinatorial interpretation of the Bernstein degree of \emph{all} unitary highest weight modules for the groups $\U(p,q)$, $\operatorname{Mp}(2n, \R)$, and $\O^*(2n)$.
The ultimate goal of the program is to extend such an interpretation beyond those groups:

\begin{problem}
    \label{problem:exceptionals}
    Extend Theorem~\ref{thm:main result} uniformly to arbitrary unitary highest weight modules (including those of the exceptional groups).
\end{problem}

Unlike Problem~\ref{problem:conjecture} (where for $\operatorname{Mp}(2n, \R)$ we had already defined $\Qsigma$ and $\PP$), and Problem~\ref{problem:O 2n outside dual pairs} (where for $\O^*(2n)$ we had already defined $\PP$), solving Problem~\ref{problem:exceptionals} requires us to define analogues of both $\Qsigma$ and $\PP$ in a uniform manner for Hermitian symmetric pairs of all types.
As a first step toward solving Problem~\ref{problem:exceptionals}, we now give a uniform generalization of the diagram $D_k$ that yields an analogous set $\PP$ for all groups of Hermitian type; see Proposition~\ref{prop:Dk general} below.

In the following discussion, suppose $G_\R$ is of Hermitian type, that is, one of the seven groups listed in Table~\ref{table:Dk} (page~\pageref{table:Dk}).
Recall from~\eqref{Phi p+} the set $\Phi(\p^+)$ of positive noncompact roots of the complexified Lie algebra $\g$.
This set $\Phi(\p^+)$ is a graded poset, where the grading is given by the height of the roots in $\Phi(\p^+)$.
In the third column of Table~\ref{table:Dk}, we show the Hasse diagram of $\Phi(\p^+)$ for the five types in which $\g$ is simply laced.
For the remaining two types, $\g$ is not simply laced:
\begin{itemize}
   \item $G_\R = \operatorname{Mp}(2n, \R)$, in which case $\g_\R = \sp(2n, \R)$ and $\g = \sp(2n, \C)$.
    \item $G_\R = \SO(2, 2n-1)$, in which case $\g_\R = \so(2, 2n-1)$ and $\g = \so(2n+1)$.
\end{itemize}
In each of these two types, however, there is a certain Hermitian symmetric pair $(\widetilde{\g}_\R, \widetilde{\k}_\R)$ with $\widetilde{\g}$ simply laced, admitting a natural embedding $\g_\R \hookrightarrow \widetilde{\g}_\R$, as shown below (taking $n=7$):

\begin{equation}
    \label{Phi tilde p diagrams}
    \begin{array}{ccccc}
     \g_\R & \Phi(\p^+) & & \widetilde{\g}_\R & \Phi( \, \widetilde{\p}^+) \\[2ex]
     \sp(2n, \R) & \begin{tikzpicture}[scale=.2, baseline=(current bounding box.center),rotate=45]
        \node at (0,7) {};
        \foreach \x in {1,...,7}{\foreach \y in {\x,...,7}{\node [dot] at (\x,\y) {};}}
        \foreach \x in {1,...,7}{\draw (\x,7) -- (\x,\x) (1,\x) -- (\x,\x);}
\end{tikzpicture} 
&
&
\mathfrak{u}(n,n)
&
\begin{tikzpicture}[scale=.2, baseline=(current bounding box.center),rotate=135]
     \draw [gray, densely dotted, thick] (0,0) -- (8,8);
        \foreach \x in {1,...,7}{\foreach \y in {1,...,7}{\node [dot] at (\x,\y) {};}}
        \draw (1,1) grid (7,7);
    \end{tikzpicture}
    \\
    & & & & \\
\so(2, 2n-1) &
\begin{tikzpicture}[scale=.2, baseline=(current bounding box.center),rotate=45]
        \node at (0,8) {};
        \foreach \x in {1,...,7}{\node [dot] at (\x,7) {};
        \node [dot] at (1,\x) {};}
        \draw (1,1) -- (1,7) -- (7,7);
\end{tikzpicture}
&
&
\so(2,2n)
&
\begin{tikzpicture}[scale=.2, baseline=(current bounding box.center),rotate=45]
    \draw [gray, densely dotted, thick] (0,8) -- (5,3);
        \foreach \x in {1,...,7}{\node [dot] at (\x,7) {};
        \node [dot] at (1,\x) {};}
        \node [dot] at (2,6) {};
        \draw (1,1) -- (1,7) -- (7,7) (1,6) -- (2,6) -- (2,7);
\end{tikzpicture}
\end{array} 
\end{equation}

\noindent (In terms of Cartan types, $\widetilde{{\rm C}}_n = {\rm A}_{2n-1}$ and $\widetilde{{\rm B}}_n = {\rm D}_{n+1}$, respectively.)
We observe from~\eqref{Phi tilde p diagrams} that $\Phi( \, \widetilde{\p}^+)$ is self-dual as a graded poset, with height $2n-1$.
The dotted line in~\eqref{Phi tilde p diagrams} indicates the ``middle rank,'' that is, those roots $\alpha \in \Phi( \, \widetilde{\p}^+)$ such that $\operatorname{height}(\alpha) = n$.
Define the lower order ideal 
\begin{equation}
    \label{Phi <n}
    \Phi( \, \widetilde{\p}^+)_{\leqslant} \coloneqq \{\alpha \in \Phi( \, \widetilde{\p}^+) : \operatorname{height}(\alpha) \leq n \},
\end{equation}
whose Hasse diagram is the bottom ``half'' of that of $\Phi( \, \widetilde{\p}^+)$; see the third column of Table~\ref{table:Dk}.

\begin{dfn}
    \label{def:Dk}
    Let $G_\R$ be of Hermitian type, with real rank $r$.
    Let
    \[
    D_0 \coloneqq \text{the $135^\circ$ clockwise rotation of the Hasse diagram of } \begin{cases}
        \Phi(\p^+) & \text{if $\g$ is simply laced},\\
        \Phi( \, \widetilde{\p}^+)_{\leqslant} & \text{if $\g$ is not simply laced},
    \end{cases}
    \]
    with each point in the Hasse diagram replaced by a box in $D_0$.
    The \emph{interior} of a diagram $D$ is the set of all boxes that have a box of $D$ directly to their southwest.
    For $1 \leq k \leq r$, recursively define
    \[
    D_k \coloneqq \text{interior of $D_{k-1}$}.
    \]
\end{dfn}

In Table~\ref{table:Dk} (page~\pageref{table:Dk}), we present the case-by-case details from Definition~\ref{def:Dk}.
Note that for $G_\R = \U(p,q)$ or $\operatorname{Mp}(2n, \R)$ or $\O^*(2n)$, we recover the diagrams $D_k$ given in Definition~\ref{def:R}.
Therefore, the following is a true generalization of Definition~\ref{def:R}:
\begin{dfn}
    In the setting of Definition~\ref{def:Dk}, we define
    \[
    \PP \coloneqq \{ \textup{plane partitions bounded by $k$ and contained in $D_k$} \}.
    \]
\end{dfn}

For $P \in \PP$, let $P_{ij}$ denote its entry in row $i$ (from the top) and column $j$ (from the left).
If no such entry exists, then set $P_{ij} = 0$.
Define the statistic
\begin{equation}
    \label{c(P)}
    c(P) \coloneqq \sum_{i,j} \big(P_{ij} - \max\{P_{i+1,j}, P_{i,j-1}\}\big),
\end{equation}
where the sum ranges over all pairs $(i,j)$ corresponding to a box in $D_k$.

\begin{prop}
    \label{prop:Dk general}
    Let $(G_\R, K_\R)$ be an irreducible Hermitian symmetric pair of noncompact type, and let $r$ denote the real rank of $G_\R$.
    For all $1 \leq k \leq r$, we have the Hilbert series
    \[
    H_{ \C[\overline{\scrO}_k]}(t) = \frac{\sum_{P \in \PP} t^{c(P)}}{(1-t)^{\dim \p^+ - |D_k|}},
    \]
    and therefore
    \[
    \# \PP = \deg \overline{\scrO}_k. 
    \]
\end{prop}

\begin{proof}
    Observe that the bijection $\Theta : \PP \longrightarrow \widehat{\mathcal{F}}_k$ (see Definition~\ref{def:Theta}) generalizes uniformly across all Hermitian symmetric types, where $\widehat{\mathcal{F}}_k$ now denotes the set of facets of the $k$th order complex of $\Phi(\p^+)$ or $\Phi( \, \widetilde{\p}^+)_{\leqslant}$.
    (From now on, to simplify our descriptions, we will consider the $135^\circ$ clockwise rotation of this Hasse diagram, meaning that it has the same shape and orientation as the diagram $D_0$.
    We write $(i,j)$ to denote the root in row $i$ and column $j$ with respect to this orientation.)
    Thus by~\eqref{F hat as paths}, if $P \in \PP$ then $\Theta(P)$ is a maximal family of nonintersecting lattice paths $\mathbf{N}_1, \ldots, \mathbf{N}_k$ in the Hasse diagram of $\Phi(\p^+)$ or $\Phi( \, \widetilde{\p}^+)_{\leqslant}$.
    For $\mathbf{F} = \mathbf{N}_1 \sqcup \cdots \sqcup \mathbf{N}_k \in \widehat{\mathcal{F}}_k$, define its \emph{corners} as follows:
    \begin{equation}
        \label{corners}
        {\rm cor}(\mathbf{F}) \coloneqq \coprod_{t=1}^k {\rm cor}(\mathbf{N}_t), \quad \text{where }   
        {\rm cor}(\mathbf{N_t}) \coloneqq \{ (i,j) \in \mathbf{N}_t : (i-1, j), (i, j+1) \in \mathbf{N}_t \}.
    \end{equation}
    In other words, a corner of $\mathbf{F}$ is a point at which one of its constituent lattice paths $\mathbf{N}_t$ turns from south to east.
    (This is a typical convention in the literature, although it differs from the definition of a corner in our paper~\cite{EricksonHunzikerMOC2024}.)
    Since all elements of $\widehat{\mathcal{F}}_k$ have the same cardinality, we may define
    \[
    d_k \coloneqq \# \mathbf{F} \text{ for any (equivalently, every) } \mathbf{F} \in \widehat{\mathcal{F}}_k,
    \]
    coinciding with the definition of $d_k$ in~\eqref{F hat E hat}.
    We observe from Definition~\ref{def:Dk} that
    \begin{equation}
        \label{dk fact}
        d_k = |D_0| - |D_k|.
    \end{equation}
    Moreover, we have
    \begin{equation}
        \label{p+ and D0}
        \dim \p^+ = \# \Phi(\p^+) = 
        \begin{cases}
        2n-1 & \text{if $G_\R = \SO(2, 2n-1)$},\\
        |D_0| & \text{otherwise.}
    \end{cases}
    \end{equation}
    
    First consider the dual pair setting (i.e., the first three rows of Table~\ref{table:Dk}).
    In this setting, the $K$-orbit closure $\overline{\scrO}_k$ is the determinantal variety
    \begin{equation}
    \label{O bar det var}
    \overline{\scrO}_k = \begin{cases}
        \{X \in {\rm M}_{p,q} : \operatorname{rk} X \leq k \} & \text{if $G_\R = \U(p,q)$},\\
        \{X \in {\rm SM}_{n} : \operatorname{rk} X \leq k \} & \text{if $G_\R = \operatorname{Mp}(2n, \R)$},\\
        \{X \in {\rm AM}_{n} : \operatorname{rk} X \leq 2k \} & \text{if $G_\R = \O^*(2n)$},\\
    \end{cases}
    \end{equation}
    with notation as in the last column of Table~\ref{table:GR}.
    Applying the well-known results from the 1990s~\cites{Sturmfels,Herzog,Conca94,GhorpadeKrattenthalerPfaffians} on Stanley decompositions and Hilbert series of these determinantal rings, we obtain
    \begin{equation}
        \label{Hilbert dual pair setting}
    H_{\C[\overline{\scrO}_k]}(t) = \frac{\sum_{\mathbf{F} \in \widehat{\mathcal{F}}_k} t^{\# {\rm cor}(\mathbf{F})}}{(1-t)^{d_k}}.
    \end{equation}
    Let $P \in \PP$, and put $\Theta(P) = \mathbf{N}_1 \sqcup \cdots \sqcup \mathbf{N}_k$ as in Definition~\ref{def:Theta}.
    Recall the subdiagrams $P_t$ in Definition~\ref{def:Theta}.
    A box $(i,j)$ in $P_t$ is said to be an \emph{outer corner} if $P_t$ contains neither $(i+1,j)$ nor $(i,j-1)$.
    By the construction of $\Theta$ in Definition~\ref{def:Theta}, for all $1 \leq t \leq k$ there is a bijective correspondence
    \[
    {\rm cor}(\mathbf{N}_t) \longleftrightarrow \{ \text{outer corners of $P_t$} \} ,
    \]
    where ${\rm cor}(\mathbf{N}_t)$ was defined in~\eqref{corners}.
    Therefore, we have
    \begin{align}
        \# {\rm cor}(\Theta(P)) = \sum_{t=1}^k \# {\rm cor}(\mathbf{N}_t) &= \sum_{t=1}^k \# \{ \text{outer corners of $P_t$} \} \nonumber \\
        &= \sum_{i,j} \#\{t : \text{$(i,j)$ is an outer corner of $P_t$}\} \nonumber \\
        &= \sum_{i,j} \big(P_{ij} - \max\{P_{i+1,j}, P_{i,j-1}\}\big) \nonumber \\
        & = c(P). \label{c(P) = cor}
    \end{align}
    Thus, starting with~\eqref{Hilbert dual pair setting} and applying the bijection $\Theta$ along with~\eqref{dk fact}--\eqref{p+ and D0} and~\eqref{c(P) = cor}, we obtain
    \begin{equation}
        \label{end of Hilbert series proof}
    H_{\C[\overline{\scrO}_k]}(t) = \frac{\sum_{\mathbf{F} \in \widehat{\mathcal{F}}_k} t^{\# {\rm cor}(\mathbf{F})}}{(1-t)^{d_k}} = \frac{\sum_{P \in \PP} t^{\# {\rm cor}(\Theta(P))}}{(1-t)^{|D_0| - |D_k|}} = \frac{\sum_{P \in \PP} t^{c(P)}}{(1-t)^{\dim \p^+ - |D_k|}},
    \end{equation}
    which completes the proof in the dual pair setting.

    Outside the dual pair setting, with the exception of the case $G_\R = \SO(2,2n-1)$, we carried out a case-by-case verification of~\eqref{Hilbert dual pair setting} in our previous paper~\cite{EricksonHunzikerMOC2024}*{\S7}.
    We stated this result in~\cite{EricksonHunzikerMOC2024} in terms of the \emph{$k$th Wallach representation} $L_{-kc\zeta}$, whose Hilbert series coincides with that of $\C[\overline{\scrO}_k]$.
    Hence for all but $G_\R = \SO(2, 2n-1)$, the result follows again from~\eqref{end of Hilbert series proof}.

    Finally, let $G_\R = \SO(2, 2n-1)$.
    Since $r=2$, we need only verify the case $k=1$.
    As above, the Hilbert series of $\C[ \overline{\scrO}_1]$ coincides with that of the first Wallach representation of $\g_\R$.    
    By~\cite{EnrightHunziker04}*{Thm.~27}, this Hilbert series is given by
    \[
    H_{\C[\overline{\scrO}_1]}(t) = \frac{1+t}{(1-t)^{2n-2}}.
    \]
    Since $D_1$ is a single box, we have $\mathcal{P}_1 = \{ \ytableaushort{0}, \ytableaushort{1} \}$, where $c(\ytableaushort{0}) = 0$ and $c(\ytableaushort{1}) = 1$, and we have $\dim \p^+ - |D_k| = 2n-2$ by~\eqref{p+ and D0}.
    This completes the proof.    
\end{proof}

\begin{remark}
    Proposition~\ref{prop:Dk general} provides a combinatorial interpretation for the \emph{Wonderful Correspondence} observed in~\cite{EHP}*{\S8.1}.
    In particular, if $G_\R = \U(p,q)$ or $\Mp(2n, \R)$ or $\O^*(2n)$, then there is a ``reduced'' Hermitian symmetric pair $(\g'_\R, \k'_\R)$, such that the Hilbert series of $\C[\overline{\scrO}_k]$ (equivalently, of the $k$th Wallach representation $L_{-kc\zeta}$ of $\g_\R$) can be expressed in rational form as
    \[
    H_{\C[\overline{\scrO}_k]}(t) = \frac{H_{E}(t)}{(1-t)^{\dim \p^+ - \dim \p'^+}},
    \]
    where $E$ is the simple finite-dimensional $\g'$-module with highest weight $k \zeta'$.
    (Here $\zeta'$ is the unique fundamental weight of $\g'$ that is orthogonal to the compact roots; for details, see Table 3 on page 178 of~\cite{EHP}, and Table 3 in~\cite{EricksonHunzikerJCTA} for Types I--IIIa.)
    Thus, comparing with Proposition~\ref{prop:Dk general}, it follows that the Hilbert series of $E$ can be understood as the generating function for the $c$ statistic on plane partitions, namely
    \[
    H_E(t) = \sum_{P \in \PP} t^{c(P)},
    \]
    or equivalently, the generating function for the number of turns in maximal families of $k$ nonintersecting lattice paths in $\mathbf{P}$.
    Moreover, we have $|D_k| = \dim \p'^+$.
    The Wonderful Correspondence is a special case of \emph{Enright--Shelton reduction}~\cites{ES87,ES89}.
    In fact, for all seven types in Table~\ref{table:Dk}, the subdiagram $D_k$ can be interpreted as the Hasse diagram of $\Phi(\p'^+)$ for a certain pair $(\g'_\R, \k'_\R)$ obtained via Enright--Shelton reduction.
\end{remark}

\subsection*{Beyond the Wallach representations}

In Table~\ref{table: NOT exceptional}, we list certain families of unitary highest weight modules for the exceptional groups $\operatorname{E}_{6(-14)}$ and $\operatorname{E}_{7(-25)}$, excluding the Wallach representations mentioned above.
For certain families among these modules $L_\la$, we have observed an analogue of the Nishiyama--Ochiai--Taniguchi result~\eqref{NOT main result}, which we recall below:
\begin{equation}
    \label{NOT again}
    \Deg L_\la = \dim U_\sigma \cdot \deg \overline{\scrO}_k,
\end{equation}
where 
$k$ is the positive integer such that $\overline{\scrO}_k = \mathcal{AV}(L_\la)$, and $U_\sigma$ is the irreducible representation of the group $H_{\R}(k)$ (corresponding to $G_\R$ via Howe duality) with highest weight $\sigma$.
By an analogue of this result for $G_\R = \operatorname{E}_{6(-14)}$ or $\operatorname{E}_{7(-25)}$, we mean a family of unitary highest weight modules $L_\la$ for $G_\R$, along with a corresponding group $H_\R(k)$, and a correspondence $\la \mapsto \sigma$ that satisfies~\eqref{NOT again}.
In this sense, the findings in Table~\ref{table: NOT exceptional} suggest an analogue of the dual reductive pairs $(G_\R, H_\R(k))$ arising from Howe duality.
The one-parameter family for $(\operatorname{E}_7, \operatorname{B}_4)$ was observed by Kato--Ochiai~\cite{KO}; see their Proposition 5.2 and Remark 5.5.
However, the family for $\operatorname{E}_6$ and the two-parameter family for $(\operatorname{E}_7, \operatorname{F}_4)$ seem to be new.

\bigskip

\begin{table}[h]
    \centering
    \resizebox{\linewidth}{!}
{\begin{tblr}{colspec={|c|c|c|c|c|c|c|},stretch=1.5}

\hline

$G_\R$ & $\la$ & $k$ & $\deg \overline{\scrO}_k$ & $H_\R(k)$ & $\sigma$ & $\dim U_\sigma$ \\

\hline[2pt]


\SetCell[r=2]{} $\operatorname{E}_6$ & \SetCell[r=2]{} {$[0,0,0,0,a,-2a\!-\!6]$, \\ for $a \geq 0$} & \SetCell[r=2]{} $2$ & \SetCell[r=2]{} $1$ & $\operatorname{B}_3$ & $[a,0,0]$ & \SetCell[r=2]{} $\begin{aligned} &(a+1)(a+2)(a+3) \\ &(a+4)(2a+5) \\ & \div 120 \end{aligned}$ \\ \cline{5-6}

& & & & $\operatorname{G}_2$ & $[a,0]$ & \\ \hline



\vspace{-13pt} $\operatorname{E}_7$ & {$[0,0,0,0,0,a,-2a\!-\!8]$, \\ for $a \geq 0$} & $2$ & $3$ & $\operatorname{B}_4$ & $[a,0,0,0]$ & $\begin{aligned} & (a+1)(a+2)(a+3) \\ & (a+4)(a+5) \\ & (a+6)(2a+7) \\ & \div 5040 \end{aligned}$ \\ \hline

$\operatorname{E}_7$ & {$[0,0,0,0,a,b,-3a\!-\!2b\!-\!12]$, \\ for $a,b \geq 0$} & $3$ & $1$ & $\operatorname{F}_4$ & $[0,0,a,b]$ & {$\begin{aligned} \quad &(a+1)(a+2)(a+3)^2 \\ &(a+4)(a+5)(b+1) \\ &(a+b+2)(a+b+3) \\ &(a+b+4)^2 (a+b+5) \\ &(a+b+6) (2a+b+5) \\ & (2a\!+\!b\!+\!6) (2a\!+\!b\!+\!7)^2 \\ & (2a\!+\!b\!+\!8) (2a\!+\!b\!+\!9) \\& (3a\!+\!b\!+\!10)(3a\!+\!2b\!+\!11) \\ & \div 12070840320000 \end{aligned}$} \\ \hline

\end{tblr}}
    \caption{Certain families of unitary highest weight modules for $\operatorname{E}_{6(-14)}$ and $\operatorname{E}_{7(-25)}$, outside the Wallach representations.
    In the cases where $H_\R(k)$ is given, we have $\Deg L_\la = \dim U_\sigma \cdot \deg \overline{\scrO}_k$.
    The highest weights $\la$ and $\sigma$ (for $\g$ and $H(k)$, respectively) are written in fundamental weight coordinates, following Bourbaki.
    The Bernstein degree was obtained using the Hilbert series given in~\cite{EH04exceptional}, in \S\S6--7.
    The formulas for $\dim U_\sigma$ were obtained using the Weyl dimension formula in Mathematica.}
    \label{table: NOT exceptional}
\end{table}

\begin{table}[t]
    \ytableausetup{boxsize=1ex}
\[
\begingroup
\renewcommand{\arraystretch}{1.5}
\begin{array}{|c|c|c|c|c|c|}
\hline
G_\R & r & \Phi(\p^+) \text{ or } \Phi(\,\widetilde{\p}^+)_{\leqslant} & D_0 & k & D_k \\ \Xhline{2pt}
\U(p,q) & \min\{p,q\} &
\begin{tikzpicture}[scale=.2, baseline=(current bounding box.center),rotate=135]
        \foreach \x in {1,...,10}{\foreach \y in {1,...,7}{\node [dot] at (\x,\y) {};}}
        \draw (1,1) grid (10,7);
        \node at (0,8) {};
        \node at (11,0) {};
    \end{tikzpicture}
&
\begin{array}{c}
q \\
p \;\;
\ydiagram{10,10,10,10,10,10,10} \;\;\phantom{p} \\
\phantom{.}
\end{array}

&
k \leq r
&
\begin{array}{c} 
\phantom{\times}(p-k) \phantom{\times} \\[-1.3ex]
\times (q-k) \phantom{\times} \\[1.5ex]
\ydiagram{7,7,7,7,0,0}
\end{array}
\\ \hline
\operatorname{Mp}(2n, \R) & n & 
\begin{tikzpicture}[scale=.2, baseline=(current bounding box.center),rotate=135]
        \foreach \x in {1,...,10}{\foreach \y in {\x,...,10}{\node [dot] at (\x,\y) {};}}
        \foreach \x in {1,...,10}{\draw (\x,10) -- (\x,\x) (1,\x) -- (\x,\x);}
        \node at (0,11) {};
        \node at (3,1) {};
\end{tikzpicture}
&
\begin{array}{c} 
n \\
\ydiagram{10,9,8,7,6,5,4,3,2,1,0,0}
\end{array}
&
k \leq r
&
\begin{array}{c} 
n-k \\
\ydiagram{7,6,5,4,3,2,1,0,0}
\end{array}
\\ \hline
\O^*(2n) & \lfloor n/2 \rfloor & 
\begin{tikzpicture}[scale=.2, baseline=(current bounding box.center),rotate=45]
        \foreach \x in {1,...,9}{\foreach \y in {\x,...,9}{\node [dot] at (\x,\y) {};}}
        \foreach \x in {1,...,9}{\draw (\x,9) -- (\x,\x) (1,\x) -- (\x,\x);}
        \node at (0,0) {};
        \node at (10,10) {};
\end{tikzpicture}
&
\begin{array}{c} 
n-1 \\
\ydiagram{9,1+8,2+7,3+6,4+5,5+4,6+3,7+2,8+1,0}
\end{array}
&
k \leq r
&
\begin{array}{c} 
n-2k-1 \\
\ydiagram{3,1+2,2+1,0,0}
\end{array}
\\ \hline
\SO(2, 2n-2) & 2 & \begin{tikzpicture}[scale=.2, baseline=(current bounding box.center),rotate=45]
        \foreach \x in {1,...,9}{\node [dot] at (\x,9) {};
        \node [dot] at (1,\x) {};}
        \node [dot] at (2,8) {};
        \draw (1,1) -- (1,9) -- (9,9) (1,8) -- (2,8) -- (2,9);
        \node at (0,0) {};
        \node at (10,10) {};
\end{tikzpicture} 
&
\begin{array}{c} 
n-1 \\
\ydiagram{9,7+2,8+1,8+1,8+1,8+1,8+1,8+1,8+1}
\end{array}
&
1
&
\ydiagram{1} \\ \hline
\SO(2,2n-1) & 2 &
\begin{tikzpicture}[scale=.2, baseline=(current bounding box.center),rotate=45]
        \foreach \x in {1,...,10}{      \node [dot] at (1,\x) {};}
        \node [dot] at (2,9) {};
        \draw (1,1) -- (1,10) (1,9) -- (2,9);
        \node at (0,0) {};
        \node at (3,11) {};
\end{tikzpicture} 
 &
 \begin{array}{c} 
 n \\[-1ex]
\ydiagram{10,8+1}
\end{array}
 &
 1
 & \ydiagram{1} \\ \hline
{\rm E}_{6(-14)} & 2 &
\begin{tikzpicture}[scale=.2,baseline=(current bounding box.center),rotate=135]
\node at (1,1) [dot] {};
\node at (2,1) [dot] {};
\node at (3,1) [dot] {};
\node at (4,1) [dot] {};
\node at (5,1) [dot] {};
\node at (3,0) [dot] {};
\node at (4,0) [dot] {};
\node at (5,0) [dot] {};
\node at (4,-1) [dot] {};
\node at (5,-1) [dot] {};
\node at (6,-1) [dot] {};
\node at (4,-2) [dot] {};
\node at (5,-2) [dot] {};
\node at (6,-2) [dot] {};
\node at (7,-2) [dot] {};
\node at (8,-2) [dot] {};

\draw (1,1) -- (5,1) (3,0) -- (5,0) (4,-1) -- (6,-1) (4,-2) -- (8,-2) (3,1) -- (3,0) (4,1) -- (4,-2) (5,1) -- (5,-2) (6,-1) -- (6,-2);

\node at (0,1) {};
\node at (9,-2) {};
    \end{tikzpicture}
&
\ydiagram{0,5,2+3,3+3,3+5}
&
1
&
\ydiagram{2,1+1,1+2} \\ \hline
\multirow{2}{*}{
${\rm E}_{7(-25)}$
} 
&
\multirow{2}{*}{3} 
&
\multirow{2}{*}{
\begin{tikzpicture}[scale=.2,baseline=(current bounding box.center), rotate=135]
\node at (0,1) [dot] {};
\node at (1,1) [dot] {};
\node at (2,1) [dot] {};
\node at (3,1) [dot] {};
\node at (4,1) [dot] {};
\node at (5,1) [dot] {};
\node at (3,0) [dot] {};
\node at (4,0) [dot] {};
\node at (5,0) [dot] {};
\node at (4,-1) [dot] {};
\node at (5,-1) [dot] {};
\node at (6,-1) [dot] {};
\node at (4,-2) [dot] {};
\node at (5,-2) [dot] {};
\node at (6,-2) [dot] {};
\node at (6,-3) [dot] {};
\node at (4,-3) [dot] {};
\node at (5,-3) [dot] {};
\node at (7,-2) [dot] {};
\node at (8,-2) [dot] {};
\node at (7,-3) [dot] {};
\node at (8,-3) [dot] {};
\node at (7,-4) [dot] {};
\node at (8,-4) [dot] {};
\node at (8,-5) [dot] {};
\node at (8,-6) [dot] {};
\node at (8,-7) [dot] {};

\draw (0,1) -- (5,1) (3,0) -- (5,0) (4,-1) -- (6,-1) (4,-2) -- (8,-2) (4,-3) -- (8,-3) (7,-4) -- (8,-4) (8,-2) -- (8,-7) (7,-2) -- (7,-4) (6,-1) -- (6,-3) (5,1) -- (5,-3) (4,1) -- (4,-3) (3,1) -- (3,0);

\end{tikzpicture}
}
&
\multirow{2}{*}{
\ydiagram{0,0,6,3+3,4+3,4+5,4+5,7+2,8+1,8+1,8+1}
}
& 
1
&
\ydiagram{0,0,2,1+1,1+2,1+4,4+1,0,0} \\ \cline{5-6}
& & & & 2
& \ydiagram{0,0,0,0,1,0,0,0,0} \\ \hline
\end{array}
\endgroup
\]
    \caption{The diagram $D_k$ in Definition~\ref{def:Dk}.
    In the two non-simply laced cases, namely $\operatorname{Mp}(2n, \R)$ and $\SO(2, 2n-1)$, the third column shows the Hasse diagram of $\Phi(\,\widetilde{\p}^+)_{\leqslant}$ defined in~\eqref{Phi <n}. 
    For our diagrams we take $n=10$, with $k=3$ in the first three rows.
    Note that $D_r = \varnothing$ in all types.}
    \label{table:Dk}
\end{table}


\clearpage

\bibliographystyle{amsplain}
\bibliography{references}

@article{Andrews,
 author = {Andrews, G.},
 title = {Mac{Mahon}'s conjecture on symmetric plane partitions},
 fjournal = {Proceedings of the National Academy of Sciences of the United States of America},
 journal = {Proc. Natl. Acad. Sci. USA},
 volume = {74},
 pages = {426--429},
 year = {1977},
 doi = {10.1073/pnas.74.2.426},
 url = {europepmc.org/articles/pmc392301},
}

@article {BaiHunziker,
    AUTHOR = {Bai, Z. and Hunziker, M.},
     TITLE = {The {G}elfand--{K}irillov dimension of a unitary highest weight
              module},
   JOURNAL = {Sci. China Math.},
  FJOURNAL = {Science China. Mathematics},
    VOLUME = {58},
      YEAR = {2015},
    NUMBER = {12},
     PAGES = {2489--2498},
      ISSN = {1674-7283,1869-1862},
   MRCLASS = {22E47 (17B10 22E46)},
  MRNUMBER = {3429262},
MRREVIEWER = {William\ M.\ McGovern},
       DOI = {10.1007/s11425-014-4968-y},
       URL = {https://doi.org/10.1007/s11425-014-4968-y},
}

@book{BrunsHerzog,
 author = {Bruns, W. and Herzog, J.},
 title = {Cohen-{M}acaulay rings.},
 edition = {Rev. ed.},
 fseries = {Cambridge Studies in Advanced Mathematics},
 series = {Camb. Stud. Adv. Math.},
 volume = {39},
 year = {1998},
 publisher = {Cambridge: Cambridge University Press}
}

@article {BKU,
    AUTHOR = {Bruns, W. and Krattenthaler, C. and Uliczka, J.},
     TITLE = {Stanley decompositions and {H}ilbert depth in the {K}oszul
              complex},
   JOURNAL = {J. Commut. Algebra},
  FJOURNAL = {Journal of Commutative Algebra},
    VOLUME = {2},
      YEAR = {2010},
    NUMBER = {3},
     PAGES = {327--357},
      ISSN = {1939-0807,1939-2346},
   MRCLASS = {13F20 (13D02)},
  MRNUMBER = {2728147},
MRREVIEWER = {Adam\ L.\ Van Tuyl},
       DOI = {10.1216/JCA-2010-2-3-327},
       URL = {https://doi.org/10.1216/JCA-2010-2-3-327},
}

@article {Conca94,
    AUTHOR = {Conca,
    A.},
    TITLE = {Gr\"{o}bner bases of ideals of minors of a symmetric matrix},
    JOURNAL = {J. Algebra},
    FJOURNAL = {Journal of Algebra},
    VOLUME = {166},
    YEAR = {1994},
    NUMBER = {2},
    PAGES = {406--421},
    ISSN = {0021-8693},
    MRCLASS = {13C40 (13P10)},
    MRNUMBER = {1279266},
    MRREVIEWER = {Rafael H. Villarreal},
    DOI = {10.1006/jabr.1994.1160},
    URL = {https://doi.org/10.1006/jabr.1994.1160},
    }

@article {DES,
    AUTHOR = {Davidson, M. and Enright, T. and Stanke, R.},
     TITLE = {Differential operators and highest weight representations},
   JOURNAL = {Mem. Amer. Math. Soc.},
  FJOURNAL = {Memoirs of the American Mathematical Society},
    VOLUME = {94},
      YEAR = {1991},
    NUMBER = {455},
     PAGES = {iv+102},
      ISSN = {0065-9266,1947-6221},
   MRCLASS = {22E47},
  MRNUMBER = {1081660},
MRREVIEWER = {D.\ Mili\v{c}i\'{c}},
       DOI = {10.1090/memo/0455},
       URL = {https://doi.org/10.1090/memo/0455},
}

@article {EnrightHunziker04,
    AUTHOR = {Enright, T. and Hunziker, M.},
     TITLE = {Resolutions and {H}ilbert series of determinantal varieties
              and unitary highest weight modules},
   JOURNAL = {J. Algebra},
  FJOURNAL = {Journal of Algebra},
    VOLUME = {273},
      YEAR = {2004},
    NUMBER = {2},
     PAGES = {608--639},
      ISSN = {0021-8693},
   MRCLASS = {17B10 (22E46)},
  MRNUMBER = {2037715},
       DOI = {10.1016/S0021-8693(03)00159-5},
       URL = {https://doi.org/10.1016/S0021-8693(03)00159-5},
}

@article{EH04exceptional,
 author = {Enright, T. and Hunziker, M.},
 title = {Resolutions and {Hilbert} series of the unitary highest weight modules of the exceptional groups},
 journal = {Represent. Theory},
 issn = {1088-4165},
 volume = {8},
 pages = {15--51},
 year = {2004}
}

@incollection {EHP,
    AUTHOR = {Enright, T. and Hunziker, M. and Pruett, W.},
     TITLE = {Diagrams of {H}ermitian type, highest weight modules, and
              syzygies of determinantal varieties},
 BOOKTITLE = {Symmetry: {R}epresentation {T}heory and its {A}pplications},
    SERIES = {Progr. Math.},
    VOLUME = {257},
     PAGES = {121--184},
 PUBLISHER = {Birkh\"{a}user/Springer, New York},
      YEAR = {2014},
   MRCLASS = {17B22 (05E15 13D02 17B10)},
  MRNUMBER = {3363009},
MRREVIEWER = {William M. McGovern},
       DOI = {10.1007/978-1-4939-1590-3\_6},
       URL = {https://doi.org/10.1007/978-1-4939-1590-3_6},
}

@article {ES87,
    AUTHOR = {Enright, T. and Shelton, B.},
     TITLE = {Categories of highest weight modules: applications to classical {H}ermitian symmetric pairs},
   JOURNAL = {Mem. Amer. Math. Soc.},
  FJOURNAL = {Memoirs of the American Mathematical Society},
    VOLUME = {67},
      YEAR = {1987},
    NUMBER = {367},
     PAGES = {iv+94},
      ISSN = {0065-9266},
   MRCLASS = {22E47 (17B10)},
  MRNUMBER = {888703},
MRREVIEWER = {Floyd L. Williams},
       DOI = {10.1090/memo/0367},
       URL = {https://doi.org/10.1090/memo/0367},
}

@article {ES89,
    AUTHOR = {Enright, Thomas J. and Shelton, Brad},
     TITLE = {Highest weight modules for {H}ermitian symmetric pairs of
              exceptional type},
   JOURNAL = {Proc. Amer. Math. Soc.},
  FJOURNAL = {Proceedings of the American Mathematical Society},
    VOLUME = {106},
      YEAR = {1989},
    NUMBER = {3},
     PAGES = {807--819},
      ISSN = {0002-9939},
   MRCLASS = {17B10 (22E46 22E47)},
  MRNUMBER = {961404},
MRREVIEWER = {Brian D. Boe},
       DOI = {10.2307/2047440},
       URL = {https://doi.org/10.2307/2047440},
}

@article {EW,
    AUTHOR = {Enright, T. and Willenbring, J.},
     TITLE = {Hilbert series, {H}owe duality and branching for classical
              groups},
   JOURNAL = {Ann. of Math. (2)},
  FJOURNAL = {Annals of Mathematics. Second Series},
    VOLUME = {159},
      YEAR = {2004},
    NUMBER = {1},
     PAGES = {337--375},
      ISSN = {0003-486X},
   MRCLASS = {22E47 (17B10)},
  MRNUMBER = {2052357},
MRREVIEWER = {Mark R. Sepanski},
       DOI = {10.4007/annals.2004.159.337},
       URL = {https://doi.org/10.4007/annals.2004.159.337},
}

@article {EricksonHunzikerMOC2024,
    AUTHOR = {W. Erickson and M. Hunziker},
     TITLE = {Stanley decompositions of modules of covariants},
     journal = {Algebr.~Combin.},
     NOTE = {To appear. \url{	arXiv:2312.16749}}
}

@article{EricksonHunzikerJCTA,
 author = {Erickson, W. and Hunziker, M.},
 title = {Dimension identities, almost self-conjugate partitions, and {BGG} complexes for {Hermitian} symmetric pairs},
 journal = {J. Combin. Theory Ser. A},
 issn = {0097-3165},
 volume = {219},
 pages = {106118},
 year = {2026},
doi = {10.1016/j.jcta.2025.106118}
}

@article{GV85,
 author = {Gessel, I. and Viennot, G.},
 title = {Binomial determinants, paths, and hook length formulae},
 fjournal = {Advances in Mathematics},
 journal = {Adv. Math.},
  volume = {58},
 pages = {300--321},
 year = {1985},
  doi = {10.1016/0001-8708(85)90121-5},
}

@incollection {GhorpadeKrattenthalerPfaffians,
    AUTHOR = {Ghorpade,
    S. and Krattenthaler,
    C.},
    TITLE = {The {H}ilbert series of {P}faffian rings},
    BOOKTITLE = {Algebra,
    arithmetic and geometry with applications ({W}est
    {L}afayette,
    {IN},
    2000)},
    PAGES = {337--356},
    PUBLISHER = {Springer,
    Berlin},
    YEAR = {2004},
    ISBN = {3-540-00475-0},
    MRCLASS = {13C40 (13D40)},
    MRNUMBER = {2037100},
    MRREVIEWER = {Ngo Viet Trung},
    }

@article{Giambelli,
 author = {Giambelli, G.},
 title = {Sulle variet{\`a} rappresentate coll' annullare determinanti minori contenuti in un determinante simmetrico od emisimmetrico generico di forme.},
 fjournal = {Atti della Accademia delle Scienze di Torino},
 journal = {Torino Atti},
 volume = {41},
 pages = {102--125},
 year = {1906},
}

@book {GW,
    AUTHOR = {Goodman, R. and Wallach, N.},
     TITLE = {Symmetry, {R}epresentations, and {I}nvariants},
    SERIES = {Graduate Texts in Mathematics},
    VOLUME = {255},
 PUBLISHER = {Springer, Dordrecht},
      YEAR = {2009},
     PAGES = {xx+716},
      ISBN = {978-0-387-79851-6},
   MRCLASS = {20G05 (14L35 17B10 20C30 20G20 22E46)},
  MRNUMBER = {2522486},
MRREVIEWER = {Vladimir V. Shchigolev},
       DOI = {10.1007/978-0-387-79852-3},
       URL = {https://doi.org/10.1007/978-0-387-79852-3},
}

@article{HC55,
author = {Harish-Chandra},
title = {Representations of semisimple {L}ie groups. {IV}},
journal = {Amer. J. Math.},
volume = {77},
year = {1955},
pages = {743-777}
}

@article{HC56,
author = {Harish-Chandra},
title = {Representations of semisimple {L}ie groups. {V}},
journal = {Amer. J. Math.},
volume = {78},
year = {1956},
pages = {1-41}
}

@article{HarrisTu,
 author = {Harris, J. and Tu, L.},
 title = {On symmetric and skew-symmetric determinantal varieties},
 fjournal = {Topology},
 journal = {Topology},
 issn = {0040-9383},
 volume = {23},
 pages = {71--84},
 year = {1984},
 language = {English},
 doi = {10.1016/0040-9383(84)90026-0},
 url = {hdl.handle.net/2027.42/24955}
}

@article {Herzog,
    AUTHOR = {Bruns,
    W. and Herzog,
    J.},
    TITLE = {On the computation of {$a$}-invariants},
    JOURNAL = {Manuscripta Math.},
    FJOURNAL = {Manuscripta Mathematica},
    VOLUME = {77},
    YEAR = {1992},
    NUMBER = {2-3},
    PAGES = {201--213},
    ISSN = {0025-2611},
    MRCLASS = {13H10 (13F50)},
    MRNUMBER = {1188581},
    MRREVIEWER = {Gert Naud\'{e}},
    DOI = {10.1007/BF02567054},
    URL = {https://doi.org/10.1007/BF02567054},
    }

@article {Howe89,
    AUTHOR = {Howe, R.},
     TITLE = {Remarks on classical invariant theory},
   JOURNAL = {Trans. Amer. Math. Soc.},
  FJOURNAL = {Transactions of the American Mathematical Society},
    VOLUME = {313},
      YEAR = {1989},
    NUMBER = {2},
     PAGES = {539--570},
      ISSN = {0002-9947},
   MRCLASS = {22E45 (11E57 15A72 20G05 22E47)},
  MRNUMBER = {986027},
MRREVIEWER = {Masato Wakayama},
       DOI = {10.2307/2001418},
       URL = {https://doi.org/10.2307/2001418},
}

@article{HVZ,
    title = {How to compute the {S}tanley depth of a monomial ideal},
    journal = {J. Algebra},
    volume = {322},
    number = {9},
    pages = {3151-3169},
    year = {2009},
    doi = {https://doi.org/10.1016/j.jalgebra.2008.01.006},
    url = {https://www.sciencedirect.com/science/article/pii/S0021869308000227},
    author = {J. Herzog and M. Vladoiu and X. Zheng},
}

@phdthesis{Jackson,
    author = {S. Jackson},
    title = {Standard monomial theory for reductive dual pairs},
    YEAR = {2003},
    INSTITUTION = {Yale University},
    MRNUMBER = {2704428},
    URL = {http://gateway.proquest.com/openurl?url_ver = Z39.88-2004&rft_val_fmt = info:ofi/fmt:kev:mtx:dissertation&res_dat = xri:pqdiss&rft_dat = xri:pqdiss:3084313},
    }

@article{Joseph,
 author = {Joseph, A.},
 title = {Annihilators and associated varieties of unitary highest weight modules},
 journal = {Ann. Sci. {\'E}c. Norm. Sup{\'e}r. (4)},
 volume = {25},
 number = {1},
 pages = {1--45},
 year = {1992},
 doi = {10.24033/asens.1642},
 url = {https://eudml.org/doc/82311}
}

@article {Jahan,
    AUTHOR = {Jahan, A.},
     TITLE = {Prime filtrations of monomial ideals and polarizations},
   JOURNAL = {J. Algebra},
  FJOURNAL = {Journal of Algebra},
    VOLUME = {312},
      YEAR = {2007},
    NUMBER = {2},
     PAGES = {1011--1032},
      ISSN = {0021-8693,1090-266X},
   MRCLASS = {13F20 (13F55)},
  MRNUMBER = {2333198},
MRREVIEWER = {Siamak\ Yassemi},
       DOI = {10.1016/j.jalgebra.2006.11.002},
       URL = {https://doi.org/10.1016/j.jalgebra.2006.11.002},
}

@article {KV,
    AUTHOR = {Kashiwara, M. and Vergne, M.},
     TITLE = {On the {S}egal--{S}hale--{W}eil representations and harmonic
              polynomials},
   JOURNAL = {Invent. Math.},
  FJOURNAL = {Inventiones Mathematicae},
    VOLUME = {44},
      YEAR = {1978},
    NUMBER = {1},
     PAGES = {1--47},
      ISSN = {0020-9910},
   MRCLASS = {22E45},
  MRNUMBER = {463359},
MRREVIEWER = {A. U. Klimyk},
       DOI = {10.1007/BF01389900},
       URL = {https://doi.org/10.1007/BF01389900},
}

@incollection{KO,
 author = {Kato, S. and Ochiai, H.},
 title = {The degrees of orbits of the multiplicity-free actions},
 booktitle = {Nilpotent {O}rbits, {A}ssociated {C}ycles and {W}hittaker {M}odels for {H}ighest {W}eight {R}epresentations},
 isbn = {2-85629-101-5},
 pages = {139--158},
 year = {2001},
 publisher = {Paris: Soci{\'e}t{\'e} Math{\'e}matique de France}
}

@article {KrattenthalerMem95,
    AUTHOR = {Krattenthaler, C.},
     TITLE = {The major counting of nonintersecting lattice paths and generating functions for tableaux},
   JOURNAL = {Mem. Amer. Math. Soc.},
  FJOURNAL = {Memoirs of the American Mathematical Society},
    VOLUME = {115},
      YEAR = {1995},
    NUMBER = {552},
     PAGES = {vi+109},
      ISSN = {0065-9266},
   MRCLASS = {05E10 (05A15)},
  MRNUMBER = {1254150},
MRREVIEWER = {Kevin W. J. Kadell},
       DOI = {10.1090/memo/0552},
       URL = {https://doi.org/10.1090/memo/0552},
}

@incollection{Krattenthaler,
 author = {Krattenthaler, C.},
 title = {Lattice path enumeration},
 booktitle = {Handbook of {E}numerative {C}ombinatorics},
 pages = {589--678},
 year = {2015},
 publisher = {Boca Raton, FL: CRC Press}
}

@misc{MacMahon,
 author = {MacMahon, P.},
 title = {Combinatory {A}nalysis. {Vol}. {I}, {II}},
 year = {1960},
  note = {New {York}: {Chelsea} {Publishing} {Company}. xix, 340 p.},
}

@article{NOT,
    AUTHOR = {Nishiyama, K. and Ochiai, H. and Taniguchi, K.},
     TITLE = {Bernstein degree and associated cycles of {H}arish-{C}handra
              modules---{H}ermitian symmetric case},
   JOURNAL = {Ast\'{e}risque},
  FJOURNAL = {Ast\'{e}risque},
    NUMBER = {273},
      YEAR = {2001},
     PAGES = {13--80},
      ISSN = {0303-1179},
   MRCLASS = {22E46 (14L30 32M15)},
  MRNUMBER = {1845714},
}

@misc{Porteous,
 author = {Porteous, I.},
 title = {Simple singularities of maps},
 year = {1971},
 note = {Proc. {Liverpool} {Singularities}-{Sympos}. {I}, {Dept}. {Pure} {Math}. {Univ}. {Liverpool} 1969-1970, {Lect}. {Notes} {Math}. 192, 286-307 (1971).},
 doi = {10.1007/bfb0066829},
}

@article {ProctorRSK,
    AUTHOR = {Proctor, R.},
     TITLE = {A {S}chensted algorithm which models tensor representations of
              the orthogonal group},
   JOURNAL = {Canad. J. Math.},
  FJOURNAL = {Canadian Journal of Mathematics. Journal Canadien de
              Math\'{e}matiques},
    VOLUME = {42},
      YEAR = {1990},
    NUMBER = {1},
     PAGES = {28--49},
  MRNUMBER = {1043509},
       DOI = {10.4153/CJM-1990-002-1},
       URL = {https://doi.org/10.4153/CJM-1990-002-1},
}

@article{Proctor94,
 author = {Proctor, R.},
 title = {Young tableaux, {Gelfand} patterns, and branching rules for classical groups},
 fjournal = {Journal of Algebra},
 journal = {J. Algebra},
 issn = {0021-8693},
 volume = {164},
 number = {2},
 pages = {299--360},
 year = {1994},
 doi = {10.1006/jabr.1994.1064}
}

@article{ProctorPPs,
title = {New symmetric plane partition identities from invariant theory work of {D}e {C}oncini and {P}rocesi},
journal = {European J. Combin.},
volume = {11},
number = {3},
pages = {289-300},
year = {1990},
doi = {https://doi.org/10.1016/S0195-6698(13)80128-X},
url = {https://www.sciencedirect.com/science/article/pii/S019566981380128X},
author = {R. Proctor}
}

@incollection {SchmidNotes,
    AUTHOR = {Schmid,
    W.},
    TITLE = {Geometric methods in representation theory},
    BOOKTITLE = {Poisson geometry,
    deformation quantisation and group
    representations},
    SERIES = {London Math. Soc. Lecture Note Ser.},
    VOLUME = {323},
    PAGES = {273--323},
    NOTE = {Lecture notes taken by Matvei Libine},
    PUBLISHER = {Cambridge Univ. Press,
    Cambridge},
    YEAR = {2005},
    ISBN = {978-0-521-61505-1; 0-521-61505-4},
    MRCLASS = {17B10 (22E46)},
    MRNUMBER = {2166454},
    }

@article {Stanley82,
    AUTHOR = {Stanley, R.},
     TITLE = {Linear {D}iophantine equations and local cohomology},
   JOURNAL = {Invent. Math.},
  FJOURNAL = {Inventiones Mathematicae},
    VOLUME = {68},
      YEAR = {1982},
    NUMBER = {2},
     PAGES = {175--193},
      ISSN = {0020-9910},
   MRCLASS = {10B05 (13D03 90C10)},
  MRNUMBER = {666158},
MRREVIEWER = {Hiroaki Terao},
       DOI = {10.1007/BF01394054},
       URL = {https://doi.org/10.1007/BF01394054},
}

@article {Stembridge,
    AUTHOR = {Stembridge, J.},
     TITLE = {Rational tableaux and the tensor algebra of {$\mathfrak{gl}_n$}},
   JOURNAL = {J. Combin. Theory Ser. A},
  FJOURNAL = {Journal of Combinatorial Theory. Series A},
    VOLUME = {46},
      YEAR = {1987},
    NUMBER = {1},
     PAGES = {79--120},
      ISSN = {0097-3165,1096-0899},
   MRCLASS = {05A15 (05A17 20C07 20C15)},
  MRNUMBER = {899903},
MRREVIEWER = {Bruce\ Sagan},
       DOI = {10.1016/0097-3165(87)90077-X},
       URL = {https://doi.org/10.1016/0097-3165(87)90077-X},
}

@article {Sturmfels,
    AUTHOR = {Sturmfels,
    B.},
    TITLE = {Gr\"{o}bner bases and {S}tanley decompositions of determinantal
    rings},
    JOURNAL = {Math. Z.},
    FJOURNAL = {Mathematische Zeitschrift},
    VOLUME = {205},
    YEAR = {1990},
    NUMBER = {1},
    PAGES = {137--144},
    ISSN = {0025-5874},
    MRCLASS = {14M12 (13P10 68Q40)},
    MRNUMBER = {1069489},
    MRREVIEWER = {Piotr Pragacz},
    DOI = {10.1007/BF02571229},
    URL = {https://doi.org/10.1007/BF02571229},
    }

@article {Vogan,
    AUTHOR = {Vogan, D.},
     TITLE = {Gel'fand-{K}irillov dimension for {H}arish-{C}handra
              modules},
   JOURNAL = {Invent. Math.},
  FJOURNAL = {Inventiones Mathematicae},
    VOLUME = {48},
      YEAR = {1978},
    NUMBER = {1},
     PAGES = {75--98},
      ISSN = {0020-9910,1432-1297},
   MRCLASS = {17B35},
  MRNUMBER = {506503},
MRREVIEWER = {Anthony\ Joseph},
       DOI = {10.1007/BF01390063},
       URL = {https://doi.org/10.1007/BF01390063},
}

@incollection{Vogan91,
 author = {Vogan, D.},
 title = {Associated varieties and unipotent representations},
 booktitle = {Harmonic analysis on reductive groups.},
 pages = {315--388},
 year = {1991},
 publisher = {Boston, MA etc.: Birkh{\"a}user}
}

@incollection{Vogan81,
 author = {Vogan, D.},
 title = {Singular unitary representations},
 year = {1981},
 language = {English},
pages={506-535},
booktitle={Noncommutative harmonic analysis and {Lie} groups},
city={{Marseille}- {Luminy}},
}

\end{document}